\documentclass[11pt, DIV10,a4paper]{article}

\usepackage{natbib}
\newcommand {\ctn}{\citet} 
\usepackage{graphicx,subfigure,amsmath,latexsym,amssymb}
\usepackage{float,epsfig,multirow,rotating,times}
\usepackage{upgreek,wrapfig}
\usepackage{comment}
\usepackage[title]{appendix}
\usepackage{algorithm}

\newcommand{\btheta}{\boldsymbol{\theta}}

\newcommand{\bbeta}{\boldsymbol{\beta}}

\newcommand{\bTheta}{\boldsymbol{\Theta}}

\newcommand{\bpsi}{\boldsymbol{\psi}}

\newcommand{\bt}{\boldsymbol{t}}

\newtheorem{theorem}{Theorem}

\newtheorem{corollary}[theorem]{Corollary}

\newtheorem{lemma}[theorem]{Lemma}

\newtheorem{remark}[theorem]{Remark}

\newenvironment{proof}[1][Proof]{\textbf{#1.} }{\ \rule{0.5em}{0.5em}}

\newcommand{\topline}{\hrule height 1pt width \textwidth \vspace*{2pt}}
\newcommand{\botline}{\vspace*{2pt}\hrule height 1pt width \textwidth \vspace*{4pt}}
\newtheorem{algo}{Algorithm} 

\numberwithin{equation}{section}
\numberwithin{algo}{section}
\numberwithin{table}{section}
\numberwithin{figure}{section}

\usepackage{setspace}
\usepackage[mathscr]{euscript}
\usepackage[margin=1in]{geometry}
\usepackage{xr-hyper}
\begin{document}
\normalsize

\title{\vspace{-0.8in}
{\bf Increasing Domain Infill Asymptotics for Stochastic Differential Equations Driven by Fractional Brownian Motion}}
\author{Trisha Maitra and Sourabh Bhattacharya\thanks{
	Trisha Maitra is an Assistant Professor in Prasanta Chandra Mahalanobis Mahavidyalaya, 111/3, B. T. Road, Kolkata 700108 and Sourabh Bhattacharya 
	is a Professor in
	Interdisciplinary Statistical Research Unit, Indian Statistical
	Institute, 203, B. T. Road, Kolkata 700108.
	Corresponding e-mail: sourabh@isical.ac.in.}}
\date{\vspace{-0.5in}}
\maketitle%

\begin{abstract}
Although statistical inference in stochastic differential equations ($SDE$s) driven by Wiener process has received significant attention in the literature, 
inference in those driven by fractional Brownian motion seem to have seen much less development in comparison, despite their importance in modeling long 
range dependence. In this article, we consider both classical and Bayesian inference in such fractional Brownian motion based $SDE$s, observed
on the time domain $[0,T]$. In particular, we consider
asymptotic inference for two parameters in this regard; a multiplicative parameter $\beta$ associated with the drift function, and the so-called ``Hurst parameter" 
$H$ of the fractional Brownian motion, when $T\rightarrow\infty$. For unknown $H$, the likelihood does not lend itself amenable to the 
popular Girsanov form, rendering usual asymptotic development difficult. As such, we develop increasing domain infill asymptotic theory, by discretizing
the $SDE$ into $n$ discrete time points in $[0,T]$, and letting $T\rightarrow\infty$, $n\rightarrow\infty$, such that either $n/T^2$ or $n/T$
tends to infinity. In this setup, we establish consistency and asymptotic normality of the maximum likelihood estimators, as well as consistency and
asymptotic normality of the Bayesian posterior distributions. However, classical or Bayesian asymptotic normality with respect to the Hurst parameter
could not be established. We supplement our theoretical investigations with simulation studies in a non-asymptotic setup, prescribing suitable methodologies
for classical and Bayesian analyses of $SDE$s driven by fractional Brownian motion. Applications to a real, close price data, along with comparison with
standard $SDE$ driven by Wiener process, is also considered. As expected, it turned out that our Bayesian fractional $SDE$ outperformed the other model and methods,
in both simulated and real data applications. 
\\[2mm]
\end{abstract}

\section{Introduction}
\label{sec:intro}

Stochastic differential equations ($SDE$s) are well-placed in both probability and statistics literature for their probabilistic importance as 
well as in statistical modeling of random phenomena varying continuously with time. The most important development in this regard have been for 
$SDE$s driven by Brownian motion (Wiener process). Accessible books on the theory of $SDE$s driven by Brownian motion can be found in \ctn{Arnold74}, 
\ctn{Ikeda81}, \ctn{Lipster01}, \ctn{Oksendal03}, \ctn{Mao11}, \ctn{Karatzas19} 
while statistical inference for such $SDE$s have been considered in \ctn{Basawa80}, \ctn{Iacus08}, \ctn{Bishwal08} and \ctn{Fuchs13}, the last containing interesting
applications from life sciences. Perhaps the most important applications of $SDE$s have been developed in the fields of 
biology and financial modeling; see, for example, \ctn{Allen10}, \ctn{Donnet13}, \ctn{Crepey13} and \ctn{Braumann19}. 

However, Wiener processes, that form the core of the $SDE$ literature, are inadequate for modeling phenomena that exhibit long-range dependence.
Such long range dependence is very much evident in asset returns and financial time series such as stock prices, foreign exchange rates, market indices and macroeconomics
(see, for example, \ctn{Henry03}). Geophysical time series (\ctn{Witt13}), internet traffic (\ctn{Park11}), etc. are also associated with long-range dependence. 
Theory and applications of long-range dependence can be found in \ctn{Doukhan03}. 
Note that long-range dependence is associated with scale-invariant self-similar processes with
stationary increments, as under suitable conditions, the latter can be shown to exhibit long-range dependence.

Continuous-time phenomena with long-range dependence can be more appropriately handled by $SDE$s driven by fractional Brownian motion, which is a self-similar
Gaussian process with stationary increments, consisting of a fractal index $H\in(0,1)$, usually termed as the ``Hurst parameter". Long-range dependence results
when $H>1/2$. Details of fractional Brownian motion and that of $SDE$s driven by fractional Brownian motion can be found in \ctn{Norros}, \ctn{Oks08}, \ctn{Mish29}, 
\ctn{Zien09}, \ctn{Xu18}. Statistical inference, including asymptotic theory in such $SDE$s, assuming that $H>1/2$ and known, is detailed in \ctn{BLSP10}; 
see also \ctn{Mish14}, \ctn{Neu14}, \ctn{Hu17}, where again $H$ is assumed to be known. For further related works with known $H$ and deterministic diffusion
coefficient, see \ctn{Dai21}, \ctn{Xiao19a}, \ctn{Xiao19b}, \ctn{Hu10}, \ctn{Tanaka20}, \ctn{Xiao11}.
In \ctn{Wang23}, a discretized version of the $SDE$ associated with fractional Ornstein-Uhlenbeck process (driven by fractional Brownian motion) is proposed,
where a two-stage estimation procedure is considered to estimate $H$ and the other parameters; closed form expressions for the estimators are available in this setup.
Asymptotic theory for the estimators are also developed.

As we shall subsequently discuss, although many attempts have been made to estimate $H$, most are at best heuristic, and in general there does not seem to exist any
theory for validation of such heuristic procedures. For instance, we are unaware of any asymptotic theory of such heuristic inference regarding $H$. 
Moreover, to our knowledge, in the general setup, unless $H$ is known, there does not exist any asymptotic theory of inference for the other parameters. 
Our goal in this article is to establish asymptotic theory of both classical and Bayesian inference for fractional Brownian motion based $SDE$s, 
considering $H$ and other relevant
parameters to be unknown. Among the other relevant parameters, here we focus on a single multiplicative real-valued parameter, $\beta$, associated with the
``drift function" of the underlying $SDE$. 

As we shall point out subsequently, for unknown $H$, the likelihood function does not admit expression
in terms of the convenient Girsanov formula that is usually employed for inference when $H$ is assumed to be known (see \ctn{Norros}, \ctn{Klept} and \ctn{BLSP10} for example).
In any case, in reality, continuous trajectories of the underlying stochastic process is not available, and the data are observed only at discrete time points. 
Hence, we construct a likelihood by discretizing the time domain $[0,T]~(T>0)$ of the fractional $SDE$ into $n$ time points, 
and quantifying the joint distribution associated with the $n$ time points. Our asymptotic setup is an ``increasing domain infill asymptotics" framework
where both $T\rightarrow\infty$ and $n\rightarrow\infty$ such that $n$ grows faster than $T$. In particular, we consider two cases, namely,
$n/T\rightarrow\infty$ and $n/T^2\rightarrow\infty$. Thus, the time domain is increased, but the number of discrete points in the time domain is increased faster
than $T$ or even $T^2$, attempting to fill up the domain first, so that asymptotic theory makes sense. This is the essence of increasing domain infill asymptotics 
that we employ in this article. As expected, when $n$ grows faster than even $T^2$, that is, when $n/T^2\rightarrow\infty$, our asymptotic results hold with fewer 
assumptions compared to the case $n/T\rightarrow\infty$. 

In our asymptotics framework, we establish strong consistency of the maximum likelihood estimators ($MLE$s) of $\beta$ and $H$ and asymptotic normality of the
$MLE$ of $\beta$, under mild assumptions. We also prove strong consistency of the posterior distribution of $(\beta,H)$ and asymptotic normality of the posterior of $\beta$,
in the same increasing domain infill asymptotics framework.  

For convenience, in the rest of our paper, we shall often denote ``$SDE$ driven by fractional Brownian motion" simply by ``fractional $SDE$". 
The roadmap of our paper is structured as follows. In Section \ref{sec:overview}, we provide an overview of the $SDE$ driven by fractional Brownian motion 
that we undertake in this work and in Section \ref{sec:discussion} we include a brief discussion of the main issues involved in 
the existing attempts of statistical inference in fractional $SDE$s. The modeled and the true likelihoods under discretization are presented in
Section \ref{sec:likelihoods}. We establish strong consistency of the $MLE$ of $(\beta,H)$ in Section \ref{sec:strong_consistency_MLE}, assuming compact parameter space.
Asymptotic normality of the $MLE$ of $\beta$ is established in Section \ref{sec:normality_beta}, when $H$ is estimated by its corresponding $MLE$.
In the same section it is shown that compact parameter space for $\beta$ is unnecessary for either strong consistency or asymptotic normality of the $MLE$ of $\beta$.
Inconsistency of the $MLE$ of $\beta$ when $H$ is estimated by any inconsistent estimator, is established in Section \ref{sec:incons_beta}.
In Section \ref{sec:bayesian_consistency}, we prove strong consistency of the posterior distribution of $(\beta,H)$, which is shown to hold 
even when the parameter space of $\beta$ is non-compact. In the same section we also show the types of neighborhoods the posterior distribution of $\beta$ asymptotically
concentrates on when $H$ is estimated by its $MLE$. 
Asymptotic normality of the posterior of $\beta$ when $H$ is estimated by its $MLE$, is established in Section \ref{sec:posterior_normaility}. 
In Section \ref{sec:simstudy} we illustrate classical and Bayesian inference in a simulation study
where the time domain is relatively small, a situation important for practical considerations, while applications of classical and Bayesian inference
for both fractional $SDE$ and standard $SDE$ driven by standard Brownian motion to a real, close price data is detailed in Section \ref{sec:realdata}, along with comparisons.
Finally, in Section \ref{sec:conclusion}, we summarize our contributions, make concluding remarks and provide directions for future research.
The proofs of our results, some associated technical lemmas and their proofs, as well as some technical illustrations and overview, are relegated to a supplement, whose sections, 
have the prefix ``S-" when referred to in this paper.

\section{Overview of our $SDE$ driven by fractional Brownian motion}
\label{sec:overview}
We assume that data has been generated from the following $SDE$ driven by fractional Brownian motion:
\begin{equation}
	d\tilde X_t=\beta B(t, \tilde X_t)dt+C (t, \tilde X_t)d\tilde W_t^H
\label{eq:fsde0}
\end{equation}
where the Hurst parameter $H\in (0,1)$, and $\beta\in\mathbb R$, the real line, and $\tilde X_0=\tilde x$ is the initial value of the stochastic process $\tilde X_t$, which 
is assumed to be continuously observed on the time interval $[0,T]$.
The function $B(t,\cdot)$ is a known, real-valued function on $\mathbb R$; this function is known as the drift function.
The function $C(t,\cdot):\mathbb R\mapsto\mathbb R$ is the known diffusion coefficient.

According to Theorem D.1.3 of \ctn{Oks08} we consider the following assumptions under which (\ref{eq:fsde0}) will have a unique solution.
\begin {itemize}
\item[(H1)] $B(t, u), C(t, u)$ be Lipschitz continuous in $u$, uniformly with respect to $t$, that is, for some $L>0$
$$|B(t, u_1) - B(t, u_2)| + |C(t, u_1) - C(t,u_2)|\leq L|u_1 - u_2|.$$
\item [(H2)] The initial values $x_0, y_0 \in \mathbb R$ and $\alpha_b>\alpha_1, \alpha_{\sigma}>2\alpha_2$ such that
$$B(\cdot, x_0)\in L^{\alpha_b}([0,T]), \quad\quad\quad\quad C(\cdot, y_0)\in L^{\alpha_{\sigma}}([0,T])$$
where $\alpha_i=(1-i|H-1/2|)^{-1}$ for $i=1,2$.
\end{itemize}
Let $(\Omega,\mathcal X, (\mathcal X_t ), P)$ be a stochastic basis satisfying the usual conditions. The natural filtration of a process is understood as the $P$-completion of the filtration generated by this process. In our considered $SDE$ (\ref{eq:fsde1}), $B =\{B(t,\tilde X_t), t \geq 0\}$ is an $(\mathcal X_t )$-adapted process and $C(t,\tilde X_t)$ is a 
non-vanishing, random function.

\subsection{True, data generating model and our inferential model}
We assume that data $\tilde x=\left\{\tilde x_t:t\in[0,T]\right\}$ is generated from (\ref{eq:fsde0}) with unknown true parameters $\btheta_0=(\beta_0,H_0)$.
We denote this true model, as well as the probability distribution associated with it by $P_{\btheta_0}$.

However, given data $\tilde x$, our inference regarding $\btheta=(\beta,H)$
will proceed by first creating a new fractional $SDE$ which is the same as (\ref{eq:fsde0}) but $B(t,\tilde X_t)$ and $C(t,\tilde X_t)$ replaced with 
the deterministic quantities $B(t,\tilde x_t)$ and $C(t,\tilde x_t)$, respectively. In other words,  
given data $\tilde x$, for inference we consider the following fractional $SDE$
\begin{equation}
	dX_t=\beta B(t, \tilde x_t)dt+C (t, \tilde x_t)dW_t^H.
\label{eq:fsde1}
\end{equation}
In the above, $\{W_t^H:t\in[0,T]\}$ is independent of $\{\tilde W_t^{\tilde H}:t\in[0,T]\}$, for any $H,\tilde H\in(0,1)$.

The fundamental difference between (\ref{eq:fsde0}) and (\ref{eq:fsde1}) is that the former has no semimartingale associated with it, making inferential procedures difficult, 
while, thanks to non-randomness of $B(t,\tilde x_t)$ and $C(t,\tilde x_t)$, the latter admits a Gaussian martingale (see Section \ref{sec:gaussian_martingale}), 
thereby significantly simplifying inference. Further, we argue below that given the observed data, both the equations lead to the
same likelihood and hence the same likelihood-based inference on $\btheta$ given the data. This vindicates the great generality of the formulation (\ref{eq:fsde1}).

It is important to mention that even if the data $\tilde x$ is assumed to be continuously observed, in spite of the Gaussian martingale, the likelihood based
on the so-called Girsanov formula is unavailable for (\ref{eq:fsde1}) when $H$ is unknown; see Section \ref{sec:discussion}. 
For (\ref{eq:fsde0}) such likelihood is not even conceivable, even for
known $H$, since no semimartingale is associated with it. Thus inference based on discretized data is inevitable. 
This is of course the most practical 
approach anyway, since
although from the theoretical perspective, the entire set $\tilde x$ is the data (solution of (\ref{eq:fsde0})), in practice, $\tilde x_t$ are observed 
only at discrete time points $t=t_i$,
for $i=1,\ldots,n$, for some positive integer $n$. 
Now let us consider inference by discretization of (\ref{eq:fsde0}) and (\ref{eq:fsde1}). With respect to (\ref{eq:fsde0}) and (\ref{eq:fsde1}), the discretizations are given,
respectively, by
\begin{equation}
\tilde X_{t_{i+1}}-\tilde X_{t_i}=\beta B(t_i, \tilde X_{t_i})(t_{i+1}-t_i)+C (t_i, \tilde X_{t_i})(\tilde W_{t_{i+1}}^H-\tilde W_{t_i}^H)
\label{eq:fsde0_discrete}
\end{equation}
and
\begin{equation}
X_{t_{i+1}}-X_{t_i}=\beta B(t_i, \tilde x_{t_i})(t_{i+1}-t_i)+C (t_i, \tilde x_{t_i})(W_{t_{i+1}}^H-W_{t_i}^H).
\label{eq:fsde1_discrete}
\end{equation}
Now, the distributions of $\{\tilde W_{t_{i+1}}^H-\tilde W_{t_i}^H;i=1,\ldots,n-1\}$ and $\{W_{t_{i+1}}^H- W_{t_i}^H;i=1,\ldots,n-1\}$ are the same,
albeit the increments within the respective sets are not independent. Thus, from (\ref{eq:fsde0_discrete}) and (\ref{eq:fsde1_discrete}), it follows that the distributions
of $\left\{\frac{\tilde X_{t_{i+1}}-\tilde X_{t_i}-\beta B(t_i, \tilde X_{t_i})(t_{i+1}-t_i)}{C (t_i, \tilde X_{t_i})};i=1,\ldots,n-1\right\}$ and
$\left\{\frac{X_{t_{i+1}}-X_{t_i}-\beta B(t_i, \tilde x_{t_i})(t_{i+1}-t_i)}{C (t_i, \tilde x_{t_i})};i=1,\ldots,n-1\right\}$ are the same.
Thus, the likelihoods corresponding to (\ref{eq:fsde0_discrete}) and (\ref{eq:fsde1_discrete}) have the same form. Now, for the likelihood associated with
(\ref{eq:fsde0_discrete}) we must substitute $\tilde X_{t_i}$ with the observed data $\tilde x_{t_i}$, for $i=1,\ldots,n$. Also, since we also fit (\ref{eq:fsde1_discrete})
to the observed data, we must substitute $X_{t_i}$ with the observed data $\tilde x_{t_i}$. Thus, the likelihoods associated with (\ref{eq:fsde0_discrete}) and
(\ref{eq:fsde1_discrete}) are the same. Hence, given the observed data, likelihood-based inference of $\btheta$ with respect to (\ref{eq:fsde0_discrete}) and
(\ref{eq:fsde1_discrete}) are the same.

The fundamental role of the Gaussian martingale associated with (\ref{eq:fsde1_discrete}) is to amply simplify the likelihood computation and the resultant asymptotic theory
by virtue of the corresponding independent increments.
The issue of discretization with respect to the Gaussian martingale for (\ref{eq:fsde1}) will be made precise in Section \ref{sec:likelihoods} 
as we formally introduce our premise. We establish classical and Bayesian consistency of $(\beta,H)$ and classical and Bayesian 
asymptotic normality of $\beta$ when the data are generated from
$P_{\btheta_0}$ but for inferential purposes (\ref{eq:fsde1}) is fitted to the data. Unfortunately, so far we have not been able to establish 
asymptotic normality associated with $H$.

In Section \ref{sec:illustration}, we illustrate our concept with the Ornstein-Uhlenbeck process, demonstrating that consistency can be achieved while
fitting (\ref{eq:fsde1}), even when the data has been generated from (\ref{eq:fsde0}).

\section{Discussion on existing attempts regarding inference in fractional $SDE$s}
\label{sec:discussion}
As already mentioned, consistency of the drift parameter has been studied by various authors assuming that 
the Hurst parameter $H$ is known. In relatively simple specific situations, such as in a linear model and Ornstein-Uhlenbeck process driven by fractional Brownian motion, 
\ctn{Breton98} and \ctn{Klept02} investigated statistical inference, and the latter proved consistency of the $MLE$ of the parameter,
assuming that $H$ is known. \ctn{Chiba20} studied asymptotic properties of an $M$-estimator for the drift parameter when the diffusion coefficient is a known constant
and $H\in(1/4,1/2)$ is known. \ctn{BLSP10} considered a more general setup, namely, the class of linear fractional $SDE$s, and investigated the asymptotic theory
of $MLE$ and Bayes estimators, again assuming known $H$. In particular, he established strong consistency and asymptotic normality of the $MLE$ and the 
Bayes estimator, under appropriate assumptions. His asymptotic theory hinges upon construction of an appropriate likelihood when $H$ is known. This likelihood,
which can be perceived as the Girsanov formula for the fractional Brownian motion case, has been derived by \ctn{Klept}. 
However, this relies upon the assumption that entire $x$ is observed, which is perhaps not quite realistic. 

Indeed, when $H$ is of known value, say $H_0$, then the distribution of the driving fractional Brownian motion $W^{H_0}_t$, is known. 
This then facilitates obtaining a Radon-Nikodym derivative of the measure of the original process associated with the $SDE$ with other model 
parameters $\bpsi$, with respect to that 
of an $SDE$ with known parameters, also driven by $W^{H_0}_t$, which is treated as the likelihood function. The form of this likelihood is given by
\begin{equation}
	L_T(\bpsi)=\exp\left\{-\int_0^TQ_{\bpsi}(s)dZ_s+\frac{1}{2}\int_0^TQ^2_{\bpsi}(s)dw^{H_0}_s\right\},
	\label{eq:girsanov}
\end{equation}
which has been exploited by the authors for inference regarding $\bpsi$, assuming $H$ to be known. However, a severe drawback of this Girsanov-based approach is
that the entire theory leading to (\ref{eq:girsanov}) breaks down when $H$ is unknown. Since in reality $H$ is always unknown, alternatives must be sought.

An alternative is to obtain some good estimator of $H$ and plug-in its value in (\ref{eq:girsanov}), pretending it to be $H_0$. If the estimator
is good enough in the sense of at least being consistent, then for large $T$, such a plug-in method need not be too unrealistic. When the data arise
from a fractional Brownian motion, the literature is concerned with only various types of estimator of $H$s, such as 
the $R/S$ estimator (\ctn{Man69}), Whittle estimator (\ctn{Whi53}), log-periodogram estimator (\ctn{Gew83}), to name a few. 
The Whittle estimator, for instance, is based on maximizing an approximate likelihood of the fractional Brownian motion and has asymptotic consistency and
asymptotic normality properties when the data arise from a fractional Brownian motion. However, when the data arise from some $SDE$ driven by fractional Brownian motion,
theoretical properties of such estimators, asymptotic or otherwise, need not be encouraging. It is thus not clear if a plug-in estimator of $H$ will work
even in empirical studies. In our experience, the aforementioned estimators work only when estimated from the data on fractional Brownian motion, but unfortunately,
in the $SDE$ context, such data are not available. Simulation studies reported in \ctn{Daria07} first estimate $H$ from simulated fractional Brownian motion, 
and then estimate the other parameters given such estimate of $H$. In realistic situations this is of course not possible. However, the above estimators,
particularly, the Whittle estimator, is expected to work well in fractional $SDE$s where the drift coefficient is close to zero and the diffusion coefficient is close to one,
almost surely. Moreover, since the Whittle estimator is computationally cheap, it is likely to have an computational advantage over the general, 
but computationally intensive methodologies that we propose in this paper.

\ctn{Kub12} study asymptotic properties of two estimators for fractional $SDE$ where $H\in (1/2,1)$ and where there are no other parameters, in a fixed-domain
asymptotics setup, where the time interval $[0,T]$ is partitioned uniformly, but $T$ remains fixed.
\ctn{Brou13} consider a special case of fractional $SDE$ leading to fractional Ornstein-Uhlenbeck process and investigate asymptotic properties of estimators
in the discretized case; however, even in such special scenario, their results for the drift coefficient are obtained only when $\frac{1}{2}<H<\frac{3}{4}$.
In special cases of fractional $SDE$s when the diffusion function is equal to one, 
and when the time interval is $[0,1]$ or $[0,\infty)$ with $H\leq 1/2$, and when there are no unknown parameters other than $H$, \ctn{Gairing20} provide two estimators for $H$,
proving their strong consistency and weak convergence to normality as the discretization becomes more and more fine. 
Proofs of their results also require separate restrictions for different ranges of $H$. Thus, the approaches of \ctn{Kub12} and 
\ctn{Gairing20} do not seem to be widely applicable
in practice. \ctn{Kub17} consider book-length treatment of parameter estimation under various specific $SDE$ models driven by fractional Brownian motion, 
given various ranges of $H$. Asymptotic theories, under suitable assumptions, are also provided. However, an unified estimation theory for general class of models 
with $H\in (0,1)$, does not seem to have been considered. 
Besides, the Bayesian counterpart has not been considered at all by these authors. 

Using a method based on
conic multivariate adaptive regression spline, numerical estimation of $H$ in fractional $SDE$s when there are no
other parameters and when the diffusion coefficient is a constant, is considered by \ctn{Yer14}, but no theoretical investigation has been undertaken.

Our goal in this paper is to develop both classical and Bayesian asymptotic theories for general fractional $SDE$s. We also develop methods that are 
applicable to practical, non-asymptotic scenarios, and demonstrate their importance with simulation experiments.
Because of the futility of (\ref{eq:girsanov}) in realistic situations, we do away entirely with the Girsanov formula, and consider discretization of the 
$SDE$ in question on the domain $[0,T]$, so that it consists of $n$ observations, and establish both classical and Bayesian increasing domain infill asymptotic theory  
for the model parameters as well as for $H$. Details follow.

\section{Modeled and true likelihoods under the increasing domain infill asymptotics framework}
\label{sec:likelihoods}

As already mentioned, we fit the data with the form (\ref{eq:fsde1}); $X$
is to be converted to $Z$ using (\ref{eq:zt}) given parameter $\btheta$ of (\ref{eq:fsde1}). 
In theory, this requires the continuous trajectory $\{X_s:s\in[0,t]\}$, for $t\in [0,T]$, but for practical computation we approximate the integral (\ref{eq:zt})
numerically using the available discrete trajectory. Below we provide details on the likelihood construction, assuming 
that data are available at discrete time points.

We consider the discretized version of (\ref{eq:fsde1}) on $[0,T]$ with step length $\frac{T}{n}$. Thus, denoting $Z_t$ by $Z^{\btheta}_t$ to indicate
its dependence on $\btheta$, using (\ref{eq:zt_decompose}) and Euler's discretization scheme, we obtain 
\begin{align}
	Z^{\btheta}_{t_{i+1}} - Z^{\btheta}_{t_i} &= \int_{t_i}^{t_{i+1}}Q_{\btheta}(s)dw^H_s +  M^H_{t_{i+1}} - M^H_{t_i}\notag\\
	&\approx Q_{\btheta}(t_i)(w_{t_{i+1}}^H-w_{t_i}^H) +  M^H_{t_{i+1}} - M^H_{t_i}
\label{eq:discretize1}
\end{align}
where $t_i=T\frac{i}{n}$ for $ i=1,\ldots, n$ and either $\frac{n}{T^2}\rightarrow \infty$ or $\frac{n}{T}\rightarrow \infty$ as both 
$T\rightarrow\infty$ and $n\rightarrow\infty$. 

It is important to clarify that by virtue of (\ref{eq:zt}), given $H$ and any realization $\tilde x$ (data) of $X$, $Z^{\btheta}_{t_i}$ is given by
\begin{equation}
	Z^{\btheta}_{t_i}=\int_0^{t_i}k_H(t_i,s)[C(s,\tilde x_s)]^{-1}d\tilde x_s,
\label{eq:zt_theta}
\end{equation}
which can be approximated by discretizing the above integral. Also, recalling from (\ref{eq:wth}) that 
$w_t^H=\lambda_H^{-1}t^{2-2H}$, given $\btheta$ and $\tilde x$, using chain rule of differentiation, $Q_{\btheta}(t_i)$ given by (\ref{eq:qh}) can be approximated numerically. 


We re-write (\ref{eq:discretize1}) as
\begin{equation}
	\Delta Z^{\btheta}_{t_i} - Q_{\btheta}(t_i)\Delta w_{t_i}^H =  \Delta  M^H_{t_i},
\label{eq:discretize2}
\end{equation}
replacing ``$\approx$" in (\ref{eq:discretize1}) with ``$=$".
As clarified above, the left hand side of (\ref{eq:discretize2}) is simple to compute, given $\btheta$ and $\tilde x_t$. 
The right hand side, $\Delta  M^H_{t_i}$, are independently distributed as $N\left(0, \upsilon_{M^H_{t_i}}^2\right)$, the normal distribution with mean zero and variance 
$\upsilon_{M^H_{t_i}}^2$, where
\begin{equation}
\upsilon_{M^H_{t_i}}^2=c_2^2(t_{i+1}^{2-2H} - t_{i}^{2-2H}).
\label{eq:var1_m}
\end{equation}
The notation $\Delta  M^H_{t_i}$ aptly indicates that the distribution is independent of $\beta$. With these clarifications, the likelihood of $\btheta$ can be
expressed in terms of the distribution of the right hand side of (\ref{eq:discretize2}). Note that for numerical optimization of the likelihood, or for MCMC purposes,
it is necessary to evaluate the likelihood only for ``proposed" values of $\btheta$, which is a simple exercise, since as demonstrated, computation of the left hand side
of (\ref{eq:discretize2}) is straightforward.

The true value of $\btheta=(\beta,H)$ is denoted by $\btheta_0=(\beta_0,H_0)$. 
We write the likelihood function in terms of the distribution of $\Delta  M^H_{t_i}=\Delta Z^{\btheta}_{t_i} - Q_{\btheta}(t_i)\Delta w_{t_i}^H$. 
Thus, the true likelihood corresponds to $\Delta  M^{H_0}_{t_i}=\Delta Z^{\btheta_0}_{t_i} - Q_{\btheta_0}(t_i)\Delta w_{t_i}^{H_0}$
and the modeled likelihood is associated with $\Delta  M^H_{t_i}=\Delta Z^{\btheta}_{t_i} - Q_{\btheta}(t_i)\Delta w_{t_i}^H$.

Thus, the likelihood given the true parameter $\btheta_0$ is
\begin{equation}
L_{T,n}(\btheta_0)=\prod_{i=1}^n \frac{1}{\sqrt{2\pi}\upsilon_{M^{H_0}_{t_i}}}
	\exp\left\{- \frac{\left(\Delta  M^{H_0}_{t_i}\right)^2}{2\left(\upsilon_{M^{H_0}_{t_i}}\right)^2}\right\}
	=\prod_{i=1}^n \frac{1}{\sqrt{2\pi}\upsilon_{M^{H_0}_{t_i}}} 
	\exp\left\{- \frac{\left(\Delta Z^{\btheta_0}_{t_i} - Q_{\btheta_0}(t_i)\Delta w_{t_i}^{H_0}\right)^2}{2\left(\upsilon_{M^{H_0}_{t_i}}\right)^2}\right\}
\label{eq:likelihood_t}
\end{equation}
 and the modeled likelihood is given by 
\begin{equation}
L_{T,n}(\btheta)=\prod_{i=1}^n \frac{1}{\sqrt{2\pi}\upsilon_{M^{H}_{t_i}}}
	\exp\left\{- \frac{\left(\Delta  M^{H}_{t_i}\right)^2}{2\left(\upsilon_{M^{H}_{t_i}}\right)^2}\right\}
	=\prod_{i=1}^n \frac{1}{\sqrt{2\pi}\upsilon_{M^{H}_{t_i}}} 
	\exp\left\{- \frac{\left(\Delta Z^{\btheta}_{t_i} - Q_{\btheta}(t_i)\Delta w_{t_i}^{H}\right)^2}{2\left(\upsilon_{M^{H}_{t_i}}\right)^2}\right\}.
\label{eq:likelihood_m}
\end{equation}
%
%
The $MLE$ $\hat\btheta_{T,n}=(\hat\beta_{T,n},\hat H_{T,n})$ of $\btheta=(\beta,H)$ is the maximizer of (\ref{eq:likelihood_m}).
There is no closed form for the maximizer, however, and numerical techniques must be employed to approximate the $MLE$ in practice.
For general optimization, we prefer the simulated annealing algorithm, as in every iteration of the underlying Metropolis-Hastings algorithm, 
it can escape the attraction of local modes with positive probability and has elegant
convergence theory; see \ctn{Andrieu01} for example. 
We propose new simulated annealing algorithms, Algorithms \ref{algo:sima} and \ref{algo:sima_real} 
where we replace the Metropolis-Hastings algorithm with the far more efficient Transformation based Markov Chain Monte Carlo (TMCMC) (\ctn{Dutta14}). The essence of TMCMC is to
update all the parameters by simple deterministic transformations of a single, one-dimensional random variable. This yields extremely
fast and efficient simulation algorithms. In practice, it is not difficult to identify a (perhaps large enough) set in which the true maximizer is expected to lie,
and hence, the parameter space may be considered compact. In such cases, the underlying TMCMC algorithm is uniformly ergodic, and the simulated annealing algorithm
converges in probability to the global maximizer. 
Algorithms \ref{algo:sima} and \ref{algo:sima_real} yielded excellent results for our simulation studies as well as the real data analyses. 

Let us now return to theoretical investigation of the asymptotic properties of the $MLE$ $\hat\btheta_{T,n}$.
In what follows, for any two sequences $\left\{a_k:~k\geq 1\right\}$ and $\left\{b_k:~k\geq 1\right\}$,
$a_k\sim b_k$ as $k\rightarrow\infty$ denotes $a_k/b_k\rightarrow 1$, as $k\rightarrow\infty$.
In other words, $a_k\sim b_k$ as $k\rightarrow\infty$ denotes that the sequences $a_k$ and $b_k$ are ``asymptotically equivalent".
Also, ``$\stackrel{a.s.}{\longrightarrow}$" will denote ``almost sure" convergence.

In the context of the definition $t_i=\frac{Ti}{n}$, where $i=1,\ldots,n$, we shall need to deal with limits of the form 
$\underset{n\rightarrow\infty}{\lim}~\frac{1}{n}\sum_{i=1}^nf(i)$, or
$\underset{k\rightarrow\infty}{\lim}~\frac{1}{n_k}\sum_{i=1}^{n_k}f(i)$, where $n_k\rightarrow\infty$, as $k\rightarrow\infty$.
As shown in Lemma \ref{lemma:lemma1}, the limits are easy to deal with if an asymptotic form of $f(i)$, as $i\rightarrow\infty$, can be obtained. 
However, since $i=1,\ldots,n$ (or $n_k$), instead
of writing $i\rightarrow\infty$, we shall denote $i\rightarrow n$ (or $i\rightarrow n_k)$ and $n\rightarrow\infty$ (or $k\rightarrow\infty$). Here, the asymptotic form
of $f(i)$ as $i\rightarrow\infty$ will be the same as that of $f(n)$ (or $f(n_k)$) as $n\rightarrow\infty$ (or $k\rightarrow\infty$). Thus, in general, when
we write $\underset{i\rightarrow n;~n\rightarrow\infty}{\lim}~f(i)$ (or $\underset{i\rightarrow n_k;~k\rightarrow\infty}{\lim}~f(i)$), 
we mean $\underset{n\rightarrow\infty}{\lim}~f(n)$ (or $\underset{k\rightarrow\infty}{\lim}~f(n_k)$). The quantities can be similarly defined when we also
have $T\rightarrow\infty$ (or $T_k\rightarrow\infty$ as $k\rightarrow\infty$).

\section{Strong consistency of the $MLE$ of $\btheta$}
\label{sec:strong_consistency_MLE}
\begin{theorem}
\label{theorem:consistency}
	Let data $\tilde x$ be generated from $P_{\btheta_0}$ and given $x$ consider the $SDE$ given by (\ref{eq:fsde1}) for inference. In $P_{\btheta_0}$ and (\ref{eq:fsde1}), 
	assume $(H1)$ and $(H2)$. 
	Also, let $\bTheta=\mathfrak B\times[\eta,1-\eta]$ be the parameter space for $\btheta=(\beta,H)$, where $\eta>0$ and
	$\mathfrak B$ is a compact subset of the real line $\mathbb R$.
	Further assume  
	that there exist real valued continuous functions $C_1$ and $C_2$ such that for any $H\in[\eta,1-\eta]$,
\begin{align}
	&\int_0^1 u^{\frac{1}{2}-H}(1-u)^{\frac{1}{2}-H}\frac{B(t_iu,\tilde x_{t_{i}u})}{C(t_iu,\tilde x_{t_{i}u})}du\rightarrow C_1(H,\tilde x),\notag\\
	&\qquad\qquad~\mbox{almost surely, as}~
	i\rightarrow n;~n\rightarrow\infty;~T\rightarrow\infty; 
	\label{eq:ass_lim1_general}\\
	&\frac{d}{dt}\left.\int_0^1 u^{\frac{1}{2}- H}(1-u)^{\frac{1}{2}- H}\frac{B(tu,\tilde x_{tu})}{C(tu,\tilde x_{tu})}du\right|_{t=t_{i}}
	\rightarrow C_2(H,\tilde x),\notag\\
	&\qquad\qquad~\mbox{almost surely, as}~i\rightarrow n;~n\rightarrow\infty;~T\rightarrow\infty, 
	\label{eq:ass_lim2_general}
\end{align}
	and that the convergences (\ref{eq:ass_lim1_general}) and (\ref{eq:ass_lim2_general}) hold uniformly for $H\in[\eta,1-\eta]$.
Then the $MLE$ of $\btheta=(\beta,H)$ is strongly consistent in the sense that 
	$\hat\btheta_{T,n} \rightarrow\btheta_0$ almost surely with respect to $P_{\btheta_0}$, as $T\rightarrow\infty$, $n\rightarrow\infty$, such that $n/T^2\rightarrow\infty$.
\end{theorem}

\begin{remark}
	\label{remark:int_examples}
	It is important to provide examples of functions $B(t,u)$ and $C(t,u)$ such that (\ref{eq:ass_lim1_general}) and (\ref{eq:ass_lim2_general}) hold.
	Here we provide a few examples in this regard.
	\begin{itemize}
		\item[(E1)] $B(t,u)=K C(t,u)$ for all $t$ and $u$, where $K$ is some finite constant.
		\item[(E2)] For any $t$, $B(t,u)=K_0 C(t,u)$ if $u\in (-\infty,r_0]$, $B(t,u)=K_i C(t,u)$ if $u\in (r_{i-1},r_i]$ for $i=1,\ldots,L$ and
			$B(t,u)=K_{L+1} C(t,u)$ if $u\in (r_L,\infty)$, where $K_i$; $i=0,1,\ldots,L+1$ and $r_i$; $i=0,1,\ldots,L$, are finite constants, with $L\geq 1$
			being some positive integer.
		\item[(E3)] For sufficiently small $\epsilon>0$, $B(t,u)$ and $C(t,u)$ have different functional forms for $\epsilon< t\leq T_0$, 
			for some sufficiently large $T_0$, but for $0\leq t\leq\epsilon$ 
			and for $t>T_0$, the form of (E1) or (E2) holds.
	\end{itemize}
\end{remark}

In fact, with (E1), (E2) or (E3), we can replace the assumption $n/T^2\rightarrow\infty$ with $n/T\rightarrow\infty$, as the following result shows.
\begin{corollary}
\label{theorem:consistency2}
	Let data $\tilde x$ be generated from $P_{\btheta_0}$ and given $\tilde x$ consider the $SDE$ given by (\ref{eq:fsde1}) for inference. In $P_{\btheta_0}$ and (\ref{eq:fsde1}), 
	assume $(H1)$ and $(H2)$. Also assume either of (E1), (E2) and (E3). 
	If $\bTheta=\mathfrak B\times[\eta,1-\eta]$ is the parameter space for $\btheta=(\beta,H)$, where $\eta>0$ and
	$\mathfrak B$ is a compact subset of the real line $\mathbb R$,
	then $\hat\btheta_{T,n} \rightarrow\btheta_0$ almost surely with respect to $P_{\btheta_0}$, as $T\rightarrow\infty$, $n\rightarrow\infty$, such that $n/T\rightarrow\infty$.
\end{corollary}

\section{Asymptotic normality of $\bbeta$}
\label{sec:normality_beta}

\begin{theorem}
\label{theorem:asymp_normal1}
	Let data $\tilde x$ be generated from $P_{\btheta_0}$ and given $\tilde x$ consider the $SDE$ given by (\ref{eq:fsde1}) for inference. In $P_{\btheta_0}$ and (\ref{eq:fsde1}), 
	assume $(H1)$ and $(H2)$. 
	Also, let $\bTheta=\mathbb R\times[\eta,1-\eta]$ be the parameter space for $\btheta=(\beta,H)$. 
	Further assume that there exist real valued functions $C_1$ and $C_2$ such that for any $H\in[\eta,1-\eta]$, where $\eta>0$,
	(\ref{eq:ass_lim1_general}) and (\ref{eq:ass_lim2_general}) hold.
	Then almost surely with respect to $P_{\btheta_0}$,
\begin{equation*}
	\frac{T^{2-\hat H_{T,n}}}{\left|a_{\hat H_{T,n}}\right|}\left(\hat\beta_{T,n}-\beta_0\right)\rightarrow \tilde Z\sim N(0,1),	
	\mbox{as}~T\rightarrow\infty;~n\rightarrow\infty;~\frac{n}{T^2}\rightarrow\infty.
	\label{eq:normal7}
\end{equation*}
	where, for any $H\in[\eta,1-\eta]$, $a_{H}=\frac{k_{H}c_{H}}{2HC_2(H,x)}$.
\end{theorem}

As in the case for consistency, either of the assumptions (E1), (E2) or (E3), allows replacement of $T^{2-\hat H_{T,n}}$ in the asymptotic normality
of Theorem \ref{theorem:asymp_normal1} with $T^{1-\hat H_{T,n}}$, which is associated with $n/T\rightarrow\infty$ instead of $n/T^2\rightarrow\infty$
associated with the former rate.
\begin{theorem}
\label{theorem:asymp_normal2}
	Let data $\tilde x$ be generated from $P_{\btheta_0}$ and given $\tilde x$ consider the $SDE$ given by (\ref{eq:fsde1}) for inference. In $P_{\btheta_0}$ and (\ref{eq:fsde1}), 
	assume $(H1)$ and $(H2)$. 
	Also, let $\bTheta=\mathbb R\times[\eta,1-\eta]$ be the parameter space for $\btheta=(\beta,H)$, where $\eta>0$. 
	Further assume either of (E1), (E2) or (E3). 
	Then almost surely with respect to $P_{\btheta_0}$,
\begin{equation*}
	\frac{T^{1-\hat H_{T,n}}}{\left|\alpha_{\hat H_{T,n}}\right|}\left(\hat\beta_{T,n}-\beta_0\right)\rightarrow \tilde Z\sim N(0,1),	
	~\mbox{as}~T\rightarrow\infty;~n\rightarrow\infty;~\frac{n}{T}\rightarrow\infty.
\end{equation*}
	where, for any $H\in[\eta,1-\eta]$, $\alpha_{H}=\frac{k_{H}c_{H}}{4H(1-H)C_1(H)}$.
\end{theorem}

\begin{remark}
\label{remark:classical}
It is worth noting that Theorems \ref{theorem:asymp_normal1} and \ref{theorem:asymp_normal2} show that for asymptotic normality of $\beta$, 
the parameter space is no longer needed to be compact, unlike Theorems \ref{theorem:consistency} and \ref{theorem:consistency2}. 
Here the entire real line is allowed for $\beta$.
These theorems 
also ensure that $\hat\beta_{T,n}\stackrel{a.s.}{\longrightarrow}\beta_0$.
\end{remark}

\section{Inconsistency of the $MLE$ of $\beta$ when the estimator of $H$ is inconsistent}
\label{sec:incons_beta}
Remark \ref{remark:classical} shows that consistency of $\hat\beta_{T,n}$ continues to hold even without the assumption that $\mathfrak B$ is compact. 
However, if $H$ is estimated with an inconsistent estimator, say, $H^*_{T,n}$ which does not converge to $H_0$, then even consistency does not 
hold for the $MLE$ of $\beta$ for $\mathfrak B$ either compact or non-compact. 
Theorems \ref{theorem:incons_noncompact1} and \ref{theorem:incons_noncompact2} formalize this in the situations with and without assumption of either (E1), (E2), (E3).
Notably, the results hold simply as $T\rightarrow\infty$ and $n\rightarrow\infty$, irrespective of whether or not $n$ grows faster than $T$.

	\begin{theorem}
        \label{theorem:incons_noncompact1}
		Let data $\tilde x$ be generated from $P_{\btheta_0}$ and given $\tilde x$ consider the $SDE$ given by (\ref{eq:fsde1}) for inference. In $P_{\btheta_0}$ and (\ref{eq:fsde1}), 
		assume $(H1)$ and $(H2)$. 
	Also, let $\bTheta=\mathbb R\times[\eta,1-\eta]$ be the parameter space for $\btheta=(\beta,H)$. 
	Further assume that there exist real valued functions $C_1$ and $C_2$ such that for any $H\in[\eta,1-\eta]$, where $\eta>0$,
	(\ref{eq:ass_lim1_general}) and (\ref{eq:ass_lim2_general}) hold. 
		For each $T>0$ and $n\geq 1$, let $H$ be estimated by $H^*_{T,n}$, where $H^*_{T,n}\stackrel{a.s.}{\longrightarrow}\tilde H$ with respecto to $P_{\btheta_0}$, 
		with $\tilde H\neq H_0$.
		Then with respect to $P_{\btheta_0}$, the $MLE$ of $\beta$ is almost surely an inconsistent estimator of $\beta_0$, for $\mathfrak B$ either compact or non-compact. 
	\end{theorem}

\begin{theorem}
        \label{theorem:incons_noncompact2}
	Let data $\tilde x$ be generated from $P_{\btheta_0}$ and given $\tilde x$ consider the $SDE$ given by (\ref{eq:fsde1}) for inference. In $P_{\btheta_0}$ and (\ref{eq:fsde1}), 
	assume $(H1)$ and $(H2)$. 
	Also, let $\bTheta=\mathbb R\times[\eta,1-\eta]$ be the parameter space for $\btheta=(\beta,H)$. 
	Further assume either of (E1), (E2) or (E3). 
	For each $T>0$ and $n\geq 1$, let $H$ be estimated by $H^*_{T,n}$, where $H^*_{T,n}\stackrel{a.s.}{\longrightarrow}\tilde H$ with respect to to $P_{\btheta_0}$, 
	with $\tilde H\neq H_0$.
	Then the $MLE$ of $\beta$ is almost surely an inconsistent estimator of $\beta_0$, 
	for $\mathfrak B$ either compact or non-compact. 
	\end{theorem}

\begin{remark}
	\label{remark:inconsistency}
	Theorems \ref{theorem:incons_noncompact1} and \ref{theorem:incons_noncompact2} show that if the estimator of $H$ is inconsistent, then this induces
	inconsistency in the $MLE$ of $\beta$. Thus, these results
	caution against simple-minded usage of estimators of $H$ that lack theoretical
	validation of consistency. Hence, the heuristic estimators of $H$, discussed in Section \ref{sec:discussion} need not be automatic choices, as generally
	treated by practitioners. Since we showed strong consistency of the $MLE$ of $H$ in Theorems \ref{theorem:consistency} and 
	\ref{theorem:consistency2}, we recommend $MLE$ as the appropriate estimator of $H$.
\end{remark}

\section{Bayesian posterior consistency}
\label{sec:bayesian_consistency}

For Bayesian inference on $\btheta$, with respect to any given prior density $\pi(\btheta)$ (with respect to the Lebesgue measure on the parameter space $\bTheta$) 
for $\btheta$, it is required to obtain the posterior distribution of $\btheta$, given the data $\tilde x_{t_i}$; $i=1,\ldots,n$. In our case, the posterior density with respect to the
Lebesgue measure, which we denote, slightly
abusing notation, by $\pi(\btheta|\tilde x_{t_i};~i=1,\ldots,n)$, is of the form
\begin{equation}
	\pi(\btheta|\tilde x_{t_i};~i=1,\ldots,n)\propto\pi(\btheta)L_{T,n}(\btheta),
	\label{eq:posterior_theta}
\end{equation}
where $L_{T,n}(\btheta)$ is given by (\ref{eq:likelihood_m}).
The above posterior of $\btheta$ does not admit any closed form, and hence, for practical purposes, sampling-based inference is recommended.
Since even direct sampling from the posterior is difficult, sampling based on Markov Chain Monte Carlo (MCMC) methods may be considered. Algorithm \ref{algo:tmcmc}
of the supplement is an efficient TMCMC algorithm in this regard, that suits our purpose, providing
reliable Bayesian inference, as all our simulation experiments and real data analyses vindicate. 

From the asymptotic standpoint, assuming that the data $x$ is generated from $P_{\btheta_0}$ and that $\tilde x_{t_i}$; $i=1,\ldots,n$ is a discretized version of $\tilde x$,
we can expect that the posterior (\ref{eq:posterior_theta}) concentrates more and more around $\btheta_0$, as both $n$ and $T$ increases, with $n$ growing at a faster rate
than $T$. If this holds in some precise sense, the posterior (\ref{eq:posterior_theta}) is said to be consistent at $\btheta_0$. In this section, we formally investigate
posterior consistency.

For any set $A$ in the Borel sigma-algebra of the parameter space $\bTheta$, let $I_A(\btheta)=1$ if $\btheta\in A$ and $0$, otherwise. 
\begin{theorem}
\label{theorem:bayesian_consistency}
	Let data $\tilde x$ be generated from $P_{\btheta_0}$ and given $\tilde x$ consider the $SDE$ given by (\ref{eq:fsde1}) for inference. In $P_{\btheta_0}$ and (\ref{eq:fsde1}), 
	assume $(H1)$ and $(H2)$. 
	Also, let $\bTheta=\mathbb R\times[\eta,1-\eta]$ be the parameter space for $\btheta=(\beta,H)$, where $\eta>0$.
	Further assume (\ref{eq:ass_lim1_general}) and (\ref{eq:ass_lim2_general}). 
Then the posterior distribution of $\btheta=(\beta,H)$ is strongly consistent in the sense that for any prior on $\bTheta$ that includes the true value
$\btheta_0$ in its support,
	$\pi\left(\btheta\in A|\tilde x_{t_i};~i=1,\ldots,n\right)\stackrel{a.s.}{\longrightarrow}I_A(\btheta_0)$ with respect to $P_{\btheta_0}$, as 
	$T\rightarrow\infty$, $n\rightarrow\infty$, such that $n/T^2\rightarrow\infty$, for any set $A$ in the Borel sigma-algebra of $\bTheta$.
\end{theorem}

Under (E1), (E2) or (E3), the following holds for $n/T\rightarrow\infty$, as to be expected. 
\begin{theorem}
\label{theorem:bayesian_consistency2}
	Let data $\tilde x$ be generated from $P_{\btheta_0}$ and given $\tilde x$ consider the $SDE$ given by (\ref{eq:fsde1}) for inference. In $P_{\btheta_0}$ and (\ref{eq:fsde1}), 
	assume $(H1)$ and $(H2)$. Also assume either of (E1), (E2) and (E3).
	If $\bTheta=\mathbb R\times[\eta,1-\eta]$ is the parameter space for $\btheta=(\beta,H)$, where $\eta>0$, 
then the posterior distribution of $\btheta=(\beta,H)$ is strongly consistent in the sense that for any prior on $\bTheta$ that includes the true value
$\btheta_0$ in its support,
	$\pi\left(\btheta\in A|\tilde x_{t_i};~i=1,\ldots,n\right)\stackrel{a.s.}{\longrightarrow}I_A(\btheta_0)$ with respect to $P_{\btheta_0}$, as 
	$T\rightarrow\infty$, $n\rightarrow\infty$, such that $n/T\rightarrow\infty$, for any set $A$ in the Borel sigma-algebra of $\bTheta$.
\end{theorem}

Note that although Doob's theorem ensures consistency of the posterior distributions, it does not specify which neighborhoods the posteriors concentrate on.
The theorems below show the type of neighborhoods where the posterior of $\beta$, given $H=\hat H_{T,n}$, concentrates on. 

\begin{theorem}
\label{theorem:bayesian_consistency_KL}
	Let data $\tilde x$ be generated from $P_{\btheta_0}$ and given $\tilde x$ consider the $SDE$ given by (\ref{eq:fsde1}) for inference. In $P_{\btheta_0}$ and (\ref{eq:fsde1}), 
	assume $(H1)$ and $(H2)$. 
	Further assume (\ref{eq:ass_lim1_general}) and (\ref{eq:ass_lim2_general}). 
	Let $N_0$ be any neighborhood of $\beta_0$ containing $C_{\epsilon}=\left\{\beta:(\beta-\beta_0)^2g(H_0)<\epsilon\right\}$, for any $\epsilon>0$,
	where, for any $H\in[\eta,1-\eta]$, 
	$g(H)=\frac{4H^2}{c^2_{H}}\times\left\{k^{-1}_{H}C_2(H,\tilde x)\right\}^2$.
	Let $\pi(\beta)$ be any Lebesgue-dominated prior for $\beta$ that satisfies $\pi(C_{\epsilon})>0$, for any $\epsilon>0$. 
	Then, the posterior distribution of $\beta\in\mathbb R$  
	given $H=\hat H_{T,n}$ satisfies
	$$\pi\left(\beta\in N_0|H=\hat H_{T,n},\tilde x_{t_i};~i=1,\ldots,n\right)\stackrel{a.s.}{\longrightarrow}1,$$ with respect to $P_{\btheta_0}$, as 
	$T\rightarrow\infty$, $n\rightarrow\infty$, such that $n/T^2\rightarrow\infty$.
\end{theorem}

We now consider the case $n/T\rightarrow\infty$; the result is summarized below as Theorem \ref{theorem:bayesian_consistency2_KL}.
\begin{theorem}
\label{theorem:bayesian_consistency2_KL}
	Let data $\tilde x$ be generated from $P_{\btheta_0}$ and given $\tilde x$ consider the $SDE$ given by (\ref{eq:fsde1}) for inference. In $P_{\btheta_0}$ and (\ref{eq:fsde1}), 
	assume $(H1)$ and $(H2)$. Further assume either (E1), (E2) or (E3).
	Let $N_0$ be any neighborhood of $\beta_0$ containing $C_{\epsilon}=\left\{\beta:(\beta-\beta_0)^2g(H_0)<\epsilon\right\}$, for any $\epsilon>0$,
	where, for any $H\in[\eta,1-\eta]$, 
	$g(H)=\frac{4H^2}{c^2_{H}}\times\left\{k^{-1}_{H}(2-2H)C_1(H)\right\}^2$.
	Let $\pi(\beta)$ be any Lebesgue-dominated prior for $\beta$ that satisfies $\pi(C_{\epsilon})>0$, for any $\epsilon>0$. 
	Then, the posterior distribution of $\beta\in\mathbb R$  
	given $H=\hat H_{T,n}$ satisfies
	$\pi\left(\beta\in N_0|H=\hat H_{T,n},\tilde x_{t_i};~i=1,\ldots,n\right)\stackrel{a.s.}{\longrightarrow}1$ with respect to $P_{\btheta_0}$, as 
	$T\rightarrow\infty$, $n\rightarrow\infty$, such that $n/T\rightarrow\infty$.
\end{theorem}

\section{Asymptotic posterior normality}
\label{sec:posterior_normaility}
For asymptotic posterior normality of $\beta$ given $H=\hat H_{T,n}$, we make use of Theorem 7.89 of \ctn{Schervish95}. 
In Section \ref{sec:schervish}, we state the general setup, regularity conditions referred to as (1) -- (6), and the main result presented in \ctn{Schervish95}. 
It is worth noting that the setup, conditions and the result
are not established for the purpose of increasing domain infill asymptotics as in our setup, namely, $T\rightarrow\infty$, $n\rightarrow\infty$ such that
$n/T^2\rightarrow\infty$ or $n/T\rightarrow\infty$. Nevertheless, our setup can be adapted to that of \ctn{Schervish95}, which continues to hold
even in the increasing domain infill asymptotics situation.

\subsection{Asymptotic posterior normality for fractional $SDE$}
\label{sec:post_normal1}
In our case of increasing domain infill asymptotics when $n/T^2\rightarrow\infty$, we obtain the following result using the result in \ctn{Schervish95}:
\begin{theorem}
\label{theorem:post_normal1}
	Let data $\tilde x$ be generated from $P_{\btheta_0}$ and given $\tilde x$ consider the $SDE$ given by (\ref{eq:fsde1}) for inference. In $P_{\btheta_0}$ and (\ref{eq:fsde1}), 
	assume $(H1)$ and $(H2)$. 
Further assume (\ref{eq:ass_lim1_general}) and (\ref{eq:ass_lim2_general}). 
	Let $\Psi_{T,n}=\frac{T^{2-\hat H_{T,n}}}{\left|a_{\tilde H_{T,n}}\right|}\left(\beta-\beta_0\right)$, where, for any $H\in[\eta,1-\eta]$,
and
\begin{equation}
	a_{H}=\frac{c_{H}k_{H}}{2HC_2(H,\tilde x)}.
	\label{eq:norm2}
\end{equation}
	Let $\pi_{T,n}$ denote the posterior distribution of $\Psi_{T,n}$ given $\tilde x_{t_i};i=1,\ldots,n$, and $H=\hat H_{T,n}$, with respect to any Lebesgue-dominated prior.
Then for each compact subset $B$ of $\mathbb R$ and each $\epsilon>0$, the following holds:
\begin{equation*}
	\lim_{T\rightarrow\infty,n\rightarrow\infty,\frac{n}{T^2}\rightarrow\infty}P_{\btheta_0}
	\left(\sup_{\psi\in B}\left\vert\pi_{T,n}(\psi)-\phi(\psi)\right\vert>\epsilon\right)= 0.
\end{equation*}
where $\phi(\cdot)$ denotes the density of the one-dimensional standard normal distribution.
\end{theorem}

The following theorem is regarding asymptotic posterior normality under either of (E1), (E2) and (E3), which allows $n/T\rightarrow\infty$, 
instead of $n/T^2\rightarrow\infty$, as before.
\begin{theorem}
\label{theorem:post_normal2}
	Let data $\tilde x$ be generated from $P_{\btheta_0}$ and given $\tilde x$ consider the $SDE$ given by (\ref{eq:fsde1}) for inference. In $P_{\btheta_0}$ and (\ref{eq:fsde1}), 
	assume $(H1)$ and $(H2)$. Further assume either of (E1), (E2) and (E3).
	Let $\Psi_{T,n}=\frac{T^{1-\hat H_{T,n}}}{\left|\alpha_{\hat H_{T,n}}\right|}\left(\beta-\beta_0\right)$, where, for any $H\in[\eta,1-\eta]$,
and
\begin{equation*}
	\alpha_{H}=\frac{c_{H}k_{H}}{4H(1-H)C_1(H)}.
\end{equation*}
	Let $\pi_{T,n}$ denote the posterior distribution of $\Psi_{T,n}$ given $\tilde x_{t_i};i=1,\ldots,n$, and $H=\hat H_{T,n}$, with respect to any Lebesgue-dominated prior.
Then for each compact subset $B$ of $\mathbb R$ and each $\epsilon>0$, the following holds:
\begin{equation*}
	\lim_{T\rightarrow\infty,n\rightarrow\infty,\frac{n}{T}\rightarrow\infty}P_{\btheta_0}
	\left(\sup_{\psi\in B}\left\vert\pi_{T,n}(\psi)-\phi(\psi)\right\vert>\epsilon\right)= 0.
\end{equation*}
where $\phi(\cdot)$ denotes the density of the one-dimensional standard normal distribution.
\end{theorem}

\section{Simulation study}
\label{sec:simstudy}
We perform simulation experiments by generating data from a specific fractional $SDE$ as follows:
\begin{equation}
d\tilde X_t=\beta(1-\tilde X_t)dt + (1-\tilde X_t)d\tilde W_t^ H.
	\label{eq:sde_example}
\end{equation}
For simulating data from $SDE$ of the general form (\ref{eq:fsde0}), fixing $\beta=\beta_0$ and $H=H_0$, we consider $\tilde X_0=\tilde x_0$, where $\tilde x_0$ is any fixed, non-random quantity,
and with $t_i=\frac{Ti}{n}$; $i=0,1,\ldots,n$, set $\tilde X_{t_{i+1}}-\tilde X_{t_i}=\beta B(t_i,\tilde X_{t_i})(t_{i+1}-t_i)+C(t_i,\tilde X_{t_i})(\tilde W^H_{t_{i+1}}-\tilde W^H_{t_i})$, for $i=0,1,\ldots,n$.
Since independently, the increments $\tilde W^H_{t_{i+1}}-\tilde W^H_{t_i}\sim N\left(0,(t_{i+1}-t_i)^{2H}\right)$, we first simulate these increments, which we then plug in the
above discretization to obtain $\tilde x_{t_i}$; $i=1,\ldots,n$. 
Once the data is so generated, we fit $SDE$ of the form (\ref{eq:fsde1}) to the data. 
%

We fix $T=5$, $n=100$ for our experiments, so that the non-asymptotic flavour can be retained. 
We conduct simulation experiments in both classical and Bayesian setups, detailed in Sections \ref{subsec:classical_simstudy} and \ref{subsec:bayesian_simstudy}, respectively.

\subsection{Simulation study in the classical setup}
\label{subsec:classical_simstudy}

We perform simulated annealing aided by transformation based Markov chain Monte Carlo (TMCMC) (\ctn{Dutta14})
to obtain the $MLE$ of $\btheta=(\beta,H)$. We provide our algorithm  as Algorithm \ref{algo:sima}, where we use the additive transformation. 
\begin{algo}
\topline{Simulated annealing aided by additive TMCMC}
\label{algo:sima}
\botline
\begin{itemize}
	\item[(1)] Begin with an initial value $\btheta^0=\left(\beta^{(0)},H^{(0)}\right)$.
\item[(2)] For $t=1,\ldots,N$, do the following:
	\begin{enumerate}
		\item[(i)] Generate $\varepsilon\sim N(0,1)$, $b_j\stackrel{iid}{\sim}U(\{-1,1\})$ for $j=1,2$, and set
			$\tilde\beta=\beta^{(t-1)}+b_1a_1|\varepsilon|$ and $\tilde H=H^{(t-1)}+b_2a_2|\varepsilon|$. 
			Let $\tilde\btheta=\left(\tilde\beta,\tilde H\right)$.
		\item[(ii)] Evaluate 
			$$\alpha=\min\left\{1,\exp\left\{\frac{\ell_{T,n}\left(\tilde\btheta\right)-\ell_{T,n}\left(\btheta^{(t-1)}\right)}{\tau_t}\right\}
			I_{(0,1)}(\tilde H)\right\},$$
			where $\tau_t=\frac{1}{\log(\log(\log(t)))}$ when $t\geq 3$ and $t^{-1}$ otherwise.
			In the above, $I_{(0,1)}(\tilde H)=1$ if $\tilde H\in (0,1)$ and zero otherwise.
		\item[(iii)] Set $\btheta^{(t)}=\tilde\btheta$ with probability $\alpha$, else set $\btheta^{(t)}=\btheta^{(t-1)}$.
	\end{enumerate}
\item[(3)] From the sample $\left\{\btheta^{(j)};j=0,1,\ldots,N\right\}$, set $\hat\btheta_{T,n}=\btheta^{(j^*)}$ where 
	$\btheta^{(j^*)}=\underset{j=0,1,\ldots,N}{\arg\max}~\ell_{T,n}\left(\btheta^{(j)}\right)$.
\end{itemize}
\botline
\end{algo}
For our experiments, setting $N=20000$ and $a_1=a_2=0.0001$ produced $MLE$s that are reasonably close to the true values. For the initial value $\btheta^{(0)}$,
we conduct a grid search on a grid consisting of equispaced values of $\beta$ and $H$; the maximizer of the log-likelihood on the grid
is taken as $\btheta^{(0)}$.

It is also important to obtain the distribution of the $MLE$ so that desired confidence intervals can be obtained. For $\hat\beta_{T,n}$, at least asymptotic normality
is provided by Theorem \ref{theorem:asymp_normal2} of MB, but for $\hat H_{T,n}$ even asymptotic normality is not available. In any case, it is important to provide
the distribution of the $MLE$ for any sample size. We numerically approximate the distribution by generating $1000$ data samples from the $SDE$ (\ref{eq:sde_example})
and then computing the $MLE$ of $\btheta$ in each case using Algorithm (\ref{algo:sima}). The results of our experiments are tabulated in Table \ref{table:classical}.
The 95\% confidence intervals ($CI$s) are obtained by setting the lower and upper ends to be the 2.5\% and 97.5\% quantiles obtained from the $1000$ realizations of the $MLE$s.
\begin{table}[h]
\centering
\caption{Results of our simulation study in the classical setup.}
\label{table:classical}
	\begin{tabular}{|c|c|c||c|c|c|}
\hline
$H_0$ & $\hat H_{T,n}$ & 95\% $CI$ & $\beta_0$ & $\hat \beta_{T,n}$ & 95\% $CI$\\ 
\hline
	0.1  &  0.095  &  [0.086,0.116]  &  1.0  &  0.899  &  [0.773,1.035]\\
	0.2  &  0.191  &  [0.101,0.262]  &  0.3  &  0.290  &  [0.145,0.360]\\
	0.3  &  0.289  &  [0.232,0.296]  &  -1.0  &  -0.901  &  [-0.933,-0.725]\\
	0.4  &  0.391  &  [0.315,0.452]  &  -0.5  & -0.490  &  [-0.562,-0.374]\\   
	0.5  &  0.489  &  [0.454,0.547]  &  0.7   &  0.689  &  [0.620,0.792]\\ 
	0.6  &  0.590  &  [0.583,0.600]  &  -1.2  &  -1.191 &  [-1.225,-1.190]\\
	0.7  &  0.675  &  [0.674,0.686]  &  1.7   &  1.695  &  [1.691,1.730]\\
	0.8  &  0.773  &  [0.763,0.777]  &  -1.5  &  -1.522  &  [-1.538,-1.494]\\
	0.9  &  0.856  &  [0.851,0.865]  &  1.4   &  1.415    &  [1.396,1.439]\\
\hline
\end{tabular}
\end{table}
Note that most of the $MLE$s $\hat H_{T,n}$ and $\hat \beta_{T,n}$ are reasonably close to the true values $H_0$ and $\beta_0$, respectively. Also, most of the
95\% $CI$s contain the true values. In the third case, where $H_0=0.3$ and $\beta_0=-1$, the 95\% CI of $H$ marginally misses $H_0$, and that of $\beta$ misses
$\beta_0$ by a larger margin. In cases $8$ and $9$, where $(H_0,\beta_0)=(0.8,-1.5)$ and $(0.9,1.4)$, respectively, the variabilities of the distribution of $\hat H_{T,n}$
seem to be too small to include the true value. These apparent failures demonstrate the inadequacy of the small size ($1000$) of the number of realizations from the
distribution of the $MLE$s. However, obtaining even $1000$ realizations implemented in the $R$ software, takes about $5$ hours and $40$ minutes on our ordinary laptop, 
precluding accurate evaluation of the distributions of the $MLE$s. Moreover, the choices of $a_1=a_2=0.0001$ need not be adequate for all the $1000$ samples 
and for all the nine cases shown in Table \ref{table:classical}. This problem of missing the true value from the desired confidence intervals 
can be avoided in the Bayesian setup, as we shall demonstrate in Section \ref{subsec:bayesian_simstudy}.

\subsection{Simulation study in the Bayesian setup}
\label{subsec:bayesian_simstudy}

For the Bayesian investigation we focus on the same cases as shown in Table \ref{table:classical}, and analyze the same datasets for valid comparison between
the classical and Bayesian paradigms. 

In the Bayesian setup, we first reparameterize $H$ as $H=\exp(\nu)/(1+\exp(\nu))$, where we assume {\it a priori} that $\nu\sim N\left(\nu_{T,n},\sigma^2_{\nu}\right)$.
We set $\nu_{T,n}=\log\left(\frac{\hat H_{T,n}}{1-\hat H_{T,n}}\right)$ and $\sigma_{\nu}=0.2$ for the first seven cases and $0.02$ for the last two cases.
We also assume {\it a priori} that $\beta\sim N\left(\hat\beta_{T,n},\sigma^2_{\beta}\right)$, where $\sigma_{\beta}=0.1$ in the first seven cases and $0.01$ in the 
last two cases. These choices quantify the prior belief that the true value $\btheta_0$ is located around the $MLE$ $\hat\btheta_{T,n}$. We considered a larger variance
for $H$ so that it has lesser risk of missing the true value compared to that for $\beta$. Indeed, Table \ref{table:classical} for the classical  
paradigm demonstrates that without any prior, $H$ is more at risk of missing the true values compared to $\beta$. In the Bayesian setup we exploited this wisdom gained
in the form of the priors to counter the problem of omitting the true values. For cases $8$ and $9$, we set much smaller prior variances, so that even in these extreme cases
the posteriors have less propensity to concentrate around less extreme values. Note that even though this is a simulation study where the true values are known which makes it
easy to set the prior variances as explained, in reality the $MLE$ is expected to be close to the truth, and hence even in practice, such prior variances are easy to
envisage. It is also important to remark that a $Beta$ prior on $H$ seems to be natural, but it does not have the flexibilities of setting the mean and variance as desired,
compared to our normal prior distribution.

For Bayesian inference, we propose an additive TMCMC algorithm, enriched with the provision for further enhancement of mixing properties.
The details are provided as Algorithm \ref{algo:tmcmc}. 
%
\begin{algo}
\topline{TMCMC for $SDE$ driven by fractional Brownian motion}
\label{algo:tmcmc}
\botline 
\begin{itemize}
	\item[(1)] Begin with an initial value $\btheta^0=\left(\beta^{(0)},\nu^{(0)}\right)$; set $H^{(0)}=\exp(\nu^{(0)})/(1+\exp(\nu^{(0)}))$.
\item[(2)] For $t=1,\ldots,N$, do the following:
	\begin{enumerate}
		\item[(i)] Generate $\varepsilon_1\sim N(0,1)$, $b_j\stackrel{iid}{\sim}U(\{-1,1\})$ for $j=1,2$, and set
			$\tilde\beta=\beta^{(t-1)}+b_1a_1|\varepsilon_1|$ and $\tilde\nu=\nu^{(t-1)}+b_2a_2|\varepsilon_1|$. 
			Let $\tilde\btheta=\left(\tilde\beta,\tilde H\right)$, where $\tilde H=\exp(\tilde\nu)/(1+\exp(\tilde\nu))$.
		\item[(ii)] Evaluate 
			$$\alpha_1=\min\left\{1,\exp\left\{\ell_{T,n}\left(\tilde\btheta\right)-\ell_{T,n}\left(\btheta^{(t-1)}\right)\right\}\right\}.$$
		\item[(iii)] Set $\btheta^{(t)}=\tilde\btheta$ with probability $\alpha_1$, else set $\btheta^{(t)}=\btheta^{(t-1)}$. 
		\item[(iv)] Generate $\varepsilon_2\sim N(0,1)$ and $U\sim U(0,1)$. If $U<1/2$, create $\beta^*=\beta^{(t)}+a_1|\varepsilon_2|$
			and $H^*=\exp(\nu^*)/(1+\exp(\nu^*))$, where $\nu^*=\nu^{(t)}+b_2a_2|\varepsilon_2|$. Else, set $\beta^*=\beta^{(t)}-a_1|\varepsilon_2|$
			and $H^*=\exp(\nu^*)/(1+\exp(\nu^*))$, where $\nu^*=\nu^{(t)}-a_2|\varepsilon_2|$. Let $\btheta^*=\left(\beta^*,H^*\right)$.
		\item[(v)] Evaluate 
			$$\alpha_2=\min\left\{1,\exp\left\{\ell_{T,n}\left(\btheta^*\right)-\ell_{T,n}\left(\btheta^{(t)}\right)\right\}\right\}.$$
		\item[(vi)] With probability $\alpha_2$ accept $\btheta^*$, else remain at $\btheta^{(t)}$. Abusing notation, let us continue to denote
			either $\btheta^*$ or $\btheta^{(t)}$ as $\btheta^{(t)}$, to be understood as the final value of each iteration.
	\end{enumerate}
\item[(3)] Store the realizations $\left\{\btheta^{(t)};t=0,1,\ldots,N\right\}$ for Bayesian inference. 
\end{itemize}
\botline
\end{algo}
Note that steps (2) (iv) -- (2) (vi) is actually not required for convergence of TMCMC but often has the important advantage of improving mixing properties of
the underlying Markov chain. See \ctn{Liu00} and the supplement of \ctn{Dutta14} for details. 

\begin{figure}
\centering
\subfigure []{ \label{fig:h1}
\includegraphics[width=5.5cm,height=4.5cm]{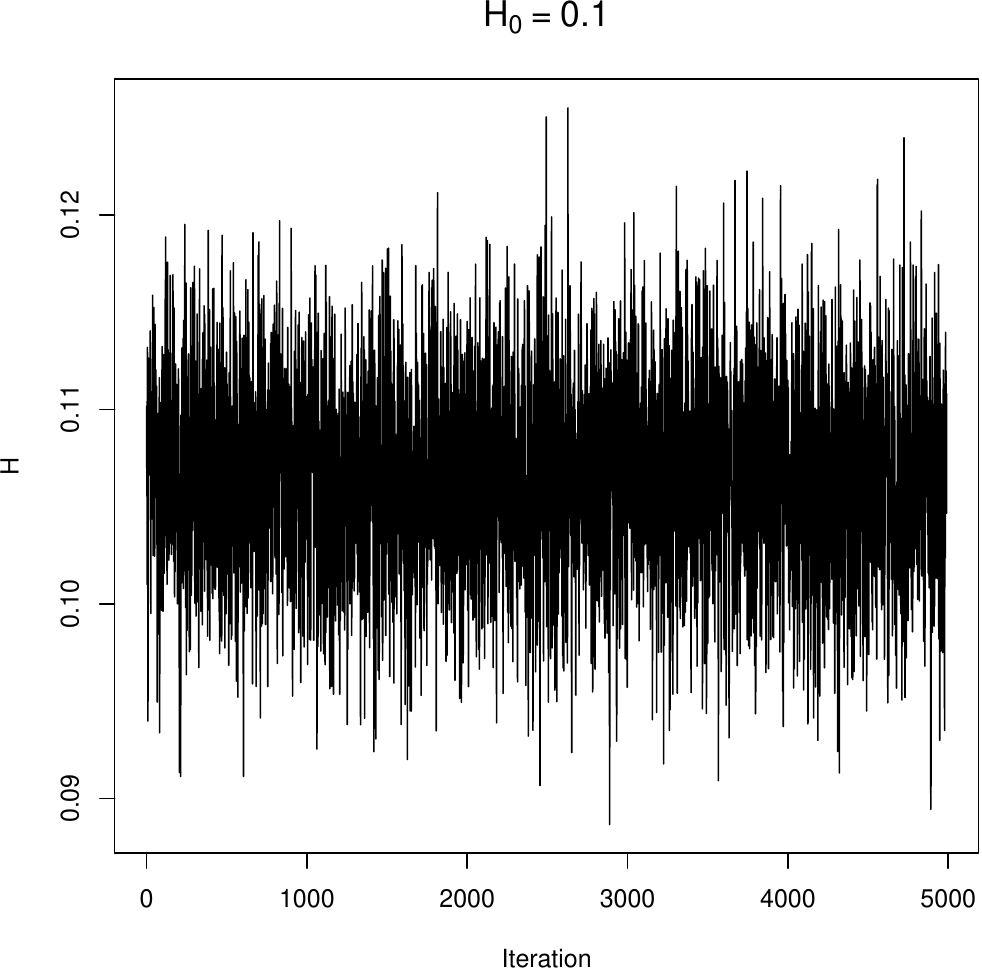}}
\hspace{2mm}
\subfigure []{ \label{fig:b1}
\includegraphics[width=5.5cm,height=4.5cm]{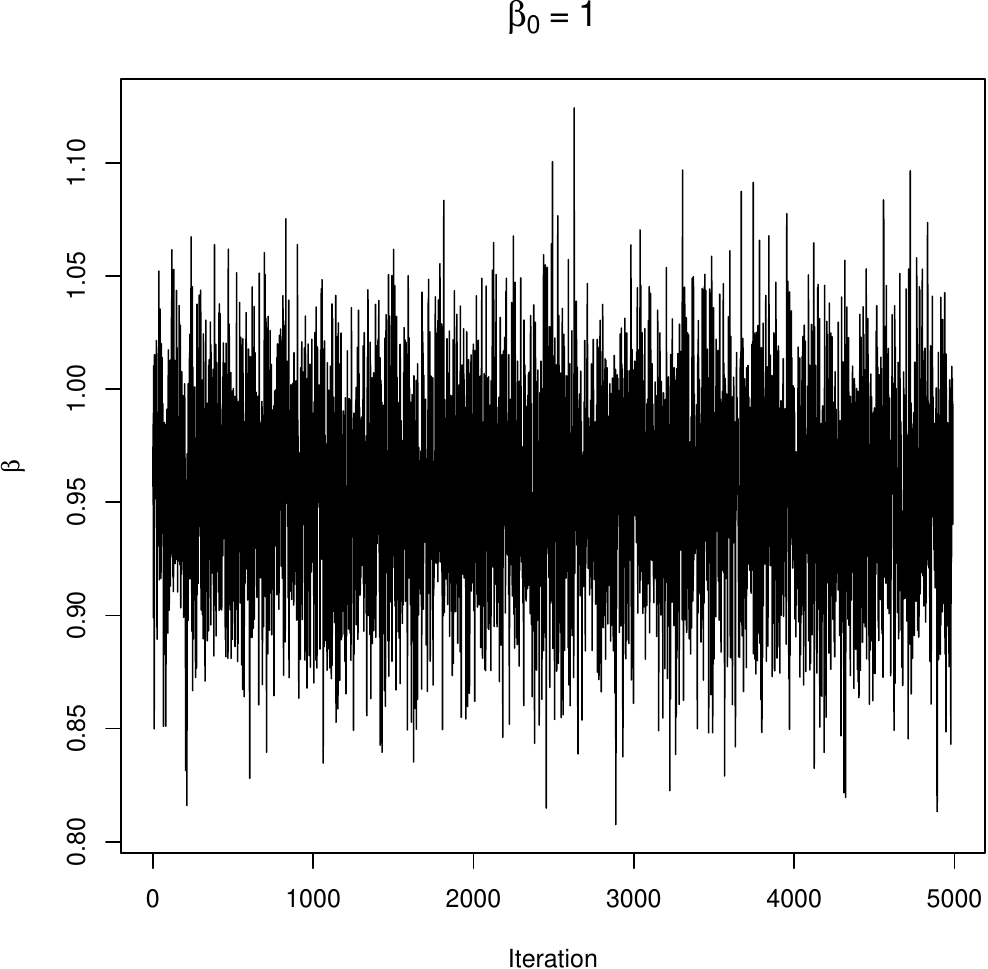}}\\
\vspace{2mm}
\subfigure []{ \label{fig:h2}
\includegraphics[width=5.5cm,height=4.5cm]{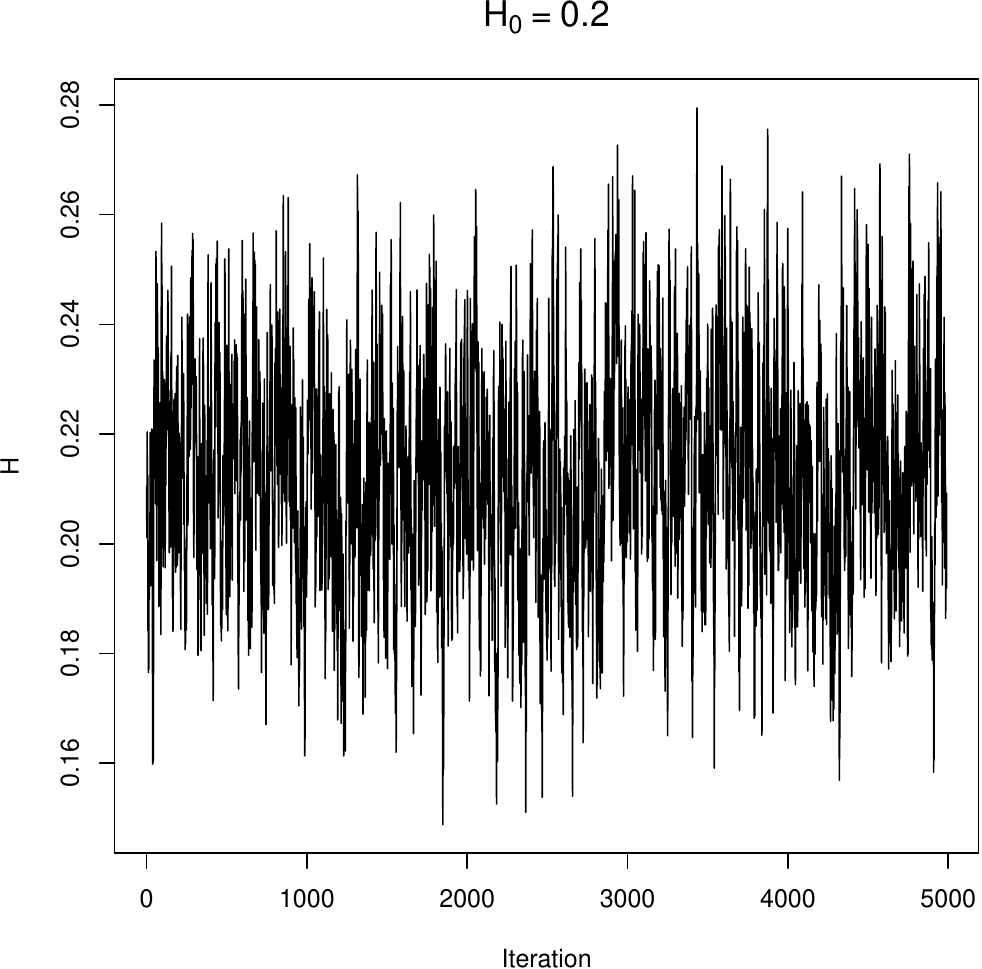}}
\hspace{2mm}
\subfigure []{ \label{fig:b2}
\includegraphics[width=5.5cm,height=4.5cm]{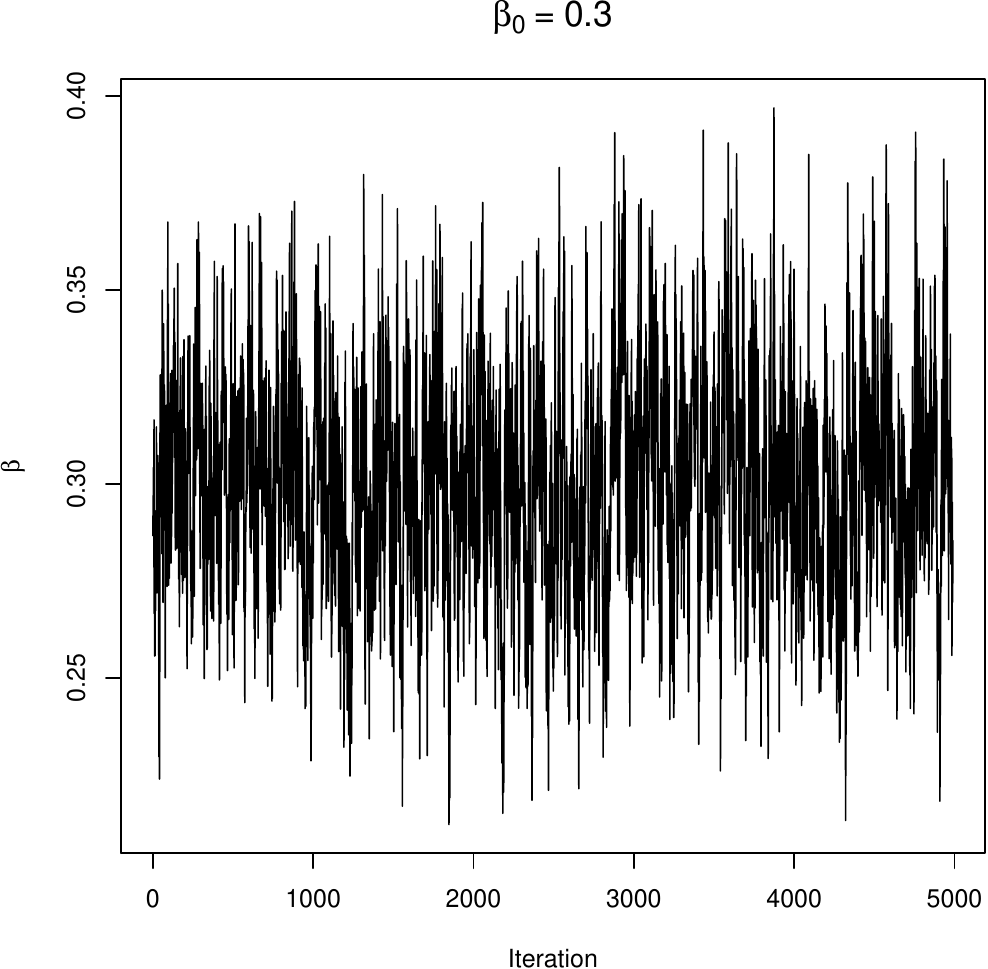}}\\
\vspace{2mm}
\subfigure []{ \label{fig:h3}
\includegraphics[width=5.5cm,height=4.5cm]{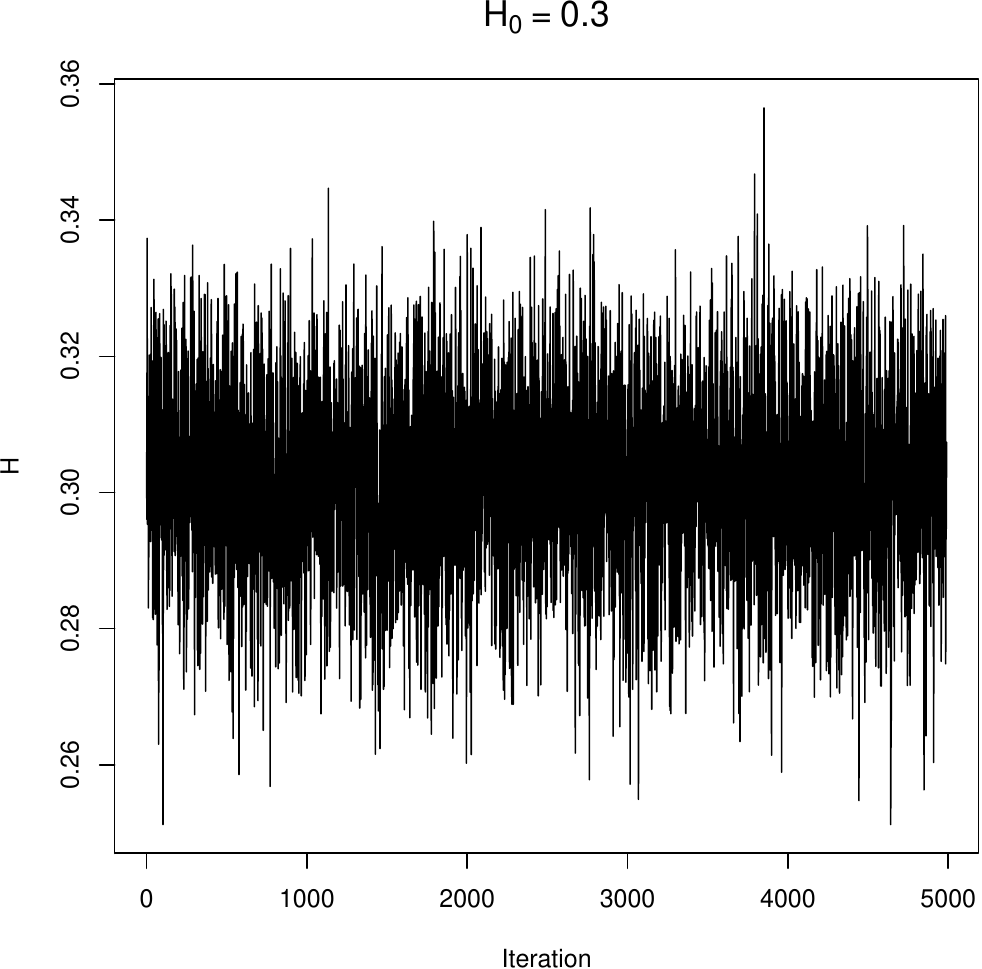}}
\hspace{2mm}
\subfigure []{ \label{fig:b3}
\includegraphics[width=5.5cm,height=4.5cm]{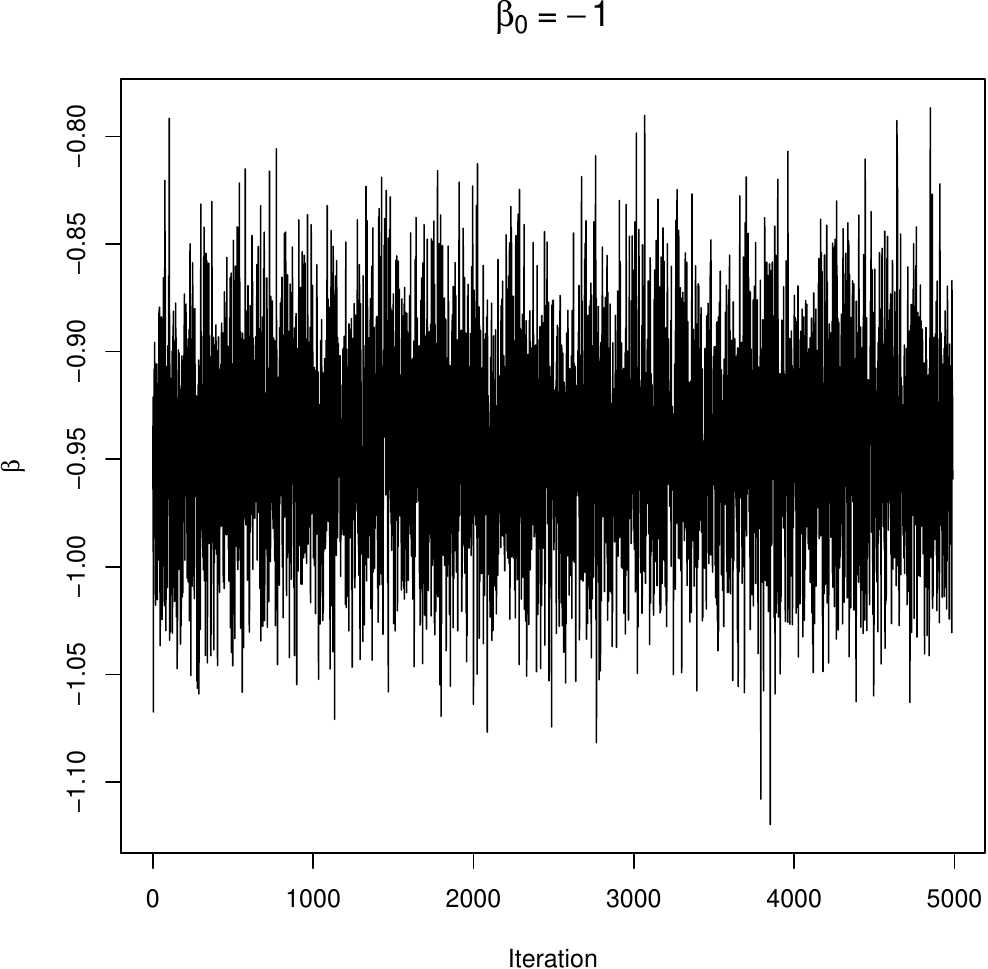}}\\
\caption{TMCMC trace plots of $H$ and $\beta$.}
\label{fig:traceplots1}
\end{figure}

\begin{figure}
\centering
\subfigure []{ \label{fig:h4}
\includegraphics[width=5.5cm,height=4.5cm]{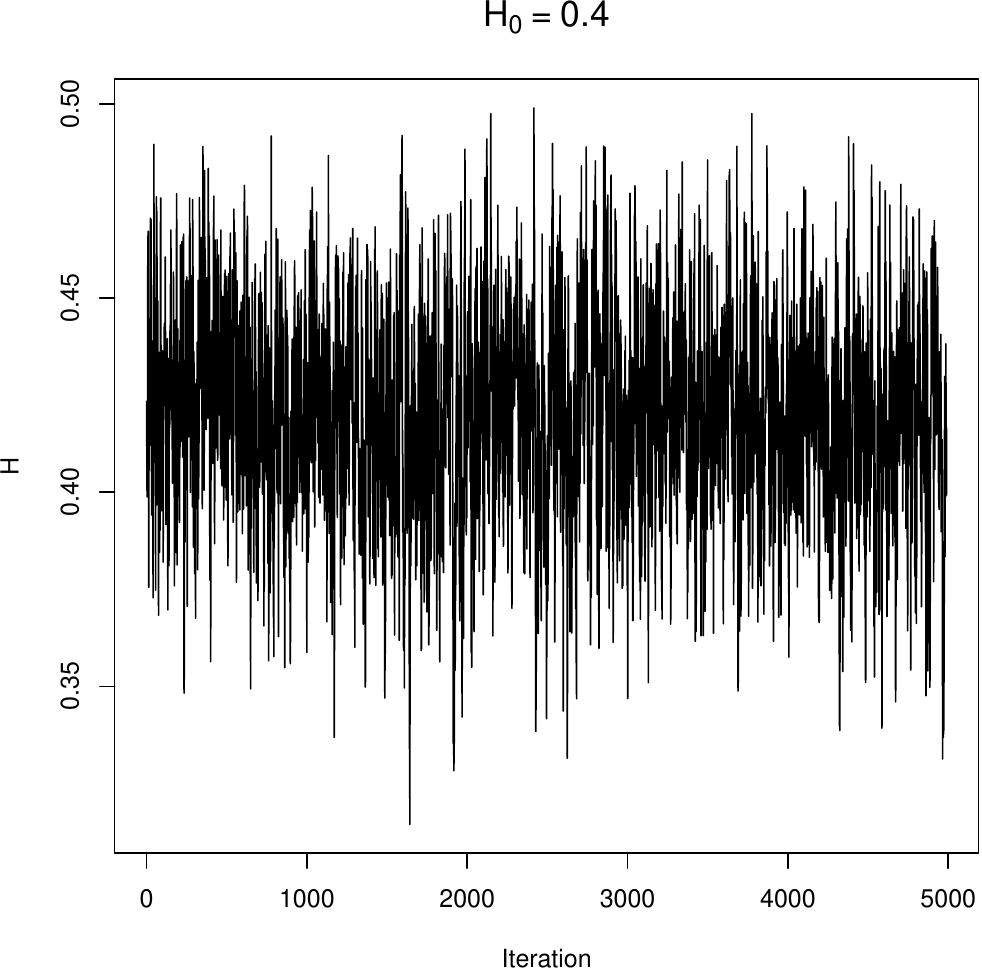}}
\hspace{2mm}
\subfigure []{ \label{fig:b4}
\includegraphics[width=5.5cm,height=4.5cm]{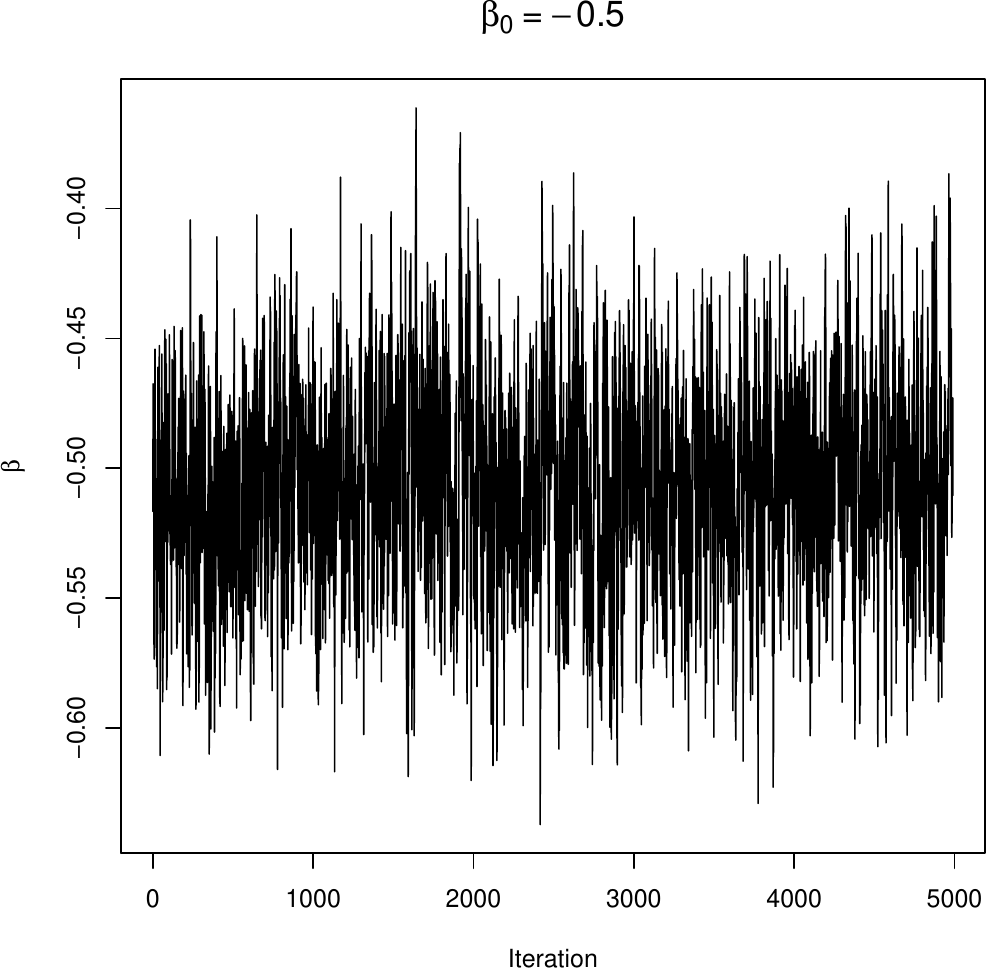}}\\
\vspace{2mm}
\subfigure []{ \label{fig:h5}
\includegraphics[width=5.5cm,height=4.5cm]{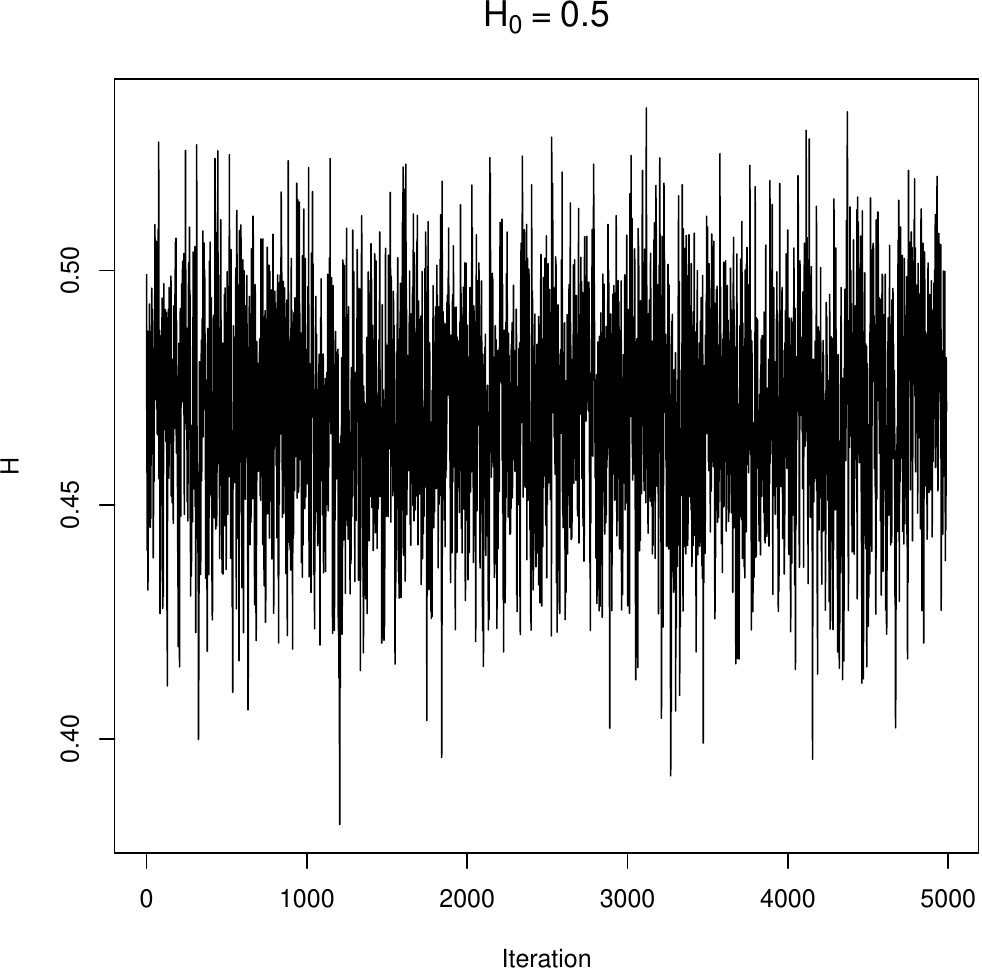}}
\hspace{2mm}
\subfigure []{ \label{fig:b5}
\includegraphics[width=5.5cm,height=4.5cm]{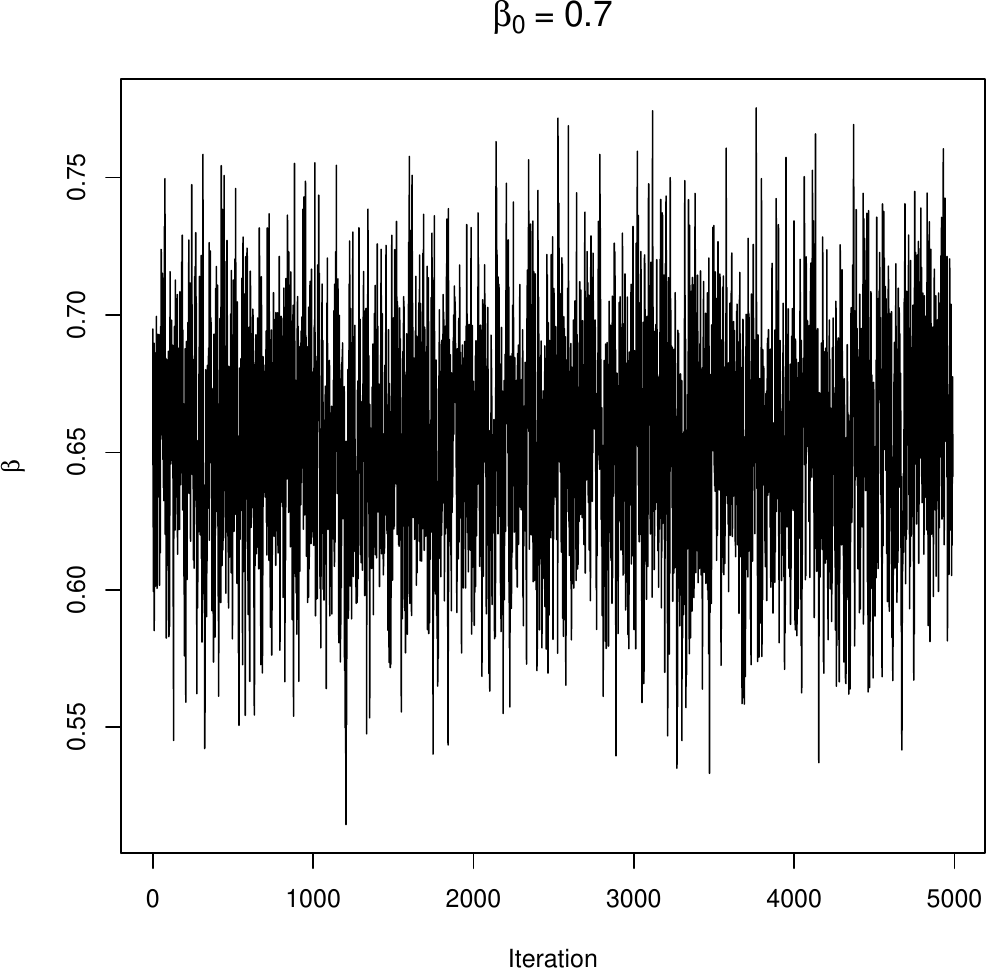}}\\
\vspace{2mm}
\subfigure []{ \label{fig:h6}
\includegraphics[width=5.5cm,height=4.5cm]{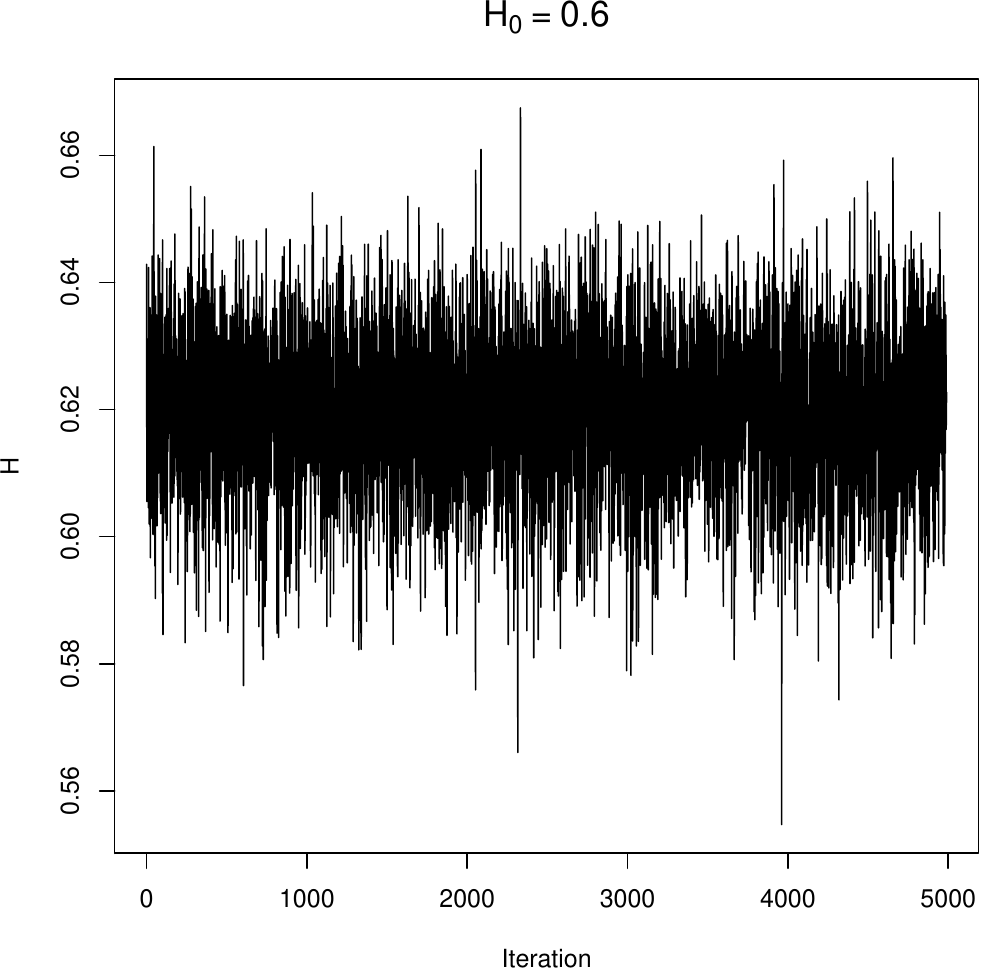}}
\hspace{2mm}
\subfigure []{ \label{fig:b6}
\includegraphics[width=5.5cm,height=4.5cm]{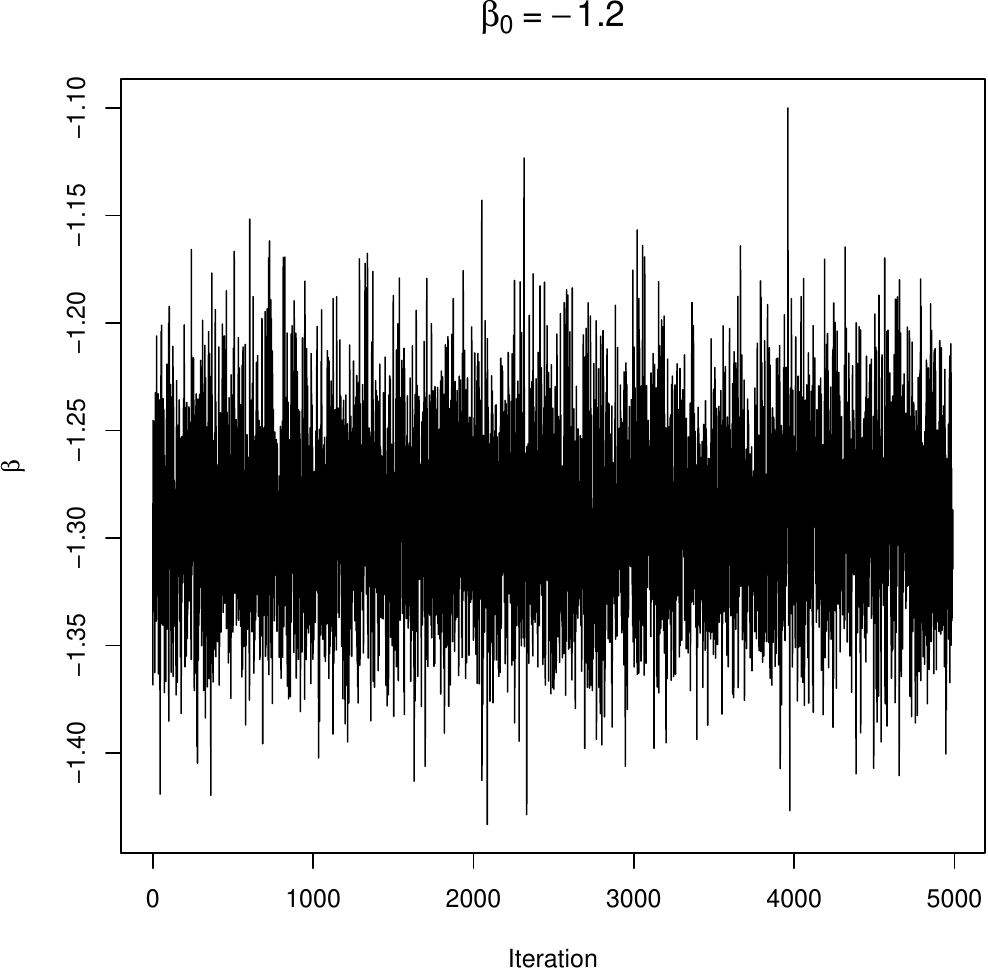}}\\
\caption{TMCMC trace plots of $H$ and $\beta$.}
\label{fig:traceplots2}
\end{figure}

\begin{figure}
\centering
\subfigure []{ \label{fig:h7}
\includegraphics[width=5.5cm,height=4.5cm]{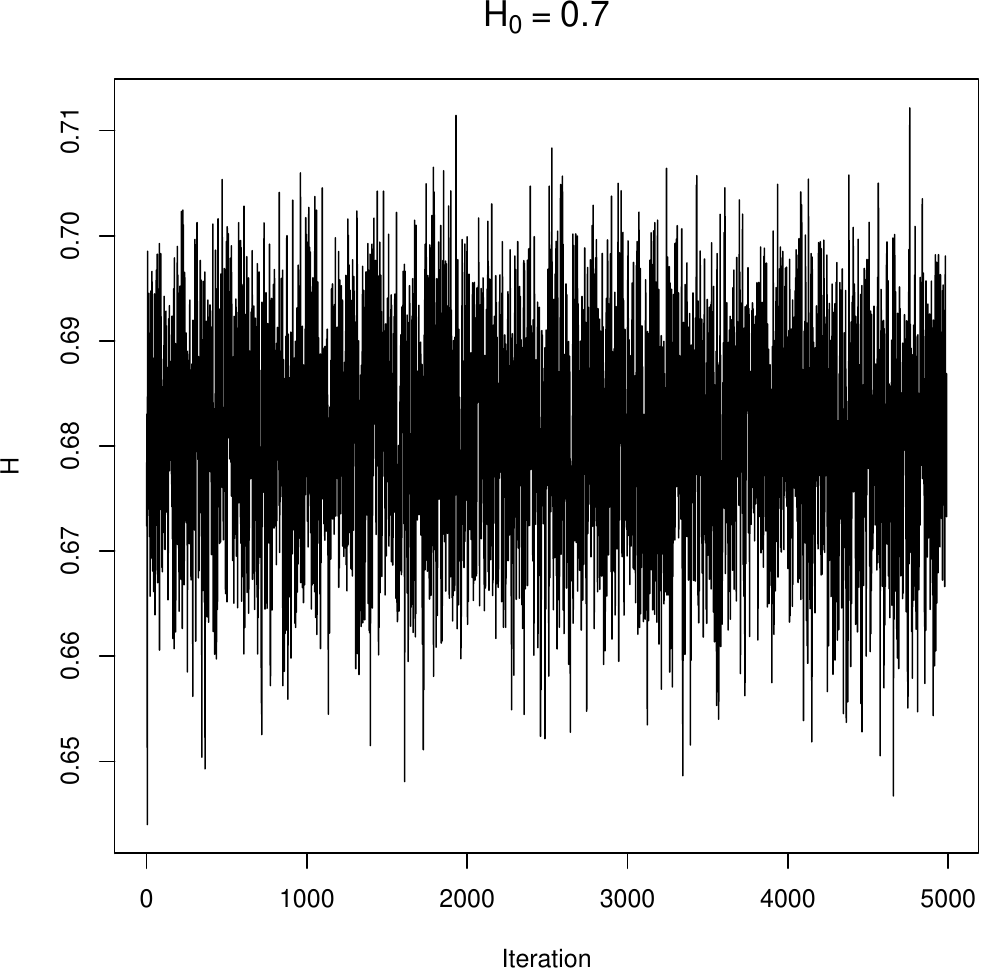}}
\hspace{2mm}
\subfigure []{ \label{fig:b7}
\includegraphics[width=5.5cm,height=4.5cm]{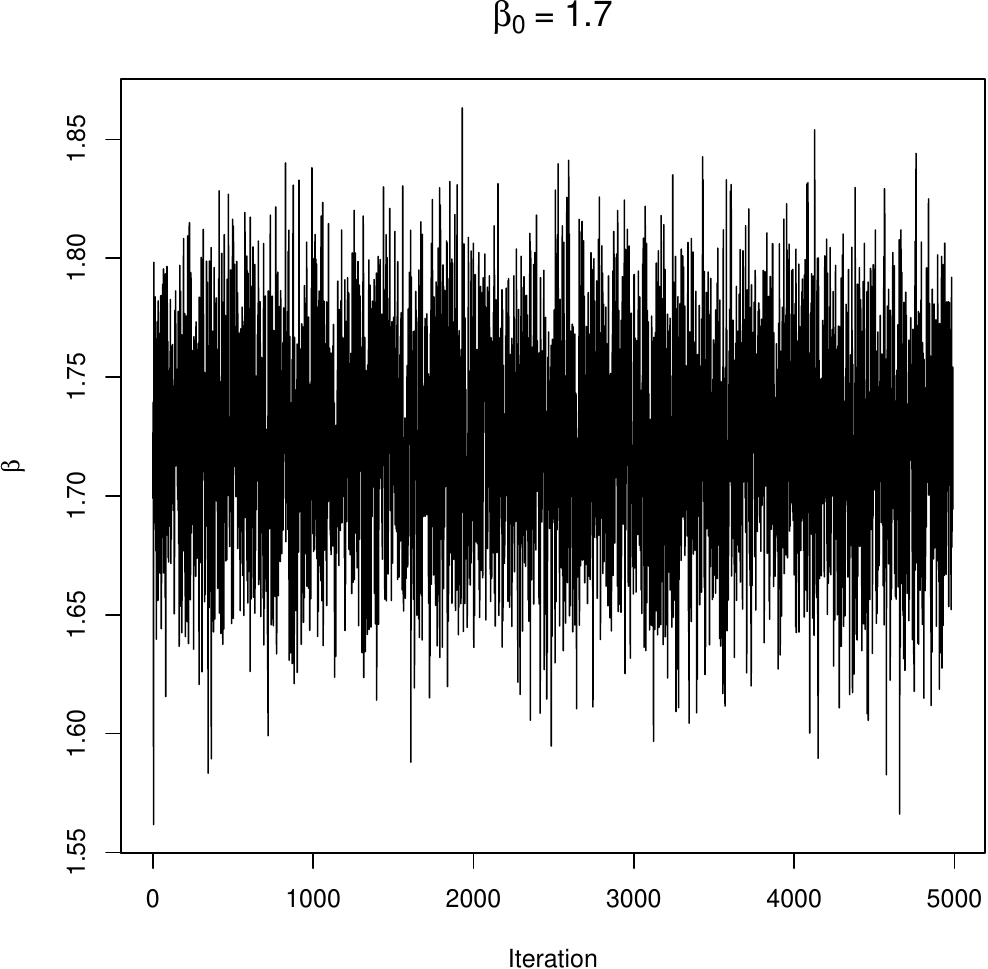}}\\
\vspace{2mm}
\subfigure []{ \label{fig:h8}
\includegraphics[width=5.5cm,height=4.5cm]{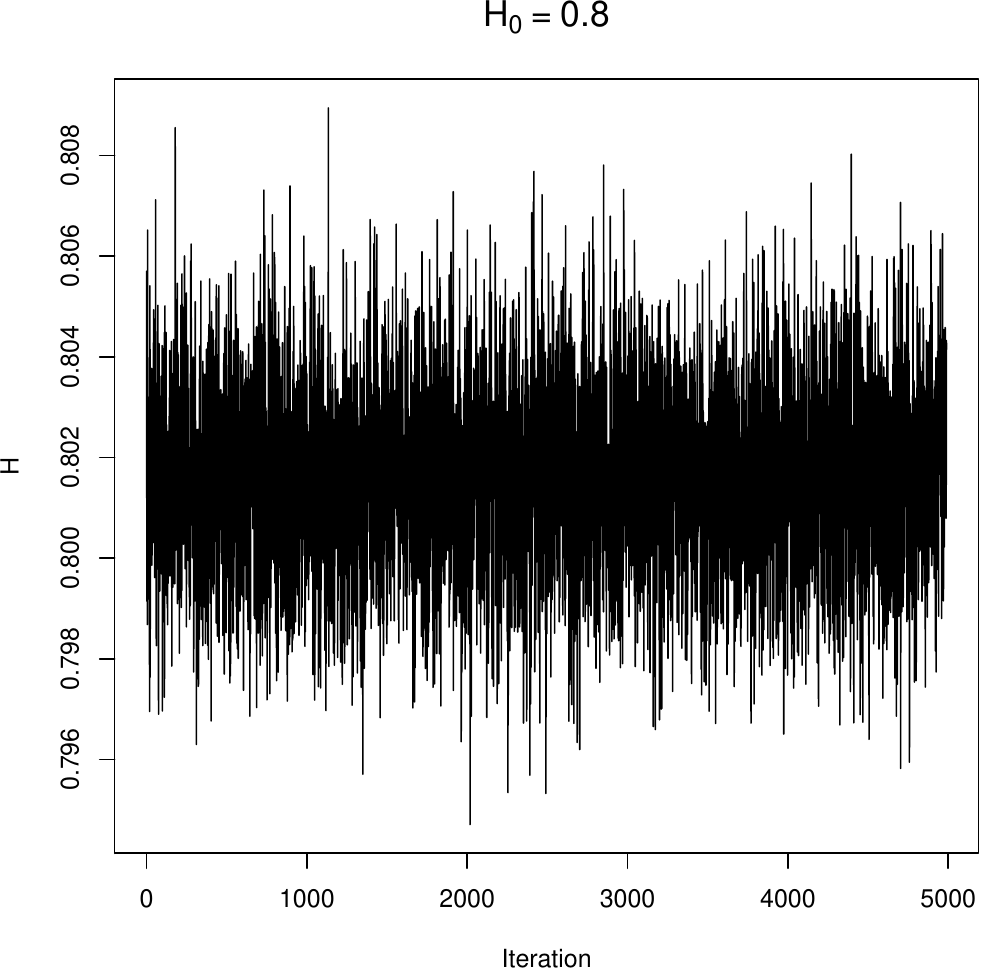}}
\hspace{2mm}
\subfigure []{ \label{fig:b8}
\includegraphics[width=5.5cm,height=4.5cm]{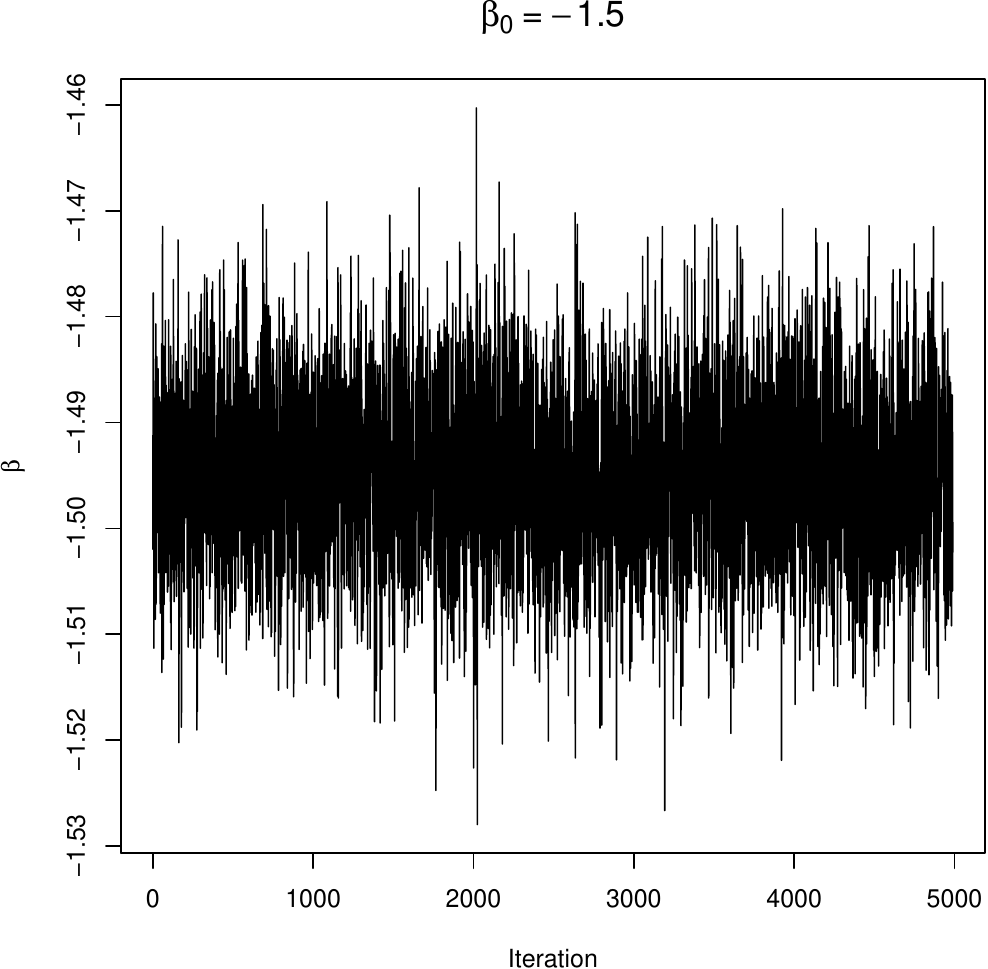}}\\
\vspace{2mm}
\subfigure []{ \label{fig:h9}
\includegraphics[width=5.5cm,height=4.5cm]{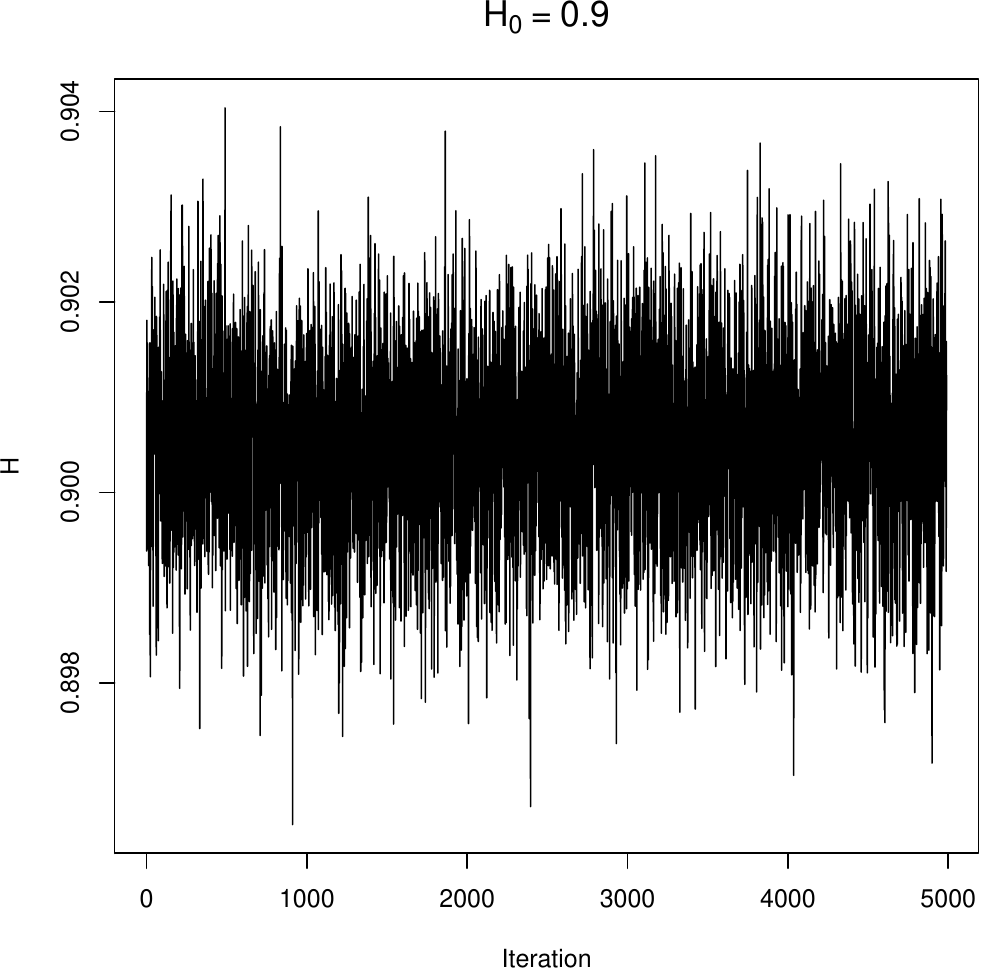}}
\hspace{2mm}
\subfigure []{ \label{fig:b9}
\includegraphics[width=5.5cm,height=4.5cm]{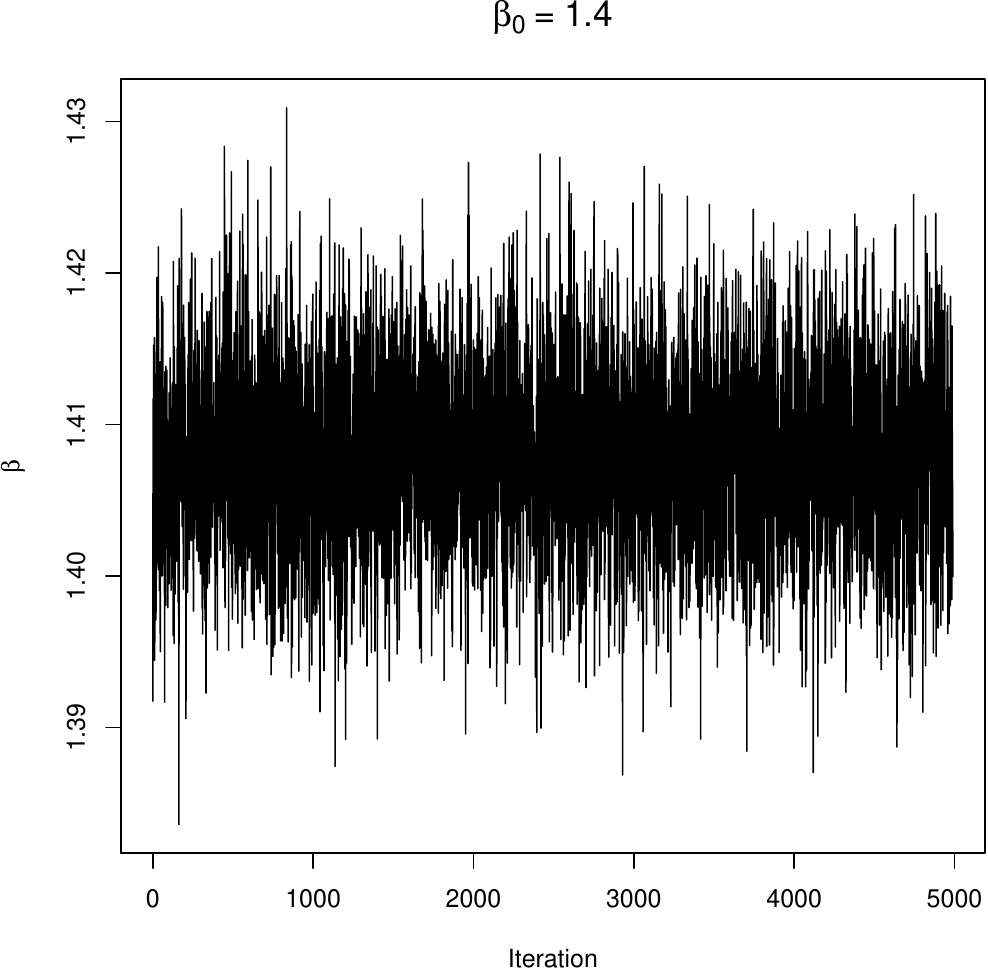}}\\
\caption{TMCMC trace plots of $H$ and $\beta$.}
\label{fig:traceplots3}
\end{figure}

In our implementation of Algorithm \ref{algo:tmcmc}, we choose $N=5\times 10^5$ for all the nine cases. Also, we set $a_1=a_2=0.05$ for all the cases except cases
$2$ and $7$, where we set $a_1=a_2=0.01$ and $a_1=a_2=1$, respectively. The initial values are set to the $MLE$s. 
We discard the first $1000$ TMCMC iterations as burn-in in all the cases and make our
Bayesian inference based on the remaining iterations. Our implementation of Algorithm \ref{algo:tmcmc} yielded excellent mixing properties
in all the cases, which is evident from the trace plots shown in Figures \ref{fig:traceplots1}, \ref{fig:traceplots2} and \ref{fig:traceplots3}  after 
discarding the first $1000$ iterations and subsequently plotting thinned samples (one in every 50 TMCMC realizations) to reduce file size.
Our implementation of Algorithm \ref{algo:tmcmc} in the $R$ software takes about $13$ minutes on our ordinary laptop to generate $5\times 10^5$ TMCMC realizations.

Table \ref{table:bayesian} summarizes the results of our Bayesian analyses for cases 1--9. 
\begin{table}[h]
\centering
	\caption{Results of our simulation study in the Bayesian setup. Here $\hat H^B_{T,n}$ and $\hat \beta^B_{T,n}$ denote the posterior expectations
	of $H$ and $\beta$, and 95\% $BCI$ stands for 95\% Bayesian credible intervals.}
\label{table:bayesian}
	\begin{tabular}{|c|c|c||c|c|c|}
\hline
		$H_0$ & $\hat H^B_{T,n}$ & 95\% $BCI$ & $\beta_0$ & $\hat \beta^B_{T,n}$ & 95\% $BCI$\\ 
\hline
	0.1  &  0.106  &  [0.097,0.116]  &  1.0  &  0.953  &  [0.870,1.036]\\
	0.2  &  0.212  &  [0.174,0.252]  &  0.3  &  0.300  &  [0.246,0.356]\\
	0.3  &  0.302  &  [0.274,0.329]  &  -1.0  &  -0.942  &  [-1.032,-0.851]\\
	0.4  &  0.419  &  [0.364,0.470]  &  -0.5  & -0.506  &  [-0.582,-0.428]\\   
	0.5  &  0.469  &  [0.426,0.509]  &  0.7   &  0.654  &  [0.579,0.727]\\ 
	0.6  &  0.620  &  [0.593,0.643]  &  -1.2  &  -1.290 &  [-1.370,-1.205]\\
	0.7  &  0.681  &  [0.661,0.700]  &  1.7   &  1.724  &  [1.640,1.806]\\
	0.8  &  0.802  &  [0.798,0.805]  &  -1.5  &  -1.495  &  [-1.511,-1.479]\\
	0.9  &  0.901  &  [0.899,0.902]  &  1.4   &  1.408    &  [1.396,1.420]\\
\hline
\end{tabular}
\end{table}
Note that in contrast to Table \ref{table:classical}, the true values
of $H$ and $\beta$ are contained in the 95\% Bayesian credible intervals ($BCI$s) in all the cases. Moreover, the lengths of such intervals are generally substantially shorter
than their classical counterpart, as shown in Table \ref{table:classical}. 
Indeed, the 95\% $BCI$ for $H$ are shorter than the corresponding classical 95\% $CI$ in all the cases, and for $\beta$, the 95\% $BCI$ are shorter than the corresponding 
95\% $CI$ in all the cases except only cases $6$ and $7$. The average lengths of the 95\% $BCI$ and 95\% $CI$ for $H$ are $0.049$ and $0.060$, respectively, and for $\beta$, 
the lengths are $0.127$ and $0.134$, respectively.

In other words, our simulation experiments vouch for the Bayesian framework for analysing
$SDE$s driven by fractional Brownian motion using the priors and the TMCMC algorithm that we propose.

\section{Application to real data}
\label{sec:realdata}
To deal with real data we collect the close price of a particular company from stock market data 
($460$ observations during the time range August $5$, 2013, to June $30$, 2015)  which is available on {\it www.nseindia.com}. 

A simple autocorrelation plot of the close price data revealed autocorrelations close to one, even for lag greater than $25$, suggesting long-range dependence. 
Thus, an $SDE$ driven by a fractional Brownian motion seems to be appropriate here. 
Moreover, for simplicity we fix the choice of the ratio $B(t,\tilde x_t)/C(t,\tilde x_t)$ as a constant; in fact, we set $B(t,\tilde x_t)=C(t,\tilde x_t)=\tilde x^{\nu}_t$, where $\nu>0$ is to be
determined appropriately, following the discussion below. 

In other words, we propose the following form of such $SDE$:
\begin{equation*}
	dX_t=\beta \tilde x_t^{\nu}dt + \gamma \tilde x_t^{\nu}dW_t^H,
\end{equation*}
where $\nu>0$ is to be chosen such that both classical and Bayesian analyses of the data using the above $SDE$ supports values of $H$ close to one. As it turned out,
$\nu=3/2$ is an appropriate choice in this regard. In other words, we fit the following $SDE$ to the close price data with respect to both classical and Bayesian setups:
\begin{equation}
d X_t=\beta \tilde x_t^{3/2}dt + \gamma \tilde x_t^{3/2}dW_t^ H
	\label{eq:sde_real}
\end{equation}
Here $\beta$, $H$ and $\gamma$ are unknown parameters which need to be learned about. However, as is the case in most statistical applications of $SDE$,
we choose to assume $\gamma$ as given. In both classical and Bayesian setups 
we choose $\gamma$ to minimize the average lengths of the 95\% confidence and credible intervals, respectively, associated with the held-out observations, subject
to inclusion of all the held-out data points.
In our application we consider the first $400$ data points as given and predict the next $60$ data points in both classical and Bayesian setups. 
The values of $\gamma$ thus obtained would be considered as given for prediction of actual future observations. 

\subsection{Real data analysis in the classical setup}
\label{subsec:classical_real}
Starting from arbitrary initial values of $\btheta = (\beta,H)$ we performed simulated annealing aided with a TMCMC algorithm 
that has provision for further enhancement of mixing properties to obtain the $MLE$ of $\btheta$. Our algorithm is provided  as Algorithm \ref{algo:sima_real}. 
\begin{algo}
\topline{Simulated annealing aided by TMCMC with enhanced mixing properties}
\label{algo:sima_real}
\botline
\begin{itemize}
	\item[(1)] Begin with an initial value $\btheta^0=\left(\beta^{(0)},H^{(0)}\right)$.
\item[(2)] For $t=1,\ldots,N$, do the following:
	\begin{enumerate}
		\item[(i)] Generate $\varepsilon_1\sim N(0,1)$, $\varepsilon_2\sim U(-1,1)$, $b_j\stackrel{iid}{\sim}U(\{-1,1\})$ for $j=1,2$, $U_1\sim U(0,1)$. If $U_1<1/2$, create $\tilde\beta=\beta^{(t-1)}+b_1a_1|\varepsilon_1|$ and $\tilde H=H^{(t-1)}+b_2a_2|\varepsilon_1|$. Else, set $\tilde\beta=\beta^{(t-1)}{\varepsilon_2}^{b_1}$ and $\tilde H=H^{(t-1)}{\varepsilon_2}^{b_2}$.  
			Let $\tilde\btheta=\left(\tilde\beta,\tilde H\right)$.
		\item[(ii)] Evaluate 
			$$\alpha_1=\min\left\{1,\exp\left\{\frac{\ell_{T,n}\left(\tilde\btheta\right)-\ell_{T,n}\left(\btheta^{(t-1)}\right)+\log (J)}{\tau_t}\right\}
			I_{(0,1)}(\tilde H)\right\},$$
			where $\tau_t=\frac{1}{\log(t)}$ when $t\geq 100$ and $t^{-1}$ otherwise and $J$ is the Jacobian of corresponding transformation. ($J = 1$ in additive case and $J={|\varepsilon_2|}^{\sum b_j}$ in multiplicative case)
		\item[(iii)] Set $\btheta^{(t)}=\tilde\btheta$ with probability $\alpha_1$, else set $\btheta^{(t)}=\btheta^{(t-1)}$.

   To provide further mixing,
                  \item[(iv)] Generate $\varepsilon_3\sim N(0,1)$,  $U_2\sim U(0,1)$. If $U_2<1/2$  ,

    Further generate $u_1\sim U(0,1)$. If $u_1<1/2$, create $\beta^*=\beta^{(t)}+a_1|\varepsilon_3|$ and $H^*=H^{(t)}+a_2|\varepsilon_3|$. 
			Else, set $\beta^*=\beta^{(t)}-a_1|\varepsilon_3|$ and $H^*=H^{(t)}-a_2|\varepsilon_3|$. 

If $U_2\geq 1/2$,
  
Generate $\varepsilon_4\sim U(-1,1)$, $u_2\sim U(0,1)$. If $u_2<1/2$, create $\beta^*=\beta^{(t)}\varepsilon_4$ and $H^*=H^{(t)}\varepsilon_4$. Else, set $\beta^*=\beta^{(t)}/\varepsilon_4$ and $H^*=H^{(t)}/\varepsilon_4$.
Let $\btheta^*=\left(\beta^*,H^*\right)$.
		\item[(v)] Evaluate 
			$$\alpha_2=\min\left\{1,\exp\left\{\frac{\ell_{T,n}\left(\tilde\btheta\right)-\ell_{T,n}\left(\btheta^{(t-1)}\right)+\log (J)}{\tau_t}\right\}
			I_{(0,1)}(\tilde H)\right\},$$
			where $\tau_t=\frac{1}{\log(t)}$ when $t\geq 100$ and $t^{-1}$ otherwise and $J$ is the Jacobian of corresponding transformation. (In the case of $U_2<1/2$, $J = 1$  and otherwise, $J={|\varepsilon_2|}^2$ or $J={|\varepsilon_2|}^{-2}$  correspondingly)
		\item[(vi)] Accept $\btheta^*$ with probability $\alpha_2$ or remain at $\btheta^{(t)}$. To avoid notational complexity whether it is $\btheta^*$ or $\btheta^{(t)}$ we write final value of each iteration as $\btheta^{(t)}$.
	\end{enumerate}
\item[(3)] From the sample $\left\{\btheta^{(j)};j=0,1,\ldots,N\right\}$, set $\hat\btheta_{T,n}=\btheta^{(j^*)}$ where 
	$\btheta^{(j^*)}=\underset{j=0,1,\ldots,N}{\arg\max}~\ell_{T,n}\left(\btheta^{(j)}\right)$.
\end{itemize}
\botline
\end{algo}
In our application, initial values are arbitrarily chosen and we set $N=10^4$ and $a_1=a_2=0.0005$. With the obtained $MLE$ given a specific value of $\gamma$ we 
predicted the close price of next $60$ days for that company in two ways. In the first way, at every step each data point is predicted with the 
consideration that data values at the previous time point is given beforehand, that is, one-step-ahead forecast is obtained. 
In the other way, starting from the 400th data point, predictions for the next 60 data points are obtained, which is the $60$-step-ahead forecast. 
To obtain desired confidence intervals for the forecasts, we generated $10,000$ realizations of size $60$ from the $SDE$ (\ref{eq:sde_real}) given the first $400$
data points.
We repeated the calculations for different values of $\gamma$, from which we obtained the 
final value of $\gamma$ as the one which yielded minimum average length of the 95\% confidence interval.
Associated with the final value of $\gamma$, we also report the number of held-out data points that fail to lie within the 95\% and the 50\% confidence intervals.


We repeated the entire procedure of classical analysis for the case where the $SDE$ is assumed to be driven by Brownian motion instead of fractional Brownian motion. 
In the latter case, as $H=0.5$ corresponds to Brownian motion, the only unknown parameter is $\beta$ and $\gamma$ is fixed in the same way as in the fractional set up. 
This analysis facilitates comparison of fractional $SDE$ model with the more usual $SDE$ models driven by Brownian motion.

For the classical analysis $\gamma$ is found to be $165$ and $0.4$ for our $SDE$s driven by fractional and standard Brownian motion, respectively. 
The results are tabulated in Tables \ref{table:classical_real_fBm} and \ref{table:classical_real_Bm} showing the $MLE$s of $H$ and $\beta$, the values
of $\gamma$ for one step ahead and $60$ steps ahead forecasts, the number of held-out data points outside the corresponding 50\% and 95\% confidence intervals $(CI$)
and the average length of the 95\% CIs of the held-out data points.
Observe that compared to fractional $SDE$, the average lengths of the 95\% CI are small for standard $SDE$, but a significantly larger number of held-out observations
are excluded from the 50\% $CI$s of the standard $SDE$. This indicates that compared to standard $SDE$, most of the held-out observations are concentrated
in the high-density regions of the fractional $SDE$. 
Hence, the fractional $SDE$ clearly outperforms the standard $SDE$.

\begin{figure}
\centering
\subfigure []{ \label{fig:c1_165}
\includegraphics[width=7.5cm,height=6.5cm]{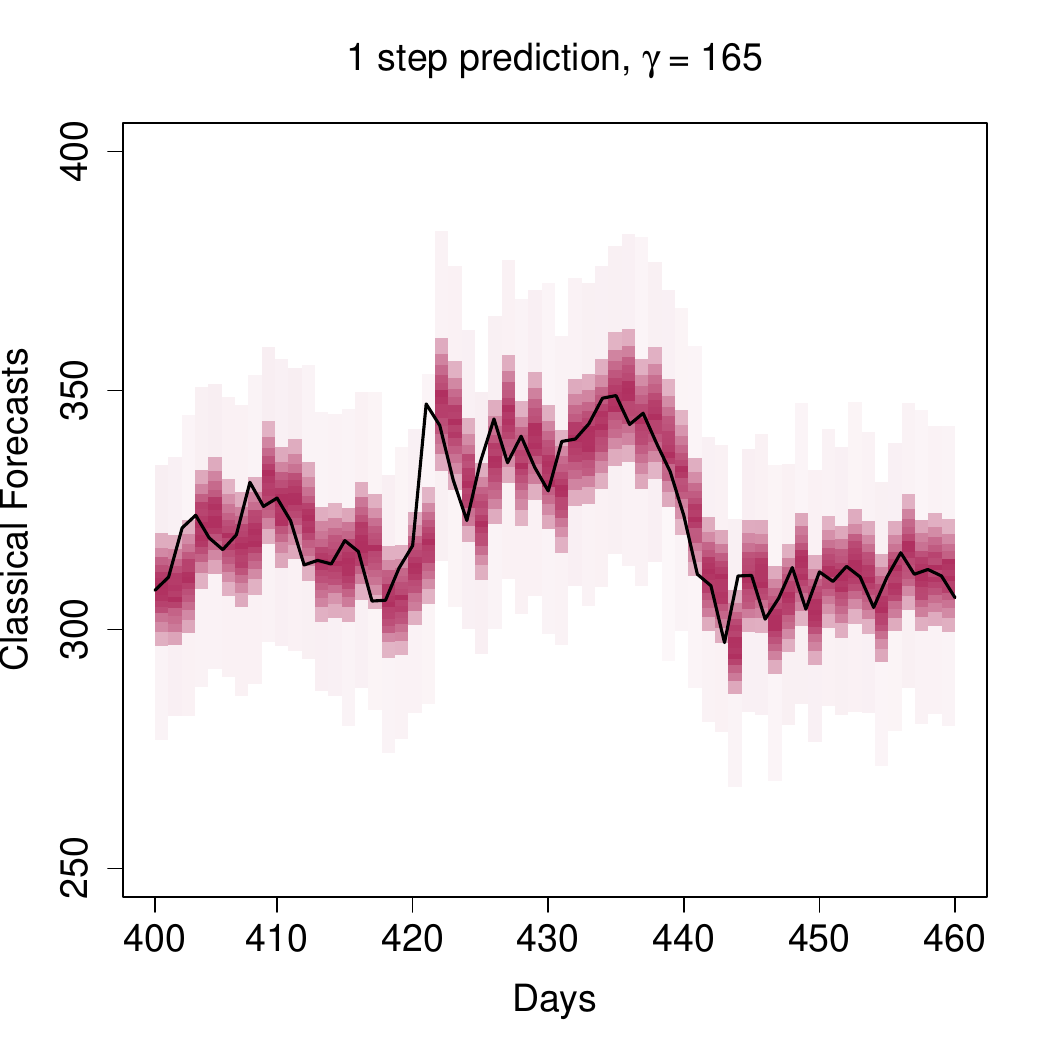}}
\hspace{2mm}
\subfigure []{ \label{fig:c60_165}
\includegraphics[width=7.5cm,height=6.5cm]{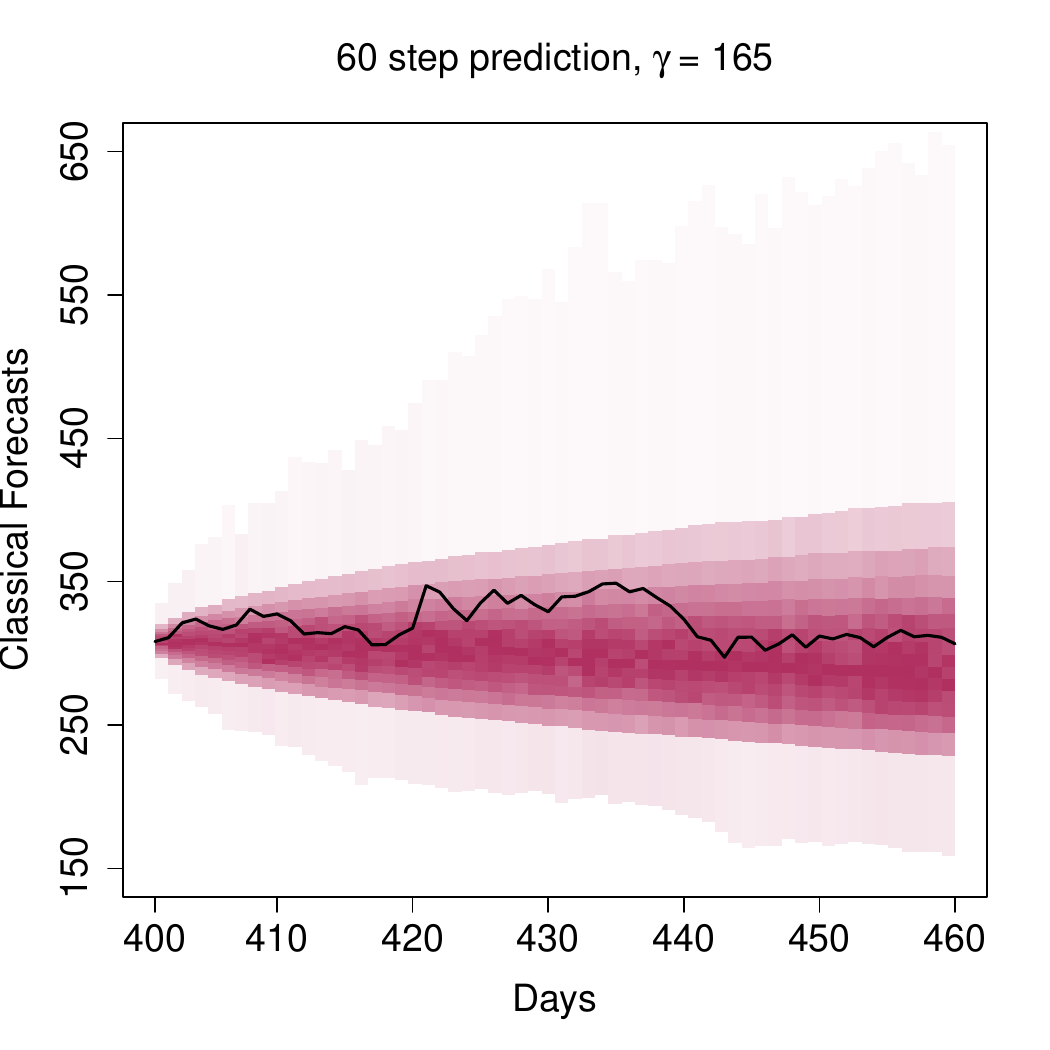}}
\caption{One step ahead and 60 steps ahead forecasts in the classical setup corresponding to fractional $SDE$.}
\label{fig:1_60_classical_fBM}
\end{figure}

\begin{figure}
\centering
\subfigure []{ \label{fig:c1_bm0.4}
\includegraphics[width=7.5cm,height=6.5cm]{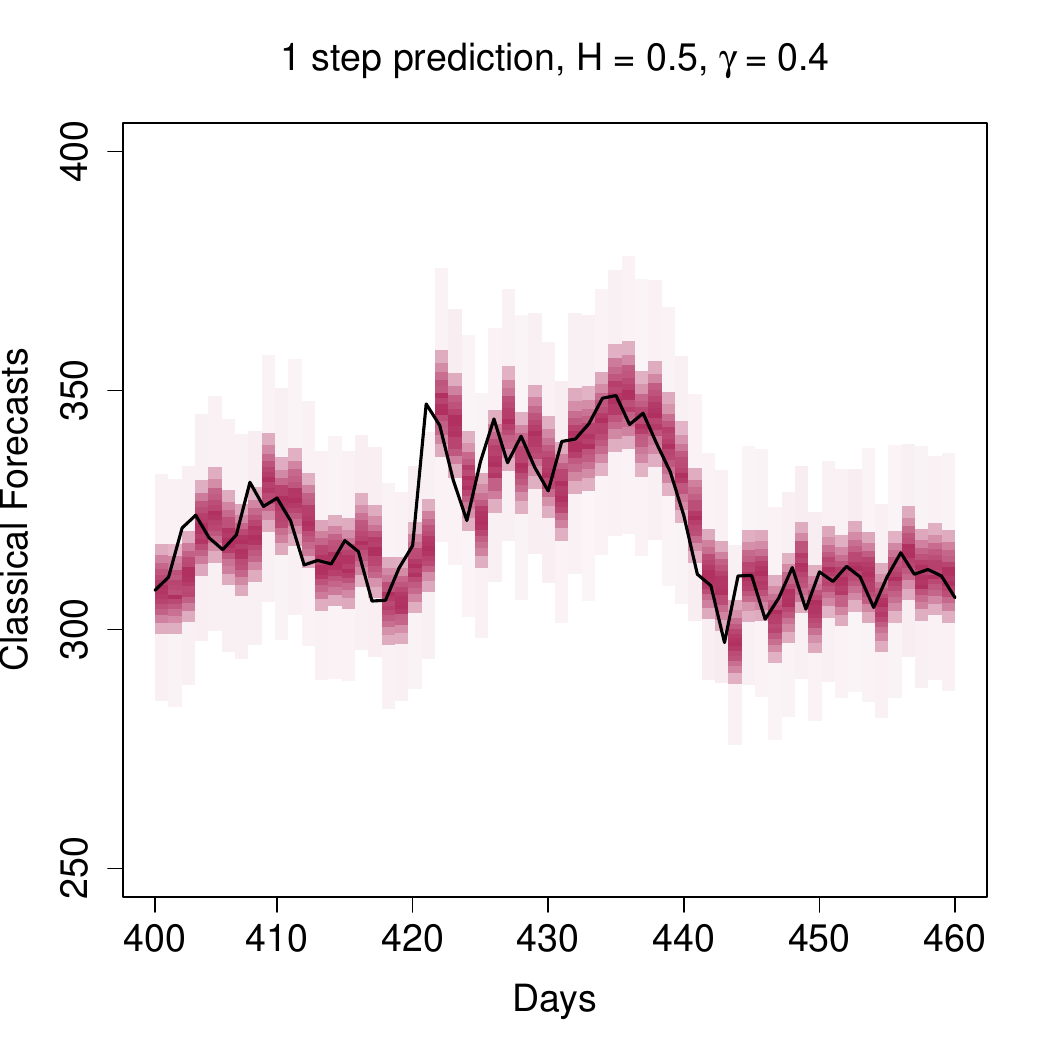}}
\hspace{2mm}
\subfigure []{ \label{fig:c60_bm0.4}
\includegraphics[width=7.5cm,height=6.5cm]{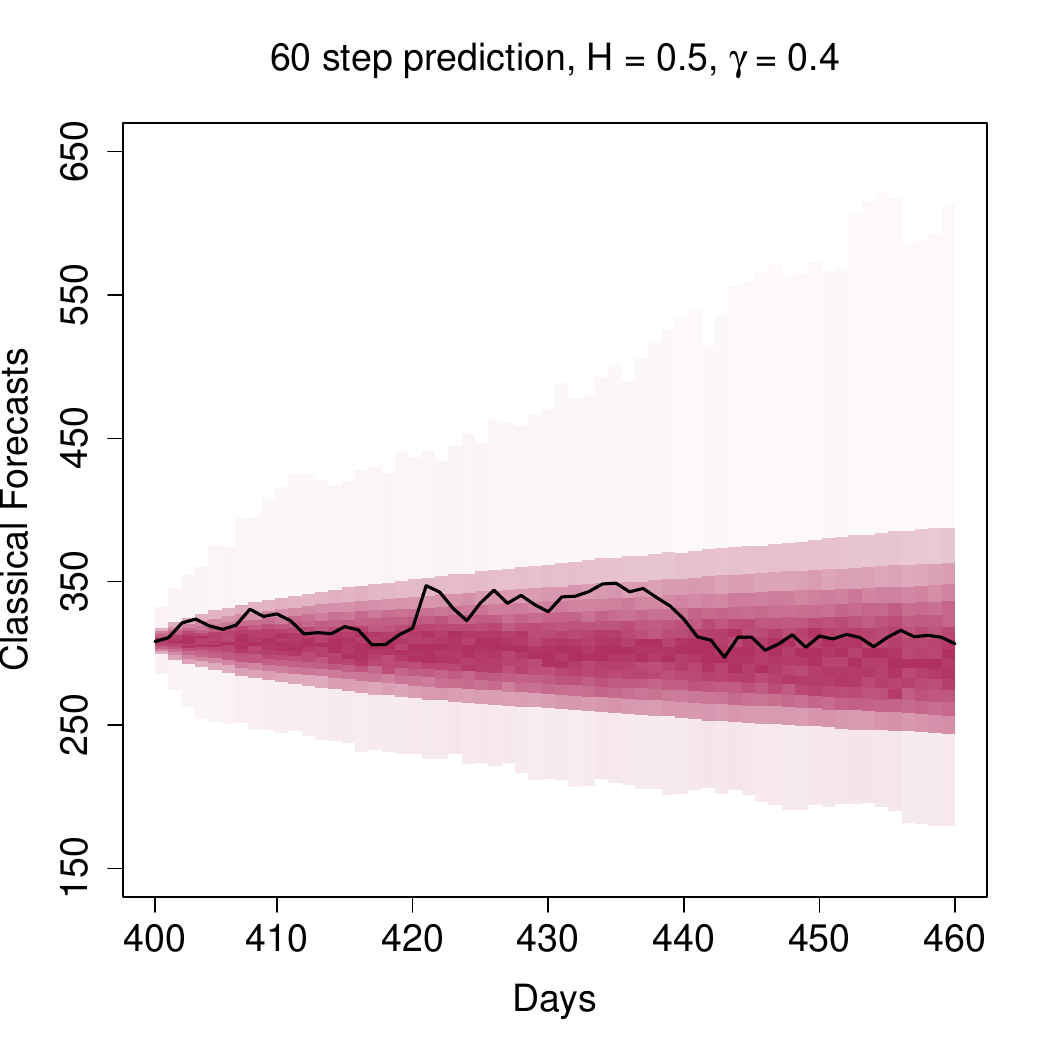}}\\
\caption{One step ahead and 60 steps ahead forecasts in the classical setup corresponding to standard $SDE$.}
\label{fig:1_60_classical_BM}
\end{figure}

The densities associated with fractional and standard $SDE$s are displayed in Figures 
\ref{fig:1_60_classical_fBM} and \ref{fig:1_60_classical_BM}, respectively, where higher densities are represented by progressively intense colours;
the thick, black line stands for the held-out observed data. These figures are the detailed versions of Tables 
\ref{table:classical_real_fBm} and \ref{table:classical_real_Bm}. 
\begin{table}[h]
\centering
\caption{Results of classical analysis of real data modeled by fractional $SDE$.}
\label{table:classical_real_fBm}
	\begin{tabular}{|c|c|c|c|c|c|c|}
\hline
$\hat H_{T,n}$ & $\hat \beta_{T,n}$ & $\gamma$ & step & \# outside 50\% $CI$ & \# outside 95\% $CI$ & Average length of 95\% $CI$\\ 
\hline
	0.995  &  -0.002578  &  165  &  1  &  24  & 1 & 31.11014\\
	0.995  &  -0.002578  &  165 &  60  &  17  &  0 & 155.5142\\
\hline
\end{tabular}
\end{table}

\begin{table}[h]
\centering
\caption{Results of classical analysis of real data modeled by standard $SDE$.}
\label{table:classical_real_Bm}
	\begin{tabular}{|c|c|c|c|c|c|c|}
\hline
$H$ & $\hat \beta_{T,n}$ & $\gamma$ & step & \# outside 50\% $CI$ & \# outside 95\% $CI$ & Average length of 95\% $CI$\\ 
\hline
	0.5  &  0.00021  &  0.4  &  1  &  38  & 2 & 25.20115\\
	0.5  &  0.00021  &  0.4 &  60  &  24  &  0 & 122.7201\\
\hline
\end{tabular}
\end{table}

\subsection{Real data analysis in the Bayesian setup}
\label{subsec:bayesian_real}

In the Bayesian setup, we first reparameterize $H$ as $H=\exp(\nu)/(1+\exp(\nu))$ as in
Section \ref{subsec:bayesian_simstudy} and assume {\it a priori} that $\nu\sim N\left(\nu_{T,n},\sigma^2_{\nu}\right)$ and 
$\beta\sim N\left(\mu_{\beta},\sigma^2_{\beta}\right)$, where $\mu_{\beta}=0.001, \sigma_{\nu}=\sigma_{\beta}=0.1$. 
As in Section \ref{subsec:bayesian_simstudy}, $\nu_{T,n}$ is set as $\log\left(\frac{\hat H_{T,n}}{1-\hat H_{T,n}}\right)$.

For Bayesian inference, a TMCMC algorithm with enhanced mixing properties is considered similar to Algorithm \ref{algo:sima_real}
where $\nu$ is updated instead of $H$ and the acceptance probability becomes \\
$\alpha=\min\left\{1,\exp\left\{\ell_{T,n}\left(\tilde\btheta\right)-\ell_{T,n}\left(\btheta^{(t-1)}\right)+\log (J)\right\}\right\}$. 
In our application, we set $N=10^4$ and $a_1=a_2=0.05$. 

Posterior predictive distributions of the $60$ held-out observations are obtained in two ways similar to that of the classical set up.
In each case, whether it is one step ahead or 60 steps ahead predictions, 95\% and 50\% Bayesian credible intervals are obtained and the number of data points lying outside 
corresponding regions are also noted. Again, $\gamma$ is treated as a fixed parameter whose value is obtained by minimizing the average length of the
95\% credible intervals subject to their inclusion of all the respective held-out observations. Once again we remind the reader that this fixed value of $\gamma$ would
be used even for actual forecasting future observations.

The whole analysis is carried out for the standard $SDE$ driven by Wiener process, that is, with $H=0.5$ 
and in this case we are interested in the posterior of $\beta$ only as $\gamma$ is fixed in the way as explained above.

In this Bayesian set up, we obtain $\gamma=158$ and $\gamma=1$ for fractional and standard $SDE$, respectively, corresponding to
which Tables \ref{table:bayesian_real_fBm} and \ref{table:bayesian_real_Bm} provide the results of our Bayesian analyses. 
\begin{table}[h]
\centering
\caption{Results of Bayesian analysis of real data modeled by fractional $SDE$.}
\label{table:bayesian_real_fBm}
	\begin{tabular}{|c|c|c|c|c|c|c|}
\hline
$ H^*$ & $\beta^*$ & $\gamma$ & step & \# outside 50\% $BCI$ & \# outside 95\% $BCI$ & Average length of 95\% $BCI$\\ 
\hline
	0.991  &  0.00121  &  158  &  1  &  25  &  1 & 31.472\\
	0.991  &  0.00121  &  158  &  60  & 8  &  0 & 201.4972\\   
\hline
\end{tabular}
\end{table}
\begin{table}[h]
\centering
\caption{Results of Bayesian analysis of real data modeled by standard $SDE$.}
\label{table:bayesian_real_Bm}
	\begin{tabular}{|c|c|c|c|c|c|c|}
\hline
$H$ & $\beta^*$ & $\gamma$ & step & \# outside 50\% $BCI$ & \# outside 95\% $BCI$ & Average length of 95\% $BCI$\\ 
\hline
	0.5  &  0.00035  &  1  &  1  &  6  &  0 & 62.98049\\
	0.5  &  0.00035 &  1  &  60  & 0 &  0 & 334.7808\\   
\hline
\end{tabular}
\end{table}
In the tables $H^*$ and $\beta^*$ stand for the respective posterior medians and Bayesian credible interval is denoted by $BCI$.
Note that the average length of 95\% $BCI$ for standard $SDE$ is much larger compared to that for fractional $SDE$ and hence it is no surprise
that even those for 50\% $BCI$s would be much larger for the former and hence excludes less number of observations. Indeed, the average lengths for
standard $SDE$ are much larger even compared to both the classical analyses reported in Tables \ref{table:classical_real_fBm} and \ref{table:classical_real_Bm},
and thus it is clear that the standard $SDE$ fails to outperform the fractional $SDE$. It is also clear, in terms of adequate average lengths of the 
95\% $BCI$ and the number of held-out observations falling in the high-density regions, that our Bayesian fractional $SDE$ is the best performer.

\begin{figure}
\centering
%
\subfigure []{ \label{fig:p1_158}
\includegraphics[width=7.5cm,height=6.5cm]{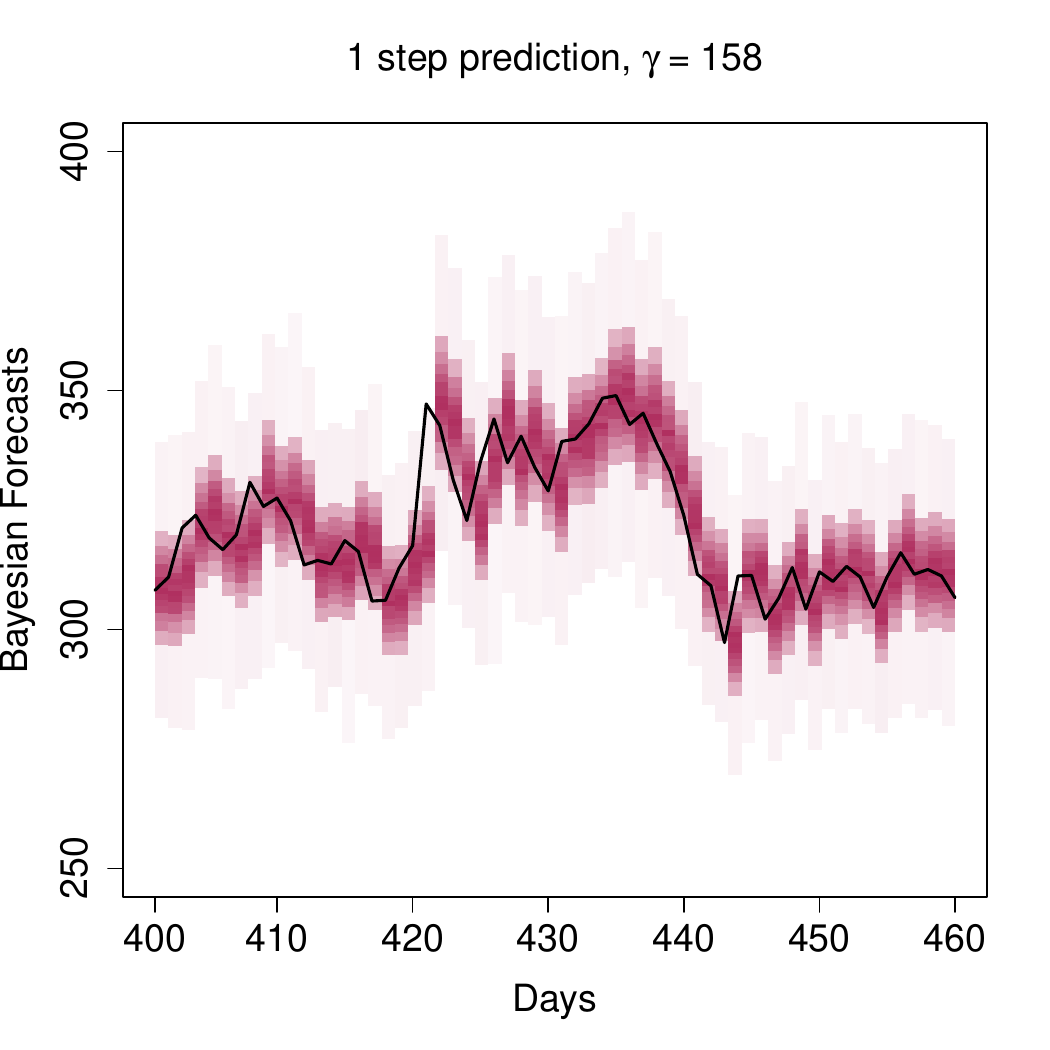}}
\hspace{2mm}
\subfigure []{ \label{fig:p60_158}
\includegraphics[width=7.5cm,height=6.5cm]{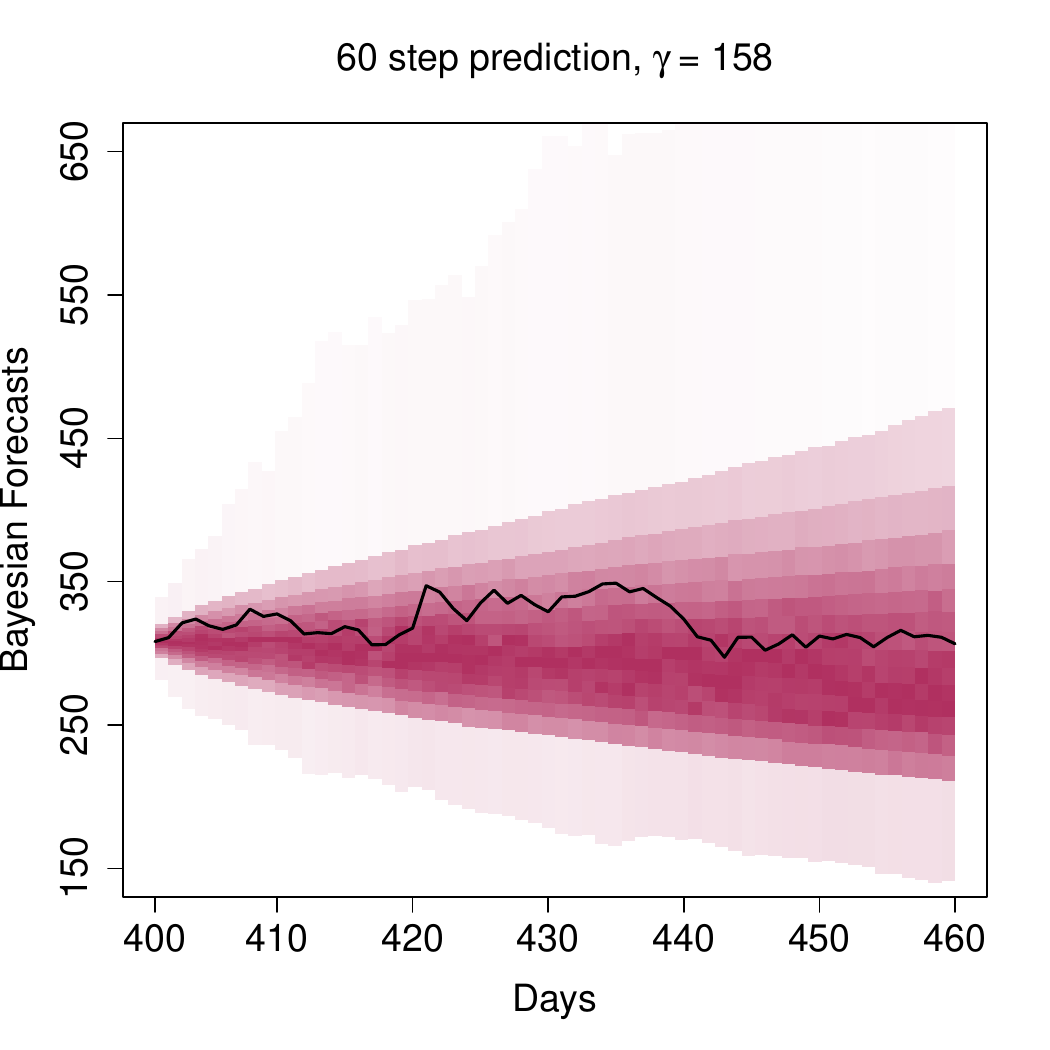}}\\
\caption{One and 60 steps ahead forecasts in the Bayesian setup corresponding to fractional $SDE$.}
\label{fig:1_60_Bayesian_fBM}
\end{figure}

\begin{figure}
\centering
\subfigure []{ \label{fig:p1_bm1}
\includegraphics[width=7.5cm,height=6.5cm]{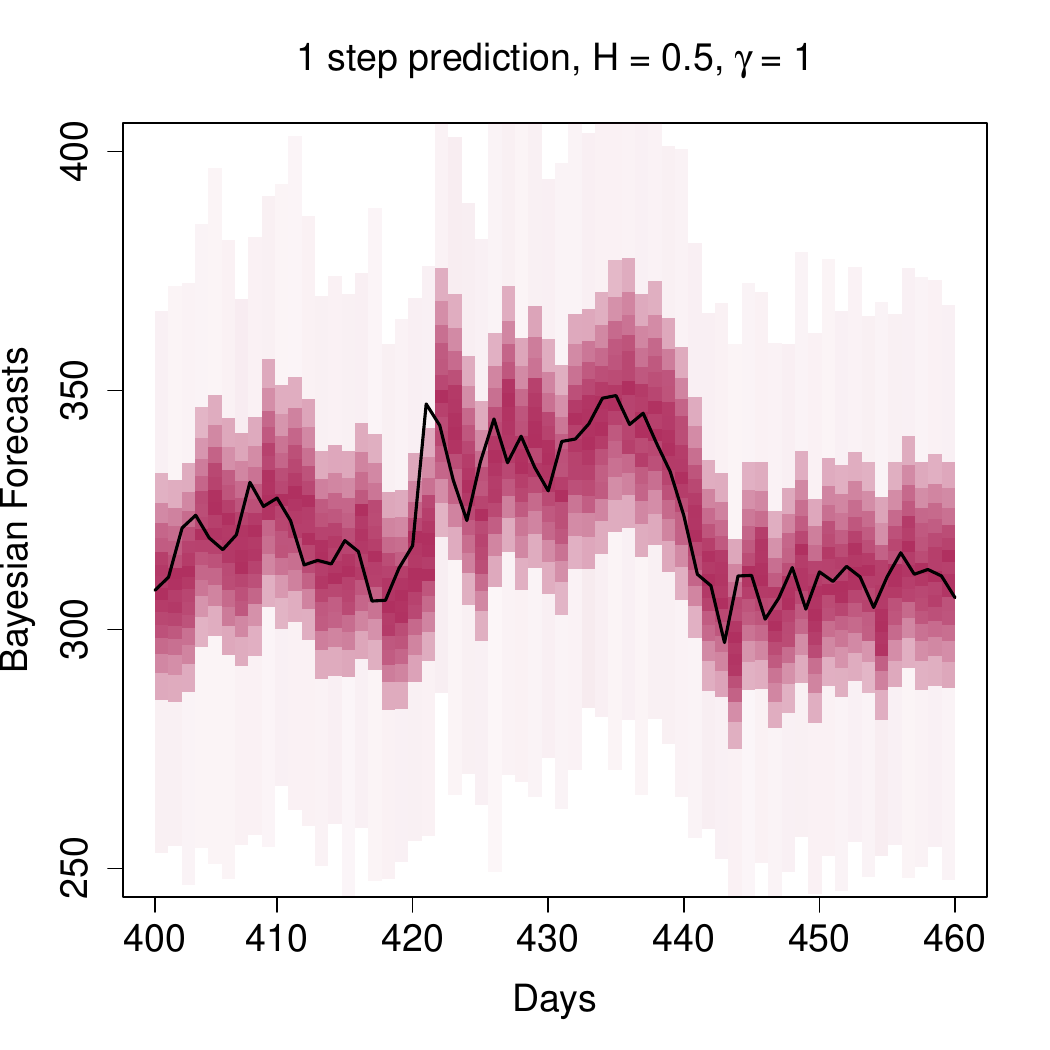}}
\hspace{2mm}
\subfigure []{ \label{fig:p60_bm1}
\includegraphics[width=7.5cm,height=6.5cm]{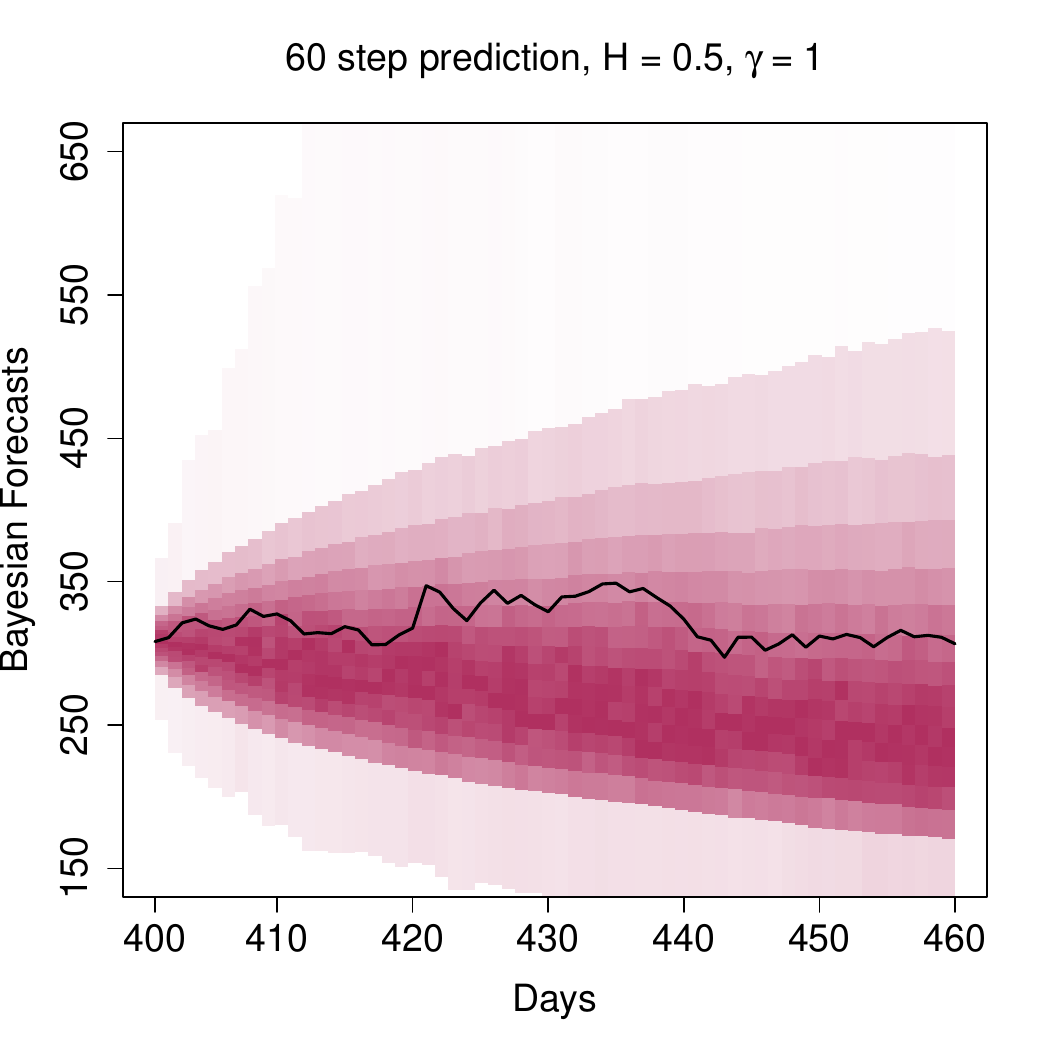}}
\caption{One and 60 steps ahead forecasts in the Bayesian setup corresponding to standard $SDE$.}
\label{fig:1_60_Bayesian_BM}
\end{figure}

The posterior predictive densities associated with fractional and standard $SDE$s are displayed in Figures 
\ref{fig:1_60_Bayesian_fBM} and \ref{fig:1_60_Bayesian_BM}, respectively, where higher posterior predictive densities are represented by progressively intense colours;
the thick, black line stands for the held-out observed data, as before.

\section{Summary, conclusion and future direction}
\label{sec:conclusion}
Although $SDE$s driven by fractional Brownian motion have undeniable importance in modeling continuous time with long-range dependence, hitherto
the statistical literature does not seem to have recognized its full potential. The theory of inference, specifically, asymptotic theory of inference,
seem to be far less developed compared to $SDE$s based on Wiener process. One reason for such negligence is probably the difficulty of statisticians in
imbibing the somewhat involved theory of fractional $SDE$s compared to traditional $SDE$s. Indeed, note that for $SDE$s driven by Wiener process, the
Girsanov formula is readily available, provided that the diffusion coefficient is known. However, this is not the case for fractional $SDE$s if $H$ is unknown.
Thus, most of the methods developed for inference in traditional $SDE$s fail to be carried over to fractional $SDE$s. Discretization is probably the only
way out in this situation, and also realistic, since in practice, the observed data will invariably be in discretized forms. Our increasing domain
infill asymptotics theory established in this article shows that statisticians can expect accurate inference for the model parameters, as well
as the Hurst parameter, both in the classical and Bayesian paradigms, provided that the time domain is sufficiently large and enough observations are collected
in the time domain. Our simulation studies testify to this even when the time domain is not large. The simulation experiments also demonstrate that the Bayesian
framework is preferable. Moreover, our classical and Bayesian applications of both fractional and standard $SDE$ to a real, close price data also brought
out the superiority of our Bayesian fractional $SDE$ over the other model and methods.
However, Theorems \ref{theorem:incons_noncompact1}, \ref{theorem:incons_noncompact2} and Remark \ref{remark:inconsistency} provide the necessary caveat 
against non-judgmental usage of the existing estimators of $H$.

We are currently in the process of extending our work to systems of $SDE$s driven by fractional Brownian motion. In the realm of $SDE$s driven by Wiener processes,
asymptotic theory based on the Girsanov formula have been investigated by these authors in a wide variety of contexts; see 
\ctn{Maitra14a}, \ctn{Maitra14b}, \ctn{Maitra18}, \ctn{Maitra19}, \ctn{Maitra20a}, \ctn{Maitra20b}, \ctn{Maitra20c}. These works deserve to be extended
to $SDE$s driven by fractional Brownian motion. Since the Girsanov formula is ineffective for unknown $H$, our increasing domain infill asymptotics 
framework seems to be a very viable alternative for such extensions. Indeed, in none of the aforementioned works, infill asymptotics has been considered.
As such, our results and their proofs for fractional $SDE$, as well as the practical applications of the relevant methods, are far more intricate, even conceptually, compared
to our previous works on $SDE$ driven by standard Brownian motion.

A somewhat disconcerting issue is that we are unable to obtain the asymptotic distribution of the $MLE$ of $H$. The problem lies in the presence of $H$ even as the
power of $T$, in the likelihood. Although this does not hinder consistency, asymptotic normality or otherwise could not be addressed because of the complications
involved that prevent even a convenient Taylor's series expansion needed for the asymptotic distribution of the $MLE$ of $H$. 
We look forward to resolving the issue in our future work, to be communicated elsewhere.

\section*{Acknowledgment}
We are sincerely grateful to the AE and two referees whose constructive comments have led to significant improvement of the presentation of our manuscript.

\setcounter{section}{0}
\setcounter{figure}{0}
\setcounter{table}{0}

\renewcommand\thefigure{S-\arabic{figure}}
\renewcommand\thetable{S-\arabic{table}}
\renewcommand\thesection{S-\arabic{section}}
\renewcommand\theequation{S-\arabic{equation}}
\renewcommand\thelemma{S-\arabic{lemma}}

\newpage 

\begin{center}
	{\LARGE \bf Supplementary Material}
\end{center}

This document is an addendum to the theory developed in the main manuscript, which we refer to as MB. 
This supplementary material is organized as follows. 

A brief overview of the concept of integration with respect to fractional Brownian motion is provided in Section \ref{sec:integrals}. Illustration of our
key concept with an example involving Ornstein-Uhlenbeck process, is presented in Section \ref{sec:illustration}.
The Gaussian martingale needed for (\ref{eq:fsde1}) of MB is discussed in Section \ref{sec:gaussian_martingale}.
In Section \ref{sec:technical} we present some results and their proofs; these results are used to prove the key results in MB.
Proofs of our results presented in Sections \ref{sec:strong_consistency_MLE} -- \ref{sec:bayesian_consistency} of MB are provided in Sections \ref{sec:theorem_consistency} -- \ref{sec:bayesian_consistency2_KL}.
In Section \ref{sec:schervish} we present the posterior asymptotic normal setup of \ctn{Schervish95}, and in Sections \ref{sec:post_normal1_supp} and \ref{sec:post_normal2}
we present the proofs of Theorems \ref{theorem:post_normal1} and \ref{theorem:post_normal2}, respectively, related to posterior asymptotic normality.

\section{Brief discussion on integrals with respect to fractional Brownian motion}
\label{sec:integrals}
In the context of (\ref{eq:fsde0}), it is important to clarify in what sense $\int_0^T C(t,\tilde X_t)d\tilde W^H_t$ is defined. Since $\tilde W^H_t$ is not a semimartingale for $H\neq 1/2$
(see, for example, \ctn{Oks08}),
the integral can not be defined in the Ito sense. Hence, more general concepts are necessary, the key being the integration theory of Gaussian processes (see, for
example, \ctn{Huang78}).
For simplicity, let us first consider Wiener integrals of deterministic integrands with respect to $\tilde W^H_t$, which can be extended to integrals of random integrands
using concepts of Skorohod integration. We follow \ctn{Norros} and \ctn{Oks08} in this regard.

For $H>1/2$, consider the integral operator
$$\mathcal I f(t)=H(2H-1)\int_0^\infty f(s)|s-t|^{2H-2}ds,$$
and the inner product
$$<f,g>_{\mathcal I}=<f,\mathcal I g>=H(2H-1)\int_0^\infty\int_0^\infty f(s)g(t)|s-t|^{2H-2}dsdt,$$
where $<\cdot,\cdot>$ denotes the usual inner product on $L^2[0,\infty)$, the space of square-integrable functions on $[0,\infty)$.
Let $L^2_{\mathcal I}$ denote the space of equivalence classes of measurable functions $f$ such that $<f,f>_{\mathcal I}<\infty$. Then the association
$\tilde W^H_t\mapsto I_{[0,t)}$, the latter being the indicator function on $[0,t)$, can be extended to an isometry between the Gaussian space generated by 
$\{\tilde W^H_t:t\geq 0\}$ as the smallest closed linear subspace of $L^2(\Omega,\mathcal X,(\mathcal X_t),P)$ containing them, and the function space $L^2_{\mathcal I}$.
For $f\in L^2_{\mathcal I}$, the integral $\int_0^{\infty}f(t)d\tilde W^H_t$ can then be defined as the image of $f$ in this isometry.
For complete details, see \ctn{Oks08}.

For $H<1/2$, the integral in $\mathcal I f(t)$ diverges, and so, $\mathcal I f(t)$ may be re-defined as
$$\mathcal I f(t)=H\int_0^{\infty}|t-s|^{2H-1}\mbox{sgn}(t-s)df(s),$$
where $\mbox{sgn}(t-s)$ denotes the sign of $t-s$. One needs to interpret $f(0-)=0$, where the left hand side denotes the left hand limit of $f$ at $0$.
This implies that $I_{[0,t)}$ is identified with the signed measure $\delta_0-\delta_t$, where $\delta$ is the Dirac delta function. As before, $\tilde W^H_t\mapsto I_{[0,t)}$
defines an isometry, which facilitates the definition of $\int_0^{\infty}f(t)d\tilde W^H_t$ as in the $H>1/2$ case.
Again, the complete details can be found in \ctn{Oks08}.

If $f$ has bounded variation, given $\omega$, the integral $\int_0^Tf(t)d\tilde W^H_t(\omega)$ can be defined as a limit
of Riemann sums, and the integral obtained coincides with the $L^2$ integral almost surely. The convergence of these Riemann sums is equivalent to that of the
usual Riemann-Stieltjes integral appearing on the right-hand side of the integration by parts
formula
$$\int_0^Tf(t)d\tilde W^H_t(\omega)=f(T)\tilde W^H_T(\omega)-f(0)\tilde W^H_0(\omega)-\int_0^T\tilde W^H_t(\omega)df(t).$$
Note that this convergence is guaranteed by the continuity of the sample paths of $\tilde W^H$. For further details, see \ctn{Norros}.

Now, when $f(t)$ is a stochastic process, extension of the above definitions of the integrals is necessary. This is given by Skorohod integration, which uses
the so-called Wiener-Ito-chaos expansion. However, Skorohod integral can also be considered as a limit of Riemann sums associated with Wick products (\ctn{Oks08}).

\section{Illustration of our concept using Ornstein-Uhlenbeck process}
\label{sec:illustration}
In this example wi regard to the OU process, discretizations to the forms (\ref{eq:fsde0_discrete}) and
(\ref{eq:fsde1_discrete}) are not needed since the Girsanov-formula based likelihood is available here since the $SDE$s are driven by ordinary Brownian motion.
However, discretization would lead to the same result.

Assume that data $\{\tilde x_t:t\in [0,T]\}$ is generated from the standard OU process given by
\begin{equation}
d\tilde X_t=-\beta_0\tilde X_t dt + d\tilde W_t,
\label{eq:fsde_eg0}
\end{equation}
and to the observed data $\{\tilde x_t:t\in[0,T]\}$ we fit the following for inferential purpose:
\begin{equation}
dX_t=-\beta\tilde x_t dt + dW_t.
\label{eq:fsde_eg1}
\end{equation}
The corresponding likelihood associated with (\ref{eq:fsde_eg1}) given data $\{\tilde x_t:t\in[0,T]\}$ (after replacing $X$ with $\tilde x$ in the Girsanov formula) is
\begin{align}
L(\beta)=&\exp \left( -\int_0^T \beta\tilde x_td\tilde x_t-\frac{1}{2}\int_0^T (\beta\tilde x_t)^2dt\right)\notag\\
=&\exp\left(\int_0^T\beta\beta_0 \tilde {x_t}^2dt-\int_0^T \beta \tilde x_td\tilde W_t-\frac{1}{2}\int_0^T (\beta\tilde x_t)^2dt \right).\notag
\end{align}
The second equality follows by using (\ref{eq:fsde_eg0}), replacing $\tilde X_t$ with $\tilde x_t$.

The maximum likelihood estimator of $\beta$, obtained by maximizing $L(\beta)$, is given by
$$\hat\beta_T=\beta_0 - \frac{\int_0^T\tilde x_td\tilde W_t}{\int_0^T\tilde {x_t}^2dt}$$
It is easy to see that $\hat\beta_T\rightarrow \beta_0$ in probability as $T\rightarrow\infty$. For, if $\beta_0>0$ (and hence $\tilde x$ is ergodic) then as $T\rightarrow\infty$
$$\frac{1}{T}\int_0^T\tilde {x_t}^2dt\longrightarrow\frac{1}{2\beta_0}$$
almost surely, 
and $$E\left(\frac{1}{\sqrt T}\int_0^T\tilde x_td\tilde W_t\right)^2=\frac{1}{T}\int_0^T E\tilde {x_t}^2dt\rightarrow \frac{1}{2\beta_0}$$
so that $$\hat \beta_T=\beta_0 - \frac{1}{\sqrt T}\frac{\frac{1}{\sqrt T}\int_0^T\tilde x_t d\tilde W_t}{\frac{1}{T}\int_0^T\tilde {x_t}^2dt}\rightarrow \beta_0$$
in probability.
Hence, $\hat\beta_T$ is a consistent estimator.

It is of course clear that working solely with (\ref{eq:fsde_eg0}) would lead to the same estimator and consistency.

\section{Gaussian martingale for (\ref{eq:fsde1})}
\label{sec:gaussian_martingale}

Even though the process $X$ is not a semimartingale, one can associate a semimartingale $Z = \{Z_t, t \geq 0\}$ which is called a fundamental semimartingale such that the natural filtration $(\mathcal Z_t )$ of the process $Z$ coincides with the natural
filtration $(\mathcal X_t )$ of the process $X$ (\ctn{Klept}). Define, for $0<s<t$,
\begin{align}
k_H&=2H\Gamma \left(\frac{3}{2}- H\right)\Gamma \left(H+\frac{1}{2}\right),\label{eq:kh}\\
k_H(t,s)&=k_H^{-1}s^{\frac{1}{2}-H}(t-s)^{\frac{1}{2}-H},\label{eq:khts}\\
\lambda_H&=\frac{2H\Gamma(3-2H)\Gamma \left(H+\frac{1}{2}\right)}{\Gamma \left(\frac{3}{2}- H\right)},\label{eq:lambda_future}\\
w_t^H&=\lambda_H^{-1}t^{2-2H},\label{eq:wth}
\end{align}
and 
\begin{align}
M_t^H&=\int_0^tk_H(t,s)dW_s^H, ~~t\geq 0.\label{eq:mth}
\end{align}

The process $M^H$ is a Gaussian martingale (\ctn{Norros}). Furthermore, the natural filtration of the martingale $M^H$ coincides with the natural filtration of 
the fractional Brownian motion $W^H$.  In fact, the stochastic integral $$\int_0^t C(s,\tilde x_s)dW_s^H$$ can be represented in terms of the stochastic integral with respect to the martingale $M^H.$ For a measurable function $f$ on $[0, T]$, let
\begin{align}
K^f_H(t,s)&=-2H\frac{d}{ds}\int_s^tf(r)r^{H-\frac{1}{2}}(r-s)^{H-\frac{1}{2}}dr,~0\leq s\leq t\label{eq:khf}
\end{align}
when the derivative exists in the sense of absolute continuity with respect to the Lebesgue measure (\ctn{Samko}).

Consider the sample $\{B(t,\tilde x_t)/C(t,\tilde x_t), t \geq 0\}$ to be smooth enough (\ctn{Samko}) such that
\begin{equation}
Q_{\btheta}(t)=\beta\frac{d}{dw^H_t}\int_0^tk_H(t,s)\frac{B(s,\tilde x_s)}{C(s,\tilde x_s)}ds,~t\in[0,T]
\label{eq:qh}
\end{equation}
is well-defined, where $w^H$ and $k_H (t, s)$ are as defined in (\ref{eq:khts}) and (\ref{eq:wth}) respectively
and the derivative is understood in the sense of absolute continuity. 
Note that, given any realization $x$, if $\frac{B(s,\tilde x_s)}{C(s,\tilde x_s)}$ is continuous in $s\in[0,T]$, or even discontinuous at finitely many points $s\in[0,T]$,
then $\int_0^tk_H(t,s)\frac{B(s,\tilde x_s)}{C(s,\tilde x_s)}ds$ exists and is continuous on $[0,T]$. Moreover, if $k_H(t,s)\frac{B(s,\tilde x_s)}{C(s,\tilde x_s)}$ and its partial derivative
with respect to $t$ are continuous on $[0,T]\times[0,T]$, then Leibnitz's rule of differentiation (here, with respect to $t$) under integration applies.
Since $\frac{dw^H_t}{dt}=\lambda^{-1}_H(2-2H)t^{1-2H}$, it then follows by the chain rule of differentiation that
$Q_{\btheta}(t)$ is well-defined. However, conditions for Leibnitz's rule are not necessary for $Q_{\btheta}(t)$ to exist. 
For instance, if $\frac{B(s,\tilde x_s)}{C(s,\tilde x_s)}$ is constant, then 
it follows by first computing the integral in (\ref{eq:qh}) and then using chain rule based differentiation, that 
$Q_{\btheta}(t)=\frac{k^{-1}_H(\frac{3}{2}-H)\tilde B\left(\frac{3}{2}-H,\frac{3}{2}-H\right)}{\lambda^{-1}_H(2-2H)}t^{H-\frac{1}{2}}$, where, for any $a>0$, $b>0$,
$\tilde B(a,b)$ is the Beta function given by $\frac{\Gamma(a)\Gamma(b)}{\Gamma(a+b)}$. 

Next we state the
following theorem by \ctn{Klept} which associates a fundamental semimartingale
$Z$ associated with the process $X$ such that the natural filtration $(\mathcal Z_t )$
coincides with the natural filtration $(\mathcal X_t )$ of $X$.
\begin{theorem}[\ctn{Klept}]
\label{theorem:blsp1}
	Suppose the sample paths of the process $Q_{\btheta}$ defined by (\ref{eq:qh})
	belong $P_{\btheta_0}$-a.s. to $L^2([0, T ], dw^H )$ where $w_H$ is as defined by (\ref{eq:wth}). Let the
process $Z = \{Z_t, t \in [0, T ]\}$ be defined by
\begin{align}
Z_t&=\int_0^tk_H(t,s)[C(s,\tilde x_s)]^{-1}dX_s
\label{eq:zt}
\end{align}
where the function $k_H (t, s)$ is as defined in (\ref{eq:khts}). Then the following results
hold:
\begin{itemize}
\item[(i)] The process $Z$ is an $(\mathcal F_t )$-semimartingale with the decomposition
\begin{align}
	Z_t&=\int_0^tQ_{\btheta}(s)dw^H_s+M_t^H
\label{eq:zt_decompose}
\end{align}
where $M^H$ is the fundamental martingale defined by (\ref{eq:mth}).
\item [(ii)] The process $X$ admits the representation
\begin{align}
X_t&=\int_0^tK^C_H(t,s)dZ_s
\label{eq:yt}
\end{align}
where the function $K^C_H (\cdot, \cdot)$ is as defined in (\ref{eq:khf}).
\item [(iii)] The natural filtrations $(\mathcal Z_t )$ and $(\mathcal X_t )$ coincide.
\end{itemize}
\end{theorem}
By Theorem 3.1 of \ctn{Norros}, the  Gaussian martingale $M^H$ has independent increments with variance function 
$$ E\left[\left(M^H_t\right)^2\right]=c_2^2t^{2-2H},$$ 
where $$c_2=\frac{c_H}{2H(2-2H)^{1/2}} \quad\mbox{with}\quad c_H=\left(\frac{2H\Gamma (\frac{3}{2}-H)}{\Gamma (H+\frac{1}{2})\Gamma (2-2H)}\right)^{1/2}.$$

\section{Some technical results}
\label{sec:technical}
\begin{lemma}
\label{lemma:lemma1}
If $a_k\sim b_k$ as $k\rightarrow\infty$, then 
it holds that $\frac{1}{r_k}\sum_{i=1}^{m_k}a_k\sim\frac{1}{r_k}\sum_{i=1}^{m_k}b_k$, as $k\rightarrow\infty$, provided that
	$m_k\rightarrow\infty$, $r_k\rightarrow\infty$ and $\frac{1}{r_k}\sum_{i=1}^{m_k}b_k$ remains non-zero, as $k\rightarrow\infty$. 
\end{lemma}	
\begin{proof}	
Note that $a_k\sim b_k$ implies that for any $\epsilon>0$, there exists $k_0\geq 1$ depending upon $\epsilon$ such that for $k> k_0$, 
$(1-\epsilon)b_k<a_k<(1+\epsilon)b_k$. Hence, for $k>k_0$, 
$(1-\epsilon)\frac{1}{r_k}\sum_{i=m_{k_0}+1}^{m_k}b_k<\frac{1}{r_k}\sum_{i=m_{k_0}+1}^{m_k}a_k<(1+\epsilon)\frac{1}{r_k}\sum_{i=m_{k_0}+1}^{m_k}b_k$.
That is,
\begin{equation}
	\frac{1}{r_k}\sum_{i=m_{k_0}+1}^{m_k}a_k\sim \frac{1}{r_k}\sum_{i=m_{k_0}+1}^{m_k}b_k,~\mbox{as}~k\rightarrow\infty.
	\label{eq:ae0}
\end{equation}
Since by hypothesis, $\frac{1}{r_k}\sum_{i=1}^{m_k}b_k$ remains non-zero, as $k\rightarrow\infty$, it follows that
	$\frac{\sum_{i=1}^{m_{k_0}}a_k}{\sum_{i=1}^{m_k}b_k}\rightarrow 0$ and $\frac{\sum_{i=1}^{m_{k_0}}b_k}{\sum_{i=1}^{m_k}b_k}\rightarrow 0$, as $k\rightarrow\infty$.
	Using these and (\ref{eq:ae0}), it follows that
	\begin{align}
        \frac{\frac{1}{r_k}\sum_{i=1}^{m_k}a_k}{\frac{1}{r_k}\sum_{i=1}^{m_k}b_k}
		=\frac{\sum_{i=1}^{m_{k_0}}a_k+\sum_{i=m_{k_0}+1}^{m_k}a_k}{\sum_{i=1}^{m_{k_0}}b_k+\sum_{i=m_{k_0}+1}^{m_k}b_k}
		\rightarrow 1,~\mbox{as}~k\rightarrow\infty.
	\end{align}
This completes the proof.
%
\end{proof}

\begin{lemma}
\label{lemma:unicon}
	Let $\left\{f_k:k\geq 1\right\}$ be a sequence of functions such that 
	$f_k(x)\rightarrow f(x)$ uniformly in $x\in\mathcal X$, as $k\rightarrow\infty$, for some $\mathcal X\subseteq \mathbb R$, where $f$ is continuous.
	Suppose that $x_k\rightarrow\tilde x$, as $k\rightarrow\infty$. Then $f_k(x_k)\rightarrow f(\tilde x)$, as $k\rightarrow\infty$.
\end{lemma}
\begin{proof}
Note that
\begin{equation}
	\left|f_k(x_k)-f(\tilde x)\right|\leq \left|f_k(x_k)-f(x_k)\right|+\left|f(x_k)-f(\tilde x)\right|.
	\label{eq:unicon1}
\end{equation}
	Now, due to uniform convergence, there exists $k_0(\epsilon)\geq 1$ such that for $k\geq k_0(\epsilon)$,
\begin{equation}
	\left|f_k(x_k)-f(x_k)\right|\leq\underset{x\in\mathcal X}{\sup}~\left|f_k(x)-f(x)\right|<\frac{\epsilon}{2}.
	\label{eq:unicon2}
\end{equation}
Also, due to continuity of $f$, there exists $k_1(\epsilon)\geq 1$ such that for $k\geq k_1(\epsilon)$,
\begin{equation}
	\left|f(x_k)-f(\tilde x)\right|<\frac{\epsilon}{2}.
	\label{eq:unicon3}
\end{equation}
	Combining (\ref{eq:unicon1}), (\ref{eq:unicon2}) and (\ref{eq:unicon3}) we see that for $k\geq \max\{k_0(\epsilon),k_1(\epsilon)\}$,
	$$\left|f_k(x_k)-f(\tilde x)\right|<\epsilon.$$
	This completes the proof.
\end{proof}


\section{Proof of Theorem \ref{theorem:consistency}}
\label{sec:theorem_consistency}
\begin{proof}
Let $\ell_{T,n}(\btheta)=-2\log L_{T,n}(\btheta)$. 
We shall show that
$\underset{\btheta\in\bTheta:\|\btheta-\btheta_0\|\geq\delta}{\inf} \left\{\ell_{T,n}(\btheta)-\ell_{T,n}(\btheta_0)\right\}\stackrel{a.s.}{\longrightarrow} \infty$ as
$T\rightarrow\infty$, $n\rightarrow\infty$ such that $T^2/n\rightarrow 0$.
Now,
$$ \ell_{T,n}(\btheta)=n \log {2\pi} + \sum_{i=1}^{n} \log {\upsilon^{2}_{M^{H}_{t_i}}}+ \sum_{i=1}^{n} \frac{(\Delta  M^{H}_{t_i})^2}{\left(\upsilon_{M^{H}_{t_i}}\right)^2}$$
Hence,
\begin{align}
\ell_{T,n}(\btheta)-\ell_{T,n}(\btheta_0)&= \sum_{i=1}^{n} \log \frac{\upsilon^2_{M^{H}_{t_i}}}{\upsilon^2_{M^{H_0}_{t_i}}}+ \sum_{i=1}^{n} \frac{(\Delta  M^{H}_{t_i})^2}{\upsilon^2_{M^{H}_{t_i}}}-\sum_{i=1}^{n} \frac{(\Delta  M^{H_0}_{t_i})^2}{\upsilon^2_{M^{H_0}_{t_i}}}\notag\\
&=-\sum_{i=1}^{n} \log \frac{\upsilon^2_{M^{H_0}_{t_i}}}{\upsilon^2_{M^{H}_{t_i}}}+ \sum_{i=1}^{n} \frac{(\Delta  M^{H}_{t_i}-\Delta  M^{H_0}_{t_i}+\Delta  M^{H_0}_{t_i})^2}{\upsilon^2_{M^{H}_{t_i}}}-\sum_{i=1}^{n} \frac{(\Delta  M^{H_0}_{t_i})^2}{\upsilon^2_{M^{H_0}_{t_i}}}\notag\\
&=-\sum_{i=1}^{n} \log \frac{\upsilon^2_{M^{H_0}_{t_i}}}{\upsilon^2_{M^{H}_{t_i}}}-\sum_{i=1}^{n} \frac{(\Delta  M^{H_0}_{t_i})^2}{\upsilon^2_{M^{H_0}_{t_i}}}+ \sum_{i=1}^{n} \frac{(\Delta  M^{H}_{t_i}-\Delta  M^{H_0}_{t_i})^2}{\upsilon^2_{M^{H_0}_{t_i}}}\frac{\upsilon^2_{M^{H_0}_{t_i}}}{\upsilon^2_{M^{H}_{t_i}}}\notag\\
&+2\sum_{i=1}^{n} \frac{(\Delta  M^{H}_{t_i}-\Delta  M^{H_0}_{t_i})\Delta  M^{H_0}_{t_i}}{\upsilon^2_{M^{H_0}_{t_i}}}\frac{\upsilon^2_{M^{H_0}_{t_i}}}{\upsilon^2_{M^{H}_{t_i}}}+\sum_{i=1}^{n} \frac{(\Delta  M^{H_0}_{t_i})^2}{\upsilon^2_{M^{H_0}_{t_i}}}\frac{\upsilon^2_{M^{H_0}_{t_i}}}{\upsilon^2_{M^{H}_{t_i}}}
\label{eq:likelihood_diff}
\end{align}
Recall that $H \in [\eta, 1-\eta]$ and $\beta\in\mathfrak B$, where $\mathfrak B$ is compact. 
	Hence for any given data sequence $\left\{x_{t_i};i=1,\ldots,n\right\}$, the infimum of $\ell_{T,n}(\btheta)-\ell_{T,n}(\btheta_0)$ with respect to $\btheta\in\bTheta$
	such that $\|\btheta-\btheta_0\|\geq \delta$, 
	will yield a sequence $\left\{\btheta_{T,n}\in\bTheta:T>0, n\geq 1\right\}$ such that $T^2/n\rightarrow 0$. Since $\bTheta$ is compact, 
	there must exist a convergent subsequence $\left\{\btheta_{T_k,n_k}=(\beta_{T_k,n_k},H_{T_k,n_k}):k=1,2,\ldots\right\}$ of the above sequence
	which converges to a limit, say $\tilde\btheta=(\tilde\beta,\tilde {H})$, as $T\rightarrow\infty$, $n\rightarrow\infty$ 
	such that $T^2/n\rightarrow 0$. 
	Thus, $\beta_{T_k,n_k}\sim \tilde \beta$ and $H_{T_k,n_k}\sim \tilde H$, almost surely, as $k\rightarrow\infty$.
	For simplicity of notation, we shall denote $\beta_{T_k,n_k}$, $H_{T_k,n_k}$, $\btheta_{T_k,n_k}$, $\ell_{T_k,n_k}(\btheta_{T_k,n_k})$
	and $\ell_{T_k,n_k}(\btheta_0)$ by $\beta_k$, $H_k$, $\btheta_k$, $\ell_k(\btheta_k)$ and $\ell_k(\btheta_0)$, respectively.
	If we can show that 
\begin{equation}
		\ell_k(\btheta_k)-\ell_k(\btheta_0)\rightarrow \infty,~\mbox{almost surely, as}~k\rightarrow\infty, 
		\label{eq:subseq1}
	\end{equation}
		then it would follow that 
	$\underset{\btheta\in\bTheta:\|\btheta-\btheta_0\|\geq \delta}{\inf} \left\{\ell_{T,n}(\btheta)-\ell_{T,n}(\btheta_0)\right\}\stackrel{a.s.}{\longrightarrow} \infty$ as
$T\rightarrow\infty$, $n\rightarrow\infty$ such that $T^2/n\rightarrow 0$.
In what follows, we show that (\ref{eq:subseq1}) indeed holds.\\
	For $t_{ik}=T_k\frac{i}{n_k}$ where $i=1,\ldots,n_k$, let $Y_{t_{ik}}=\frac{\Delta  M^{H_0}_{t_{ik}}}{\upsilon_{M^{H_0}_{t_{ik}}}}$ . 
	Note that $Y_{t_{ik}}$ for $i=1,2,\ldots, n_k$, are independent $N(0,1)$ random variables. Now observe 
	that $t_{ik}=T_k\frac{i}{n_k}$ for $i=1,\ldots,n_k$. By the mean value theorem it follows that,
\begin{align}
t_{i+1,k}^{2-2H} - t_{ik}^{2-2H}&=\left(\frac{T_k(i+1)}{n_k}\right)^{2-2H}-\left(\frac{T_ki}{n_k}\right)^{2-2H}\notag\\
	&=\frac{T_k}{n_k}(2-2H)\xi^{1-2H}_{ik}\quad\mbox{(where $\xi_{ik} \in (t_{ik}, t_{i+1,k}$))}.
\label{eq:mvt}
\end{align}
Now note that if $H\leq 1/2$, then $T^{1-2H}_k\left(\frac{i}{n_k}\right)^{1-2H}<\xi^{1-2H}_{ik}<T^{1-2H}_k\left(\frac{i+1}{n_k}\right)^{1-2H}$.
	If on the other hand $H>1/2$, then $T^{1-2H}_k\left(\frac{i+1}{n_k}\right)^{1-2H}<\xi^{1-2H}_{ik}<T^{1-2H}_k\left(\frac{i}{n_k}\right)^{1-2H}$.
In both the cases, $\underset{k\rightarrow\infty;i\rightarrow n_k}{\lim}~\left(\frac{i}{n_k}\right)^{1-2H}
=\underset{k\rightarrow\infty;i\rightarrow n_k}{\lim}~\left(\frac{i+1}{n_k}\right)^{1-2H}=1$. Hence, 
\begin{equation}
	\xi^{1-2H}_{ik}\sim T^{1-2H}_k,~\mbox{as}~i\rightarrow n_k~\mbox{and}~k\rightarrow\infty.
	\label{eq:mvt2}
\end{equation}
Hence, from (\ref{eq:mvt}) and (\ref{eq:mvt2}) we obtain
\begin{equation}
	n_k\left( t_{i+1,k}^{2-2H} - t_{ik}^{2-2H}\right)\sim T^{(2-2H)}_k(2-2H),~\mbox{as}~i\rightarrow n_k;~k\rightarrow\infty.
\label{eq:mvt3}
\end{equation}
Note that the same holds even if $H$ is replaced with $H_k$, since $H_k\in[\eta,1-\eta]$.
Using (\ref{eq:mvt3}) we obtain
\begin{align}
	\frac{\upsilon_{M^{H_0}_{t_{ik}}}^2}{\upsilon_{M^{H_k}_{t_{ik}}}^2}
	&=\frac{\frac{c^2_{H_0}}{4H^2_0(2-2H_0)}n_k(t_{i+1,k}^{2-2H_0} 
	- t_{ik}^{2-2H_0})}{\frac{c^2_{H_k}}{4H^2_k(2-2H_k)}n_k(t_{i+1,k}^{2-2H_k} - t_{ik}^{2-2H_k})}\notag\\
	&\sim \frac{\frac{c^2_{H_0}}{4H^2_0(2-2H_0)}(2-2H_0)T^{2-2H_0}_k}{\frac{c^2_{\tilde H}}{4{\tilde H}^2(2-2\tilde {H})}(2-2\tilde H)T^{2-2H_k}_k}\notag\\
	&=\frac{\frac{c^2_{H_0}}{H^2_0}}{\frac{c^2_{\tilde H}}{{\tilde H}^2}}T^{2(H_k-H_0)}_k,~\mbox{as}~i\rightarrow n_k;~k\rightarrow\infty.
\label{eq:term1}
\end{align}
Hence,
\begin{equation}
\underset{i\rightarrow n_k;~k\rightarrow\infty}{\lim}
	\frac{1}{T^{2(H_k-H_0)}_k}\frac{\upsilon_{M^{H_0}_{t_{ik}}}^2}{\upsilon_{M^{H_k}_{t_{ik}}}^2}
	=\frac{c^2_{H_0}{\tilde H}^2}{c^2_{\tilde H} H^2_0}=K ~\mbox {(say)}.
\label{eq:limit_term1}
\end{equation}
Now, equation (\ref{eq:likelihood_diff}) can be re-written as 
\begin{align}
	\ell_k(\btheta_k)-\ell_k(\btheta_0)&=-\sum_{i=1}^{n_k} \log \frac{\upsilon^2_{M^{H_0}_{t_{ik}}}}{\upsilon^2_{M^{H_k}_{t_{ik}}}}
	-\sum_{i=1}^{n_k} Y^2_{t_{ik}} + \sum_{i=1}^{n_k} Y^2_{t_{ik}}\frac{\upsilon^2_{M^{H_0}_{t_{ik}}}}{\upsilon^2_{M^{H_k}_{t_{ik}}}}+\notag\\
	&\sum_{i=1}^{n_k} \frac{(\Delta  M^{H_k}_{t_{ik}}-\Delta  M^{H_0}_{t_{ik}})^2}{\upsilon^2_{M^{H_0}_{t_{ik}}}}
	\frac{\upsilon^2_{M^{H_0}_{t_{ik}}}}{\upsilon^2_{M^{H_k}_{t_{ik}}}}
	+2\sum_{i=1}^{n_k} \frac{(\Delta  M^{H_k}_{t_{ik}}-\Delta  M^{H_0}_{t_{ik}})\Delta  M^{H_0}_{t_{ik}}}{\upsilon^2_{M^{H_0}_{t_{ik}}}}
	\frac{\upsilon^2_{M^{H_0}_{t_{ik}}}}{\upsilon^2_{M^{H_k}_{t_{ik}}}}.
\label{eq:likelihood_diff1}
\end{align}
Consider first three terms of the right hand side of above equation, namely,
\begin{align}
	&-\sum_{i=1}^{n_k} \log \frac{\upsilon^2_{M^{H_0}_{t_{ik}}}}{\upsilon^2_{M^{H_k}_{t_{ik}}}}
	-\sum_{i=1}^{n_k} Y^2_{t_{ik}}+ \sum_{i=1}^{n_k} Y^2_{t_{ik}}\frac{\upsilon^2_{M^{H_0}_{t_{ik}}}}{\upsilon^2_{M^{H_k}_{t_{ik}}}}\notag\\
	&=n_k\left(\frac{T^{2(H_k-H_0)}_k}{n_k}\sum_{i=1}^{n_k} Y^2_{t_{ik}}\left(\frac{1}{T^{2(H_k-H_0)}_k}
	\frac{\upsilon^2_{M^{H_0}_{t_{ik}}}}{\upsilon^2_{M^{H_k}_{t_{ik}}}}\right)
	-\frac{1}{n_k}\sum_{i=1}^{n_k} \log \left(\frac{1}{T^{2(H_k-H_0)}_k}\frac{\upsilon^2_{M^{H_0}_{t_{ik}}}}{\upsilon^2_{M^{H_k}_{t_{ik}}}}\right)\right.\notag\\
	&\qquad\qquad\left.-\log T^{2(H_k-H_0)}_k-\frac{1}{n_k}\sum_{i=1}^{n_k} Y^2_{t_{ik}}\right).
\label{eq:3terms}
\end{align}
Now, due to (\ref{eq:limit_term1}), $Y^2_{t_{ik}}\left(\frac{1}{T^{2(H_k-H_0)}_k} \frac{\upsilon^2_{M^{H_0}_{t_{ik}}}}{\upsilon^2_{M^{H_k}_{t_{ik}}}}\right)
\sim Y^2_{t_{ik}} K$, as $i\rightarrow n_k$; $k\rightarrow\infty$. Then Lemma \ref{lemma:lemma1} 
implies
\begin{equation}
	\frac{T^{2(H_k-H_0)}_k}{n_k}\sum_{i=1}^{n_k} Y^2_{t_{ik}}\left(\frac{1}{T^{2(H_k-H_0)}_k}
	\frac{\upsilon^2_{M^{H_0}_{t_{ik}}}}{\upsilon^2_{M^{H_k}_{t_{ik}}}}\right)
	\sim T^{2(H_k-H_0)}_kK\frac{\sum_{i=1}^{n_k} Y^2_{t_{ik}}}{n_k}
	~\mbox{as}~k\rightarrow\infty.
	\label{eq:term1_lim0}
\end{equation}
However, the strong law of large numbers applied to independent standard normal random variables $Y_{t_{ik}}$ on the right hand side of (\ref{eq:term1_lim0}) shows that
$\frac{\sum_{i=1}^{n_k} Y^2_{t_{ik}}}{n_k}\rightarrow 1$, almost surely, as $k\rightarrow\infty$.
In other words, it holds that
\begin{equation}
	\frac{T^{2(H_k-H_0)}_k}{n_k}\sum_{i=1}^{n_k} Y^2_{t_{ik}}\left(\frac{1}{T^{2(H_k-H_0)}_k}
	\frac{\upsilon^2_{M^{H_0}_{t_{ik}}}}{\upsilon^2_{M^{H_k}_{t_{ik}}}}\right)
	\sim KT^{2(H_k-H_0)}_k,~\mbox{almost surely, as}~k\rightarrow\infty.
	\label{eq:term1_lim}
\end{equation}

As for the second term of (\ref{eq:3terms}), note that due to (\ref{eq:limit_term1}) and Lemma \ref{lemma:lemma1},
\begin{equation}
	-\frac{1}{n_k}\sum_{i=1}^{n_k} \log \left(\frac{1}{T^{2(H_k-H_0)}_k}\frac{\upsilon^2_{M^{H_0}_{t_{ik}}}}{\upsilon^2_{M^{H_k}_{t_{ik}}}}\right)
	\sim -\log K,~\mbox{as}~k\rightarrow\infty.
	\label{eq:term2_lim}
\end{equation}
Since the fourth term of (\ref{eq:3terms}) tends to $1$ as $k\rightarrow\infty$ by the strong law of large numbers, then
along with (\ref{eq:term1_lim}) and (\ref{eq:term2_lim}), it follows that (\ref{eq:3terms}) has the following asymptotic equivalence: 
\begin{align}
	&-\sum_{i=1}^{n_k} \log \frac{\upsilon^2_{M^{H_0}_{t_{ik}}}}{\upsilon^2_{M^{H_k}_{t_{ik}}}}
	-\sum_{i=1}^{n_k} Y^2_{t_{ik}}+ \sum_{i=1}^{n_k} Y^2_{t_{ik}}\frac{\upsilon^2_{M^{H_0}_{t_{ik}}}}{\upsilon^2_{M^{H_k}_{t_{ik}}}}\notag\\
	&\qquad\sim n_kK T^{2(H_k-H_0)}_k-n_k\log K -n_k\log T^{2(H_k-H_0)}_k-n_k,~\mbox{as}~k\rightarrow\infty.
	\label{eq:sum0}
\end{align}
From equation (\ref{eq:qh}), in conjunction with  (\ref{eq:wth}), we obtain
\begin{align}
Q_{\btheta}(t)\Delta w_{t}^{H}&=\beta\frac{k_{H}^{-1}({(t+1)}^{2-2H}-t^{2-2H})}{t^{1-2H}(2-2H)}
	\frac{d}{dt}\int_0^{t} s^{\frac{1}{2}-H}(t-s)^{\frac{1}{2}-H}\frac{B(s,\tilde x_s)}{C(s,\tilde x_s)}ds\notag\\
	&=\beta\left[\frac{k_{H}^{-1}({(t+1)}^{2-2H}-t^{2-2H})}{(2-2H)}\right]\left[(2-2H)\int_0^1 u^{\frac{1}{2}- H}(1-u)^{\frac{1}{2}- H}
	\frac{B(tu,\tilde x_{tu})}{C(tu,\tilde x_{tu})}du\right.\notag\\
	&\qquad\left.+t\frac{d}{dt}\int_0^1 u^{\frac{1}{2}- H}(1-u)^{\frac{1}{2}- H}\frac{B(tu,\tilde x_{tu})}{C(tu,\tilde x_{tu})}du\right].
\label{eq:qhwth}
\end{align}
In the last equality we have used the substitution $u=\frac{s}{t}$.\\
Now consider the last term of equation (\ref{eq:likelihood_diff1}). Using (\ref{eq:qhwth}) we obtain 
\begin{align}
	&2\sum_{i=1}^{n_k}\frac{(\Delta  M^{H_k}_{t_{ik}}-\Delta  M^{H_0}_{t_{ik}})\Delta  M^{H_0}_{t_{ik}}}{\upsilon^2_{M^{H_0}_{t_{ik}}}}
	\frac{\upsilon^2_{M^{H_0}_{t_{ik}}}}{\upsilon^2_{M^{H_k}_{t_{ik}}}}\notag\\
	&=2\sum_{i=1}^{n_k}(\Delta  M^{H_k}_{t_{ik}}-\Delta  M^{H_0}_{t_{ik}})Y_{t_{ik}}\frac{\upsilon_{M^{H_0}_{t_{ik}}}}{\upsilon^2_{M^{H_k}_{t_{ik}}}}\notag\\
	&=2\sum_{i=1}^{n_k}\frac{(Q_{\btheta_0}(t_{ik})\Delta w_{t_{ik}}^{ H_0}-Q_{\btheta_k}(t_{ik})\Delta w_{t_{ik}}^{H_k}+\Delta Z^{\btheta_k}_{t_{ik}}
	-\Delta Z^{\btheta_0}_{t_{ik}})}
	{\frac{c^2_{H_k}}{4{H_k}^2(2-2H_k)}(t_{i+1,k}^{2-2H_k} - t_{ik}^{2-2H_k})}Y_{t_{ik}}
	\frac{c_{H_0}}{2H_0\sqrt{(2-2H_0)}}\sqrt{(t_{i+1,k}^{2-2H_0} - t_{ik}^{2-2H_0})}\notag\\
	&=\sum_{i=1}^{n_k} Y_{t_{ik}}
	\frac{c_{H_0}4{H_k}^2(2-2H_k)}{c^2_{H_k}H_0\sqrt{(2-2H_0)}}\frac{\sqrt{(t_{i+1,k}^{2-2H_0} - t_{ik}^{2-2H_0})}}{(t_{i+1,k}^{2-2H_k} 
	- t_{ik}^{2-2H_k})}\left[\beta_0\frac{k_{ H_0}^{-1}(t_{i+1,k}^{2-2 H_k}-t_{ik}^{2-2 H_0})}{t_{ik}^{1-2 H_0}(2-2H_0)}\right.\notag\\
	&\quad\quad\left.\frac{d}{dt}\left.\int_0^{t} s^{\frac{1}{2}- H_0}(t-s)^{\frac{1}{2}- H_0}\frac{B(s,\tilde x_s)}{C(s,\tilde x_s)}ds\right|_{t=t_{ik}}\right.\notag\\
	&\quad\quad\left.-\beta_k\frac{k_{H_k}^{-1}(t_{i+1,k}^{2-2H_k}
	-t_{ik}^{2-2H_k})}{t_{ik}^{1-2H_k}(2-2H_k)}\frac{d}{dt}\left.\int_0^{t} s^{\frac{1}{2}-H_k}(t-s)^{\frac{1}{2}-H_k}\frac{B(s,\tilde x_s)}{C(s,\tilde x_s)}ds\right|_{t=t_{ik}}
        +\Delta Z^{\btheta_k}_{t_{ik}}-\Delta Z^{\btheta_0}_{t_{ik}}
	\right]
	\notag\\
	&=\sum_{i=1}^{n_k}Y_{t_{ik}}\tilde w_{ik},
\label{eq:limit_term3}
\end{align}
where 
\begin{align}
	&\tilde w_{ik}
	=
	\frac{c_{H_0}4{H_k}^2(2-2H_k)}{c^2_{H_k}H_0\sqrt{(2-2H_0)}}\sqrt{(t_{i+1,k}^{2-2H_0} - t_{ik}^{2-2H_0})}
	\left[\frac{\beta_0 k_{ H_0}^{-1}}{(2-2H_0)}\left(\frac{t_{i+1,k}^{2-2 H_0}-t_{ik}^{2-2 H_0}}{t_{i+1,k}^{2-2H_k} - t_{ik}^{2-2H_k}}\right)\right.\notag\\
	&\left.\left\{(2-2H_0)\int_0^1 u^{\frac{1}{2}- H_0}(1-u)^{\frac{1}{2}- H_0}\frac{B(t_{ik}u,\tilde x_{t_{ik}u})}{C(t_{ik}u,\tilde x_{t_{ik}u})}du 
	+t_{ik}\frac{d}{dt}\int_0^1 u^{\frac{1}{2}- H_0}(1-u)^{\frac{1}{2}- H_0}\frac{B(tu,\tilde x_{tu})}{C(tu,\tilde x_{tu})}du\bigg|_{t=t_{ik}}\right\}\right.\notag\\
	&\left.-\frac{\beta_kk_{H_k}^{-1}}{2-2H_k} \left\{(2-2H_k)\int_0^1 u^{\frac{1}{2}- H_k}(1-u)^{\frac{1}{2}- H_k}\frac{B(t_{ik}u,\tilde x_{t_{ik}u})}{C(t_{ik}u,\tilde x_{t_{ik}u})}du
	\right.\right.\notag\\
	& \left.\left. + t_{ik}\frac{d}{dt} \int_0^1 u^{\frac{1}{2}- H_k}(1-u)^{\frac{1}{2}- H_k}\frac{B(tu,\tilde x_{tu})}{C(tu,\tilde x_{tu})}du\bigg|_{t=t_{ik}} \right\}
	+\Delta Z^{\btheta_k}_{t_{ik}}-\Delta Z^{\btheta_0}_{t_{ik}}\right].
	\label{eq:w0}
\end{align}
It follows from (\ref{eq:ass_lim1_general}) and (\ref{eq:ass_lim2_general}) that
\begin{align}
	&\int_0^1 u^{\frac{1}{2}-H_0}(1-u)^{\frac{1}{2}-H_0}\frac{B(t_{ik}u,\tilde x_{t_{ik}u})}{C(t_{ik}u,\tilde x_{t_{ik}u})}du\rightarrow C_1(H_0,\tilde x),~\mbox{almost surely, as}~
	i\rightarrow n_k;~k\rightarrow\infty;
	\label{eq:ass_lim1}\\
	&\frac{d}{dt}\left.\int_0^1 u^{\frac{1}{2}- H_0}(1-u)^{\frac{1}{2}- H_0}\frac{B(tu,\tilde x_{tu})}{C(tu,\tilde x_{tu})}du\right|_{t=t_{ik}}
	\rightarrow C_2(H_0,\tilde x),~\mbox{almost surely, as}~i\rightarrow n_k;~k\rightarrow\infty;
	\label{eq:ass_lim2}
\end{align}
Now, since (\ref{eq:ass_lim1_general}) and (\ref{eq:ass_lim2_general}) hold uniformly in $H\in[\eta,1-\eta]$ and $C_1$, $C_2$ are continuous in $H$, it follows 
from Lemma \ref{lemma:unicon} (see Appendix) that
\begin{align}
	&\int_0^1 u^{\frac{1}{2}-H_k}(1-u)^{\frac{1}{2}-H_k}\frac{B(t_{ik}u,\tilde x_{t_{ik}u})}{C(t_{ik}u,\tilde x_{t_{ik}u})}du\rightarrow C_1(\tilde H,\tilde x),~\mbox{almost surely, as}~
	i\rightarrow n_k;~k\rightarrow\infty;
	\label{eq:ass_lim3}\\
	&\frac{d}{dt}\left.\int_0^1 u^{\frac{1}{2}- H_k}(1-u)^{\frac{1}{2}- H_k}\frac{B(tu,\tilde x_{tu})}{C(tu,\tilde x_{tu})}du\right|_{t=t_{ik}}
	\rightarrow C_2(\tilde H,\tilde x),~\mbox{almost surely, as}~i\rightarrow n_k;~k\rightarrow\infty.
	\label{eq:ass_lim4}
\end{align}
Using (\ref{eq:ass_lim1}) -- (\ref{eq:ass_lim4}) and the result $t^{2-2H}_{i+1,k}-t^{2-2H}_{ik}\sim (2-2H)\frac{T^{2-2H}_k}{n_k}$ 
in (\ref{eq:w0}), we obtain
\begin{align}
	\tilde w_{ik}&\sim \frac{c_{H_0}}{c^2_{\tilde H}}\frac{4\tilde H^2}{H_0}(2-2\tilde H) \frac{T^{1-H_0}_k}{\sqrt{n_k}}
	\times\left[\frac{\beta_0k^{-1}_{H_0}}{(2-2H_0)}\left(\frac{1-H_0}{1-\tilde H}\right)T^{2(H_k-H_0)}_k\left\{(2-2H_0)C_1(H_0,\tilde x)+T_kC_2(H_0,\tilde x)\right\}\right.\notag\\
	&\qquad\qquad\left.-\frac{\tilde\beta k^{-1}_{\tilde H}}{(2-2\tilde H)}\left\{(2-2\tilde H)C_1(\tilde H,\tilde x)+T_kC_2(\tilde H,\tilde x)\right\}
	+\Delta Z^{\btheta_k}_{t_{ik}}-\Delta Z^{\btheta_0}_{t_{ik}}\right],\notag\\
	&\qquad\qquad~\mbox{almost surely, as}~i\rightarrow n_k;~k\rightarrow\infty.
	\label{eq:w1}
\end{align}
Now observe that due to (\ref{eq:discretize2}) and normality of $\Delta M^{H}_{t_{ik}}$, 
$\Delta Z^{\btheta}_{t_{ik}}\sim N\left(Q_{\btheta}(t_{ik})\Delta w^{H}_{t_{ik}},\upsilon^2_{M^{H}_{t_{ik}}}\right)$.
Using (\ref{eq:ass_lim1}) -- (\ref{eq:ass_lim4}) and the result $t^{2-2H}_{i+1,k}-t^{2-2H}_{ik}\sim (2-2H)\frac{T^{2-2H}_k}{n_k}$, it follows
that $\frac{\Delta Z^{\btheta_k}_{t_{ik}}-Q_{\btheta_k}(t_{ik})\Delta w^{H_k}_{t_{ik}}}{c_2\sqrt{2-2H_k}T^{1-H_k}_k/\sqrt{n_k}}\sim N(0,1)$, 
almost surely, as $k\rightarrow\infty$. Since a $N(0,1)$ random variable is almost surely finite, and since $T^{1+u_k}_k\rightarrow\infty$ as $k\rightarrow\infty$
for any non-negative-sequence $u_k$, it follows by dividing both sides of the above by $T^{1+u_k}_k$, that
$\frac{\Delta Z^{\btheta_k}_{t_{ik}}-Q_{\btheta_k}(t_{ik})\Delta w^{H_k}_{t_{ik}}}{c_2\sqrt{2-2H_k}T^{2-H_k+u_k}_k/\sqrt{n_k}}\rightarrow 0$, 
almost surely, as $k\rightarrow\infty$,
for any non-negative sequence $u_k$. Using (\ref{eq:qhwth}),
(\ref{eq:ass_lim1}) -- (\ref{eq:ass_lim4}) and the result $t^{2-2H_k}_{i+1,k}-t^{2-2H_k}_{ik}\sim (2-2H_k)\frac{T^{2-2H_k}_k}{n_k}$, it is also easy
to see that $\frac{Q_{\btheta_k}(t_{ik})\Delta w^{H_k}_{t_{ik}}}{c_2\sqrt{2-2H_k}T^{2-H_k+u_k}_k/\sqrt{n_k}}\rightarrow 0$, almost surely, as $k\rightarrow\infty$, so that
$\frac{\Delta Z^{\btheta_k}_{t_{ik}}}{T^{2-H_k+u_k}_k/\sqrt{n_k}}\rightarrow 0$, almost surely, as $k\rightarrow\infty$.
In the same way, 
$\frac{\Delta Z^{\btheta_0}_{t_{ik}}}{T^{2-H_0+u_k}_k/\sqrt{n_k}}\rightarrow 0$, almost surely, as $k\rightarrow\infty$.
Using these results, it follows from (\ref{eq:w1}) that
\begin{align}
	\tilde w_{ik}&\sim 
c\frac{T^{2-H_0+u_k}_k}{\sqrt{n_k}},~\mbox{almost surely, as}~~i\rightarrow n_k;~k\rightarrow\infty.
	\label{eq:w2}
\end{align}
where, if $\tilde H>H_0$, then $u_k=2(H_k-H_0)$ and
\begin{equation}
	c=\beta_0\frac{c_{H_0}}{c^2_{\tilde H}}\frac{4\tilde H^2}{H_0} 
	k^{-1}_{H_0}C_2(H_0,\tilde x); 
	\label{eq:w_3}
\end{equation}
if $\tilde H< H_0$, then $u_k=0$ and
\begin{equation}
	c = -\tilde\beta\frac{c_{H_0}}{c^2_{\tilde H}}\frac{4\tilde H^2}{H_0} 
	 k^{-1}_{\tilde H}C_2(\tilde H,\tilde x);
	\label{eq:w_4}
\end{equation}
and if $\tilde H=H_0$, then $u_k=0$ and
\begin{equation*}
	 c=(\beta_0-\tilde\beta)\frac{4H_0}{c_{H_0}}k^{-1}_{H_0}C_2(H_0,\tilde x).
\end{equation*}
Hence, Lemma \ref{lemma:lemma1} together with (\ref{eq:limit_term3}), (\ref{eq:w0}), (\ref{eq:w1}) and (\ref{eq:w2}), yields
\begin{equation}
	\sum_{i=1}^{n_k}Y_{t_{ik}}\tilde w_{ik}\sim c\frac{T^{2-H_0+u_k}_k}{\sqrt{n_k}}\sum_{i=1}^{n_k}Y_{t_{ik}},~\mbox{almost surely, as}~k\rightarrow\infty.
	\label{eq:sum1}
\end{equation}
Hence, combining (\ref{eq:sum1}) with (\ref{eq:sum0}) yields, almost surely, 
\begin{align}
	&-\sum_{i=1}^{n_k} \log \frac{\upsilon^2_{M^{H_0}_{t_{ik}}}}{\upsilon^2_{M^{H_k}_{t_{ik}}}}
	-\sum_{i=1}^{n_k} Y^2_{t_{ik}}+ \sum_{i=1}^{n_k} Y^2_{t_{ik}}\frac{\upsilon^2_{M^{H_0}_{t_{ik}}}}{\upsilon^2_{M^{H_k}_{t_{ik}}}}
	+2\sum_{i=1}^{n_k}\frac{(\Delta  M^{H_k}_{t_{ik}}-\Delta  M^{H_0}_{t_{ik}})\Delta  M^{H_0}_{t_{ik}}}{\upsilon^2_{M^{H_0}_{t_{ik}}}}
	\frac{\upsilon^2_{M^{H_0}_{t_{ik}}}}{\upsilon^2_{M^{H_k}_{t_{ik}}}}\notag\\
	&\qquad\sim n_kK T^{2(H_k-H_0)}_k-n_k\log K -n_k\log T^{2(H_k-H_0)}_k-n_k+c\frac{T^{2-H_0+u_k}_k}{\sqrt{n_k}}\sum_{i=1}^{n_k}Y_{t_{ik}},~\mbox{as}~k\rightarrow\infty.
	\label{eq:sum2}
\end{align}
To calculate the limit of the right hand side of (\ref{eq:sum2}), let us first consider the case $\tilde H>H_0$. Since $H_k\rightarrow \tilde H$ as $k\rightarrow\infty$,
it follows that for any $\epsilon>0$, there exists $k_0(\epsilon)\geq 1$, such that for $k\geq k_0(\epsilon)$, 
$2(\tilde H-H_0)-\epsilon<2(H_k-H_0)<2(\tilde H-H_0)+\epsilon$. Hence,
for $k\geq k_0(\epsilon)$, $T^{2(\tilde H-H_0)-\epsilon}_k<T^{2(H_k-H_0)}_k<T^{2(\tilde H-H_0)+\epsilon}_k$.
Hence, for $\tilde H>H_0$, $T^{2(H_k-H_0)}_k\rightarrow\infty$, as $k\rightarrow\infty$. 
Now, since $\frac{T^2_k}{n_k}\rightarrow 0$ 
as $k\rightarrow\infty$, it follows that
\begin{equation}
	\frac{T^{2-H_0+2(H_k-H_0)}_k}{n_kT^{2(H_k-H_0)}}=\frac{T^{2-H_0}_k}{n_k}\rightarrow 0,~\mbox{as}~k\rightarrow\infty.
	\label{eq:w3}
\end{equation}
Now note that 
$\frac{1}{\sqrt{n_k}}\sum_{i=1}^{n_k}Y_{t_{ik}}\sim N(0,1)$, for $k\geq 1$. That is, $\frac{1}{\sqrt{n_k}}\sum_{i=1}^{n_k}Y_{t_{ik}}$ remains almost surely finite,
as $k\rightarrow\infty$. Hence, it follows using $T^{2(H_k-H_0)}_k\rightarrow\infty$, as $k\rightarrow\infty$, and (\ref{eq:w3}), that almost surely,
\begin{equation*}
	n_kK T^{2(H_k-H_0)}_k-n_k\log K -n_k\log T^{2(H_k-H_0)}_k-n_k+c\frac{T^{2-H_0+u_k}_k}{\sqrt{n_k}}\sum_{i=1}^{n_k}Y_{t_{ik}}\sim n_kK T^{2(H_k-H_0)}_k
	\rightarrow\infty,~\mbox{as}~k\rightarrow\infty.
\end{equation*}
In other words, the right hand side of (\ref{eq:sum2})
tends to infinity as $k\rightarrow\infty$, when $\tilde H>H_0$.\\
Now consider the case where $\tilde H<H_0$. Then the right hand side of (\ref{eq:sum2}) is equal to 
$\frac{n_kK}{T^{2(H_0-H_k)}_k}-n_k\log K +n_k\log T^{2(H_0-H_k)}_k-n_k+c\frac{T^{2-H_0}_k}{\sqrt{n_k}}\sum_{i=1}^{n_k}Y_{t_{ik}}$.
Now observe that
\begin{align}
&\frac{n_k\log T^{2(H_0-H_k)}_k}{\frac{n_kK}{T^{2(H_0-H_k)}_k}}=\frac{1}{K}T_k^{2(H_0-H_k)}\log T^{2(H_0-H_k)}_k\rightarrow\infty,
	~\mbox{as}~k\rightarrow\infty;\label{eq:w_lim1}\\
&\frac{n_k\log T^{2(H_0-H_k)}_k}{T_k^{2-H_0}}=2(H_0-H_k)\times\frac{n_k}{T^2_k}\times T^{H_0}_k\log T_k\rightarrow\infty,~\mbox{as}~k\rightarrow\infty.\label{eq:w_lim2}
\end{align}
Using (\ref{eq:w_lim1}) and (\ref{eq:w_lim2}) we see that, almost surely
\begin{equation}
	\frac{n_kK}{T^{2(H_0-H_k)}_k}-n_k\log K +n_k\log T^{2(H_0-H_k)}_k-n_k+c\frac{T^{2-H_0}_k}{\sqrt{n_k}}\sum_{i=1}^{n_k}Y_{t_{ik}}
	\sim n_k\log T^{2(H_0-H_k)}_k\rightarrow\infty,~\mbox{as}~k\rightarrow\infty.
	\label{eq:w_lim3}
\end{equation}
In other words, for $\tilde H\neq H_0$, even if $\tilde\beta=\beta_0$, 
\begin{align}
	&-\sum_{i=1}^{n_k} \log \frac{\upsilon^2_{M^{H_0}_{t_{ik}}}}{\upsilon^2_{M^{H_k}_{t_{ik}}}}
	-\sum_{i=1}^{n_k} Y^2_{t_{ik}} + \sum_{i=1}^{n_k} Y^2_{t_{ik}}\frac{\upsilon^2_{M^{H_0}_{t_{ik}}}}{\upsilon^2_{M^{H_k}_{t_{ik}}}}\notag\\
	&\qquad+2\sum_{i=1}^{n_k} \frac{(\Delta  M^{H_k}_{t_{ik}}-\Delta  M^{H_0}_{t_{ik}})\Delta  M^{H_0}_{t_{ik}}}{\upsilon^2_{M^{H_0}_{t_{ik}}}}
	\frac{\upsilon^2_{M^{H_0}_{t_{ik}}}}{\upsilon^2_{M^{H_k}_{t_{ik}}}}\rightarrow\infty,~\mbox{almost surely, as}~k\rightarrow\infty.
\label{eq:likelihood_diff_asymp1}
\end{align}
Since the term $\sum_{i=1}^{n_k} \frac{(\Delta  M^{H_k}_{t_{ik}}-\Delta  M^{H_0}_{t_{ik}})^2}{\upsilon^2_{M^{H_0}_{t_{ik}}}}
	\frac{\upsilon^2_{M^{H_0}_{t_{ik}}}}{\upsilon^2_{M^{H_k}_{t_{ik}}}}$ of (\ref{eq:likelihood_diff1}) is non-negative, it follows that
$\ell_k(\btheta_k)-\ell_k(\btheta_0) \rightarrow\infty,~\mbox{almost surely, as}~k\rightarrow\infty$.\\
Now consider the case where $H_k=H_0$, for $k=1,2,\ldots$, but $\tilde\beta\neq\beta_0$. 
In this case, 
\begin{equation}
-\sum_{i=1}^{n_k} \log \frac{\upsilon^2_{M^{H_0}_{t_{ik}}}}{\upsilon^2_{M^{H_k}_{t_{ik}}}}
	-\sum_{i=1}^{n_k} Y^2_{t_{ik}} + \sum_{i=1}^{n_k} Y^2_{t_{ik}}\frac{\upsilon^2_{M^{H_0}_{t_{ik}}}}{\upsilon^2_{M^{H_k}_{t_{ik}}}}=0, 
	\label{eq:zero1}
\end{equation}
and it follows from (\ref{eq:w3}) and (\ref{eq:sum1}) that
\begin{align}
&2\sum_{i=1}^{n_k}\frac{(\Delta  M^{H_k}_{t_{ik}}-\Delta  M^{H_0}_{t_{ik}})\Delta  M^{H_0}_{t_{ik}}}{\upsilon^2_{M^{H_0}_{t_{ik}}}}
	\frac{\upsilon^2_{M^{H_0}_{t_{ik}}}}{\upsilon^2_{M^{H_k}_{t_{ik}}}}\notag\\
	&\sim (\beta_0-\tilde\beta)\times\frac{4H_0}{c_{H_0}}\times k^{-1}_{H_0}C_2(H_0,\tilde x)\times\frac{T^{2-H_0}}{\sqrt{n_k}}\sum_{i=1}^{n_k}Y_{t_{ik}},\notag\\
	&\qquad~\mbox{almost surely, as}~k\rightarrow\infty.
	\label{eq:zero2}
\end{align}
Also, in a similar way of obtaining (\ref{eq:sum1}), we see that
\begin{align}
	&\frac{(\Delta  M^{H_k}_{t_{ik}}-\Delta  M^{H_0}_{t_{ik}})^2}{\upsilon^2_{M^{H_0}_{t_{ik}}}}
	\frac{\upsilon^2_{M^{H_0}_{t_{ik}}}}{\upsilon^2_{M^{H_k}_{t_{ik}}}}\notag\\
	&\sim \frac{4\tilde H^2 }{c^2_{\tilde H}\frac{T^{2-2H_k}_k}{n_k}}	
	\times
	\left[\beta_0k^{-1}_{H_0}\frac{T^{2-2H_0}_k}{n_k}\left\{(2-2H_0)C_1(H_0,\tilde x)+T_kC_2(H_0,\tilde x)\right\}\right.\notag\\
	&\qquad\left.-\tilde \beta k^{-1}_{\tilde H}\frac{T^{2-2H_k}_k}{n_k}\left\{(2-2\tilde H)C_1(\tilde H,\tilde x)+T_kC_2(\tilde H,\tilde x)\right\}\right]^2,
	~\mbox{almost surely, as}~i\rightarrow n_k;~k\rightarrow\infty.\notag\\
	&= \frac{4\tilde H^2 }{n_kc^2_{\tilde H}T^{2-2H_k}_k}	
	\times
	\left[\beta_0k^{-1}_{H_0}T^{2-2H_0}_k \left\{(2-2H_0)C_1(H_0,\tilde x)+T_kC_2(H_0,\tilde x)\right\}\right.\notag\\
	&\qquad\left.-\tilde \beta k^{-1}_{\tilde H}T^{2-2H_k}_k\left\{(2-2\tilde H)C_1(\tilde H,\tilde x)+T_kC_2(\tilde H,\tilde x)\right\}\right]^2.
	\label{eq:zero3}
\end{align}
Hence, if, for $k\geq 1$, $H_k=\tilde H=H_0$ and $\tilde\beta\neq\beta_0$, then it holds due to (\ref{eq:zero3}) that
\begin{align}
\sum_{i=1}^{n_k} \frac{(\Delta  M^{H_k}_{t_{ik}}-\Delta  M^{H_0}_{t_{ik}})^2}{\upsilon^2_{M^{H_0}_{t_{ik}}}}
	&\frac{\upsilon^2_{M^{H_0}_{t_{ik}}}}{\upsilon^2_{M^{H_k}_{t_{ik}}}}\sim
	(\tilde\beta-\beta_0)^2\times\frac{4H^2_0}{c^2_{H_0}}\times\left\{k^{-1}_{H_0}C_2(H_0,\tilde x)\right\}^2\times T^{4-2H_0}_k,\notag\\
	&\qquad~\mbox{almost surely, as}~k\rightarrow\infty.
	\label{eq:zero4}
\end{align}
Combining (\ref{eq:zero1}), (\ref{eq:zero2}) and (\ref{eq:zero4}) yields in the case $H_k=H_0$, for $k=1,2,\ldots$, and $\tilde\beta\neq\beta_0$, yields
\begin{align}
	&\ell_k(\btheta_k)-\ell_k(\btheta_0)\sim
	(\beta_0-\tilde\beta)\times\frac{4H_0}{c_{H_0}}\times k^{-1}_{H_0}C_2(H_0,\tilde x)\times\frac{T^{2-H_0}_k}{\sqrt{n_k}}\sum_{i=1}^{n_k}Y_{t_{ik}}\notag\\
	&\qquad+(\tilde\beta-\beta_0)^2\times\frac{4H^2_0}{c^2_{H_0}}\times\left\{k^{-1}_{H_0}C_2(H_0,\tilde x)\right\}^2\times T^{4-2H_0}_k
	\label{eq:zero_full}\\
	&\qquad\sim(\tilde\beta-\beta_0)^2\times\frac{4H^2_0}{c^2_{H_0}}\times\left\{k^{-1}_{H_0}C_2(H_0,\tilde x)\right\}^2\times T^{4-2H_0}_k\notag\\
        &\rightarrow\infty,~\mbox{almost surely, as}~k\rightarrow\infty.
	\label{eq:zero5}
\end{align}
Now note that $\ell_{T,n}(\btheta)$ is continuous in $\btheta$, $T$ and $n$. Using continuity, in conjunction with 
(\ref{eq:ass_lim3}) and (\ref{eq:ass_lim4}), it follows that for any other subsequence $H_k$ converging to $H_0$, with $\tilde\beta\neq\beta_0$, it must hold that
$\ell_k(\btheta_k)-\ell_k(\btheta_0)\rightarrow\infty,~\mbox{almost surely, as}~k\rightarrow\infty$.\\
We thus conclude that for any $\delta>0$, $\underset{\btheta\in\bTheta:\|\btheta-\btheta_0\|\geq \delta}{\inf}~\left\{\ell_{T,n}(\btheta)-\ell_{T,n}(\btheta_0)\right\}
\rightarrow\infty,~\mbox{almost surely, as}~n\rightarrow\infty,~T\rightarrow\infty$, such that $n/T^2\rightarrow\infty$.
In other words, $\hat\btheta_{T,n}=\left(\hat\beta_{T,n},\hat H_{T,n}\right)\stackrel{a.s.}{\longrightarrow}\btheta_0=\left(\beta_0, H_0\right)$,
as $n\rightarrow\infty,~T\rightarrow\infty$, such that $n/T^2\rightarrow\infty$.
\end{proof}

\section{Proof of Corollary \ref{theorem:consistency2}}
\begin{proof}
	Under either of (E1), (E2) and (E3), we have
$$\frac{d}{dt}\left.\int_0^1 u^{\frac{1}{2}- H}(1-u)^{\frac{1}{2}- H}\frac{B(tu,\tilde x_{tu})}{C(tu,\tilde x_{tu})}du\right|_{t=t_{i}}\rightarrow 0$$
almost surely, as $i\rightarrow n_k$, $k\rightarrow\infty$,
	and $C_1(H,\tilde x)\equiv C_1(H)$.
In these cases it follows from (\ref{eq:w1}) that
\begin{align}
	\tilde w_{ik}&\sim 
c\frac{T^{1-H_0+u_k}_k}{\sqrt{n_k}},~\mbox{almost surely, as}~~i\rightarrow n_k;~k\rightarrow\infty.\notag
\end{align}
where, if $\tilde H>H_0$, then $u_k=2(H_k-H_0)$ and
\begin{equation*}
	c=\beta_0\frac{c_{H_0}}{c^2_{\tilde H}}\frac{4\tilde H^2}{H_0} 
	k^{-1}_{H_0}(2-2H_0)C_1(H_0); 
\end{equation*}
if $\tilde H<H_0$, then $u_k=0$ and
\begin{equation*}
	c = -\tilde\beta\frac{c_{H_0}}{c^2_{\tilde H}}\frac{4\tilde H^2}{H_0} 
	 k^{-1}_{\tilde H}(2-2\tilde H)C_1(\tilde H);
\end{equation*}
and if $\tilde H=H_0$, then $u_k=0$ and
\begin{equation*}
	c = (\beta_0-\tilde\beta)\frac{4H_0}{c_{H_0}} 
	 k^{-1}_{H_0}(2-2H_0)C_1(H_0);
\end{equation*}
These imply that
\begin{equation*}
2\sum_{i=1}^{n_k}\frac{(\Delta  M^{H_k}_{t_{ik}}-\Delta  M^{H_0}_{t_{ik}})\Delta  M^{H_0}_{t_{ik}}}{\upsilon^2_{M^{H_0}_{t_{ik}}}}
	\frac{\upsilon^2_{M^{H_0}_{t_{ik}}}}{\upsilon^2_{M^{H_k}_{t_{ik}}}}\notag\\
=\sum_{i=1}^{n_k}Y_{t_{ik}}\tilde w_{ik}\sim c\frac{T^{1-H_0+u_k}_k}{\sqrt{n_k}}\sum_{i=1}^{n_k}Y_{t_{ik}},~\mbox{almost surely, as}~k\rightarrow\infty.
\end{equation*}
and
\begin{align}
	&\sum_{i=1}^{n_k} \frac{(\Delta  M^{H_k}_{t_{ik}}-\Delta  M^{H_0}_{t_{ik}})^2}{\upsilon^2_{M^{H_0}_{t_{ik}}}}
	\frac{\upsilon^2_{M^{H_0}_{t_{ik}}}}{\upsilon^2_{M^{H_k}_{t_{ik}}}}\sim
	(\tilde\beta-\beta_0)^2\times\frac{4H^2_0}{c^2_{H_0}}\times\left\{k^{-1}_{H_0}(2-2H_0)C_1(H_0)\right\}^2\times T^{2-2H_0},\notag\\
	&\qquad~\mbox{almost surely, as}~k\rightarrow\infty.
\end{align}
	The rest of the proof follows in the same way as that of Theorem \ref{theorem:consistency} by 	setting $T_k/n_k\rightarrow 0$ instead of $T^2_k/n_k\rightarrow 0$, as $k\rightarrow\infty$. When $H_k=H_0$ and $\tilde\beta\neq \beta_0$, we must replace $T^{4-2H_0}_k$ with $T^{2-2H_0}_k$ in
	(\ref{eq:zero4}).
\end{proof}

\section{Proof of Theorem \ref{theorem:asymp_normal1}}
\begin{proof}
	To obtain the asymptotic distribution of $\beta$, let us minimize $\ell_{T,n}(\btheta)$ by differentiation with respect to $\beta$, fixing $H=H_0$.
	In this case, since $\ell_{T,n}(\beta,H_0)$ is a quadratic form in $\beta$,  it is easy to verify that differentiating $\ell_{T,n}(\beta,H_0)$ with respect to $\beta$ 
and solving for $\beta$ by setting the differential equal to zero, and finally applying asymptotic equivalence to the solution, is the same as first obtaining
the asymptotic equivalence of $\ell_{T,n}(\beta,H_0)-\ell_{T,n}(\beta_0,H_0)$, and then differentiating the asymptotic equivalent with respect to $\beta$ and finding the solution
by setting the differential to zero.

Indeed, in the same way as (\ref{eq:zero_full}),  
\begin{align}
	&\ell_{T,n}(\beta,H_0)-\ell_{T,n}(\beta_0,H_0)\sim
	(\beta_0-\beta)\times\frac{4H_0}{c_{H_0}}\times k^{-1}_{H_0}C_2(H_0,\tilde x)\times\frac{T^{2-H_0}}{\sqrt{n}}\sum_{i=1}^{n}Y_{t_{i}}\notag\\
	&\qquad+(\beta-\beta_0)^2\times\frac{4H^2_0}{c^2_{H_0}}\times\left\{k^{-1}_{H_0}C_2(H_0,\tilde x)\right\}^2\times T^{4-2H_0},\notag\\
	&\qquad~\mbox{almost surely, as}
	~T\rightarrow\infty;~n\rightarrow\infty;~\frac{n}{T^2}\rightarrow\infty.
	\label{eq:normal1}
\end{align}
Differentiating the asymptotic equivalent on the right hand side of (\ref{eq:normal1}) with respect to $\beta$ and setting equal to zero we obtain
the solution $\hat\beta_{T,n}$ as
\begin{equation}
	T^{2-H_0}(\hat\beta_{T,n}-\beta_0)=a_{H_0}Z_{T,n},
	\label{eq:normal2}
\end{equation}
where
\begin{align}
	& a_{H_0}=\frac{k_{H_0}c_{H_0}}{2H_0C_2(H_0,\tilde x)};\label{eq:a0}\\
	& Z_{T,n}=\frac{1}{\sqrt{n}}\sum_{i=1}^{n}Y_{t_{i}}.\label{eq:Z}
\end{align}
Note that for any $T>0$ and $n\geq 1$, $Z_{T,n}$ is equivalent in distribution to $Z$, where $Z\sim N(0,1)$. 
From (\ref{eq:normal2}), (\ref{eq:a0}) and (\ref{eq:Z}), it follows that
\begin{equation}
	T^{2-H_0}(\hat\beta_{T,n}-\beta_0)\rightarrow a_{H_0}Z\sim N\left(0,a^2_{H_0}\right),
	~\mbox{almost surely, as}~T\rightarrow\infty;~n\rightarrow\infty;~\frac{n}{T^2}\rightarrow\infty.
	\label{eq:normal3}
\end{equation}
Hence, continuity of logarithm ensures that
\begin{equation}
	(2-H_0)\log T-\log |a_{H_0}|+\log\left|\hat\beta_{T,n}-\beta_0\right|\rightarrow \log|Z|,~\mbox{almost surely, as}~T\rightarrow\infty;~n\rightarrow\infty;~\frac{n}{T^2}\rightarrow\infty.
	\label{eq:normal4}
\end{equation}
	Now note that for any fixed $\beta\in \mathbb R$, the same method of the proof of Theorem \ref{theorem:consistency} shows that 
	$\hat H_{T,n}\stackrel{a.s.}{\longrightarrow} H_0$. Thus, due to continuity of $a_H$ with respect to $H$,
$$(2-H_0)\log T-\log \left|a_{H_0}\right|+\log\left|\hat\beta_{T,n}-\beta_0\right|
\sim (2-\hat H_{T,n})\log T-\log \left|a_{\hat H_{T,n}}\right|+\log\left|\hat\beta_{T,n}-\beta_0\right|,$$
almost surely, as $T\rightarrow\infty;~n\rightarrow\infty;~\frac{n}{T^2}\rightarrow\infty$.
Hence,
\begin{equation}
	(2-\hat H_{T,n})\log T-\log \left|a_{\hat H_{T,n}}\right|+\log\left|\hat\beta_{T,n}-\beta_0\right|\rightarrow \log|Z|,
	~\mbox{almost surely, as}~T\rightarrow\infty;~n\rightarrow\infty;~\frac{n}{T^2}\rightarrow\infty.
	\label{eq:normal5}
\end{equation}
Using continuity of exponential we obtain from (\ref{eq:normal5}) that
\begin{equation}
	\frac{T^{2-\hat H_{T,n}}}{\left|a_{\hat H_{T,n}}\right|}\left|\hat\beta_{T,n}-\beta_0\right|\rightarrow \left|Z\right|,	
	~\mbox{almost surely, as}~T\rightarrow\infty;~n\rightarrow\infty;~\frac{n}{T^2}\rightarrow\infty.
	\label{eq:normal6}
\end{equation}
Now let $\tilde Z$ be a random variable that takes $\left| Z\right|$ and $-\left| Z\right|$ with probabilities $1/2$ and $1/2$. Then 
$\tilde Z\sim N(0,1)$. Since due to (\ref{eq:normal6}), $\left| Z\right|$ and $-\left| Z\right|$ are the asymptotic forms of 
$\frac{T^{2-\hat H_{T,n}}}{\left|a_{\hat H_{T,n}}\right|}\left|\hat\beta_{T,n}-\beta_0\right|$ and
$-\frac{T^{2-\hat H_{T,n}}}{\left|a_{\hat H_{T,n}}\right|}\left|\hat\beta_{T,n}-\beta_0\right|$, respectively, it follows that
\begin{equation*}
	\frac{T^{2-\hat H_{T,n}}}{\left|a_{\hat H_{T,n}}\right|}\left(\hat\beta_{T,n}-\beta_0\right)\rightarrow \tilde Z\sim N(0,1),	
	~\mbox{almost surely, as}~T\rightarrow\infty;~n\rightarrow\infty;~\frac{n}{T^2}\rightarrow\infty.
\end{equation*}
This completes the proof.
\end{proof}

\section{Proof of Theorem \ref{theorem:asymp_normal2}}
\begin{proof}
In this case we have  
\begin{align}
	&\ell_{T,n}(\beta,H_0)-\ell_{T,n}(\beta_0,H_0)\sim
	(\beta_0-\beta)\times\frac{4H_0}{c_{H_0}}\times k^{-1}_{H_0}(2-2H_0)C_1(H_0)\times\frac{T^{1-H_0}}{\sqrt{n}}\sum_{i=1}^{n}Y_{t_{i}}\notag\\
	&\qquad+(\beta-\beta_0)^2\times\frac{4H^2_0}{c^2_{H_0}}\times\left\{k^{-1}_{H_0}(2-2H_0)C_1(H_0)\right\}^2\times T^{2-2H_0},\notag\\
	&\qquad~\mbox{almost surely, as}
	~T\rightarrow\infty;~n\rightarrow\infty;~\frac{n}{T}\rightarrow\infty.\notag
\end{align}
The rest of the proof follows in the same way as that of Theorem \ref{theorem:asymp_normal1}.
\end{proof}

\section{Proof of Theorem \ref{theorem:incons_noncompact1}}
\begin{proof}
	Differentiating $\ell_{T,n}(\beta,H^*_{T,n})-\ell_{T,n}(\beta_0,H_0)$ with respect to $\beta$ and solving for $\beta$ after setting the differential equal to zero yields
	the solution $\beta=\tilde\beta_{T,n}$, where
	\begin{align}
		&\tilde\beta_{T,n}\sim \beta_0\left(\frac{d_{H_0}T^{2-2H_0}+b_{H_0}T^{3-2H_0}}{d_{\tilde H}T^{2-2H^*_{T,n}}+b_{\tilde H}T^{3-2H^*_{T,n}}}\right)
		\label{eq:solve0}\\
		&\qquad+\frac{C_3Z_{T,n}T^{2-H_0}\left[(2-2\tilde H)T^{1-2H^*_{T,n}}C_1(\tilde H,\tilde x)+T^{2-2H^*_{T,n}}C_2(\tilde H,\tilde x)\right]}{\left(d_{\tilde H}T^{2-2H^*_{T,n}}+b_{\tilde H}T^{3-2H^*_{T,n}}\right)^2},\label{eq:solve1}\\
		&\qquad~\mbox{almost surely, as}~T\rightarrow\infty~\mbox{and}
		~n\rightarrow\infty, 
		\notag
	\end{align}
	where $Z_{T,n}=n^{-1/2}\sum_{i=1}^nY_{t_i}\sim N(0,1)$ for any $T>0$ and $n\geq 1$, $C_3=\frac{c_{H_0}}{H_0}\times k^{-1}_{\tilde H}$
	and for any $H\in[\eta,1-\eta]$, $d_H=2k^{-1}_H(1-H)C_1(H,\tilde x)$; $b_H=k^{-1}_HC_2(H,\tilde x)$.
		The second term (\ref{eq:solve1}) is asymptotically equivalent to $T^{-H_0-2(1-H^*_{T,n})}$, which 
		tends to zero almost surely, as $T\rightarrow\infty$, irrespective of $\tilde H>H_0$ or $\tilde H<H_0$.
		However, the first term (\ref{eq:solve0}) goes to infinity if $\tilde H>H_0$, and to zero if $\tilde H<H_0$. 
		In other words, $\tilde\beta_{T,n}$ 
	tends to infinity if $\tilde H>H_0$ and to zero if $\tilde H<H_0$, almost surely, as $T\rightarrow\infty$.
	Hence, $\tilde\beta_{T,n}$ given above 
		is an inconsistent estimator of $\beta_0$. 
		Note that the conclusion does not depend upon compactness or non-compactness of the parameter space $\mathfrak B$, since in neither case
		$\tilde\beta_{T,n}$ can converge to $\beta_0$, for all $\beta_0\in\mathfrak B$.
	\end{proof}

\section{Proof of Theorem \ref{theorem:incons_noncompact2}}
	\begin{proof}
		Under (E1), (E2) or (E3), for sufficiently large $T$, differentiating $\ell_{T,n}(\beta,H^*_{T,n})-\ell_{T,n}(\beta_0,H_0)$ with respect to $\beta$ 
		and solving for $\beta$ after setting the differential equal to zero yields
	the solution $\beta=\tilde\beta_{T,n}$, where
	\begin{align}
		&\tilde\beta_{T,n}\sim \beta_0\left(\frac{d_{H_0}}{d_{\tilde H}}\right)T^{2(H^*_{T,n}-H_0)}
		+C_3Z_{T,n}\left(\frac{2-2\tilde H}{d^2_{\tilde H}}\right)T^{1-H_0-2(1-H^*_{T,n})}\notag\\
		&\qquad\sim \beta_0\left(\frac{d_{H_0}}{d_{\tilde H}}\right)T^{2(H^*_{T,n}-H_0)},\label{eq:solve2}\\
		&\qquad~\mbox{almost surely, as}~T\rightarrow\infty;~
		n\rightarrow\infty, 
		\notag
	\end{align}
		where the definitions of $Z_{T,n}$ and $C_3$ remain the same as in the proof of Theorem \ref{theorem:incons_noncompact1}, and
		$d_H=2k^{-1}_H(1-H)C_1(H)$. Note that (\ref{eq:solve2}) follows from the previous step since 
		$$\frac{T^{2(H^*_{T,n}-H_0)}}{T^{1-H_0-2(1-H^*_{T,n})}}=T^{1-H_0}\rightarrow\infty,~\mbox{as}~~T\rightarrow\infty;~n\rightarrow\infty,$$
		since $H_0<1$.

		Hence if $\tilde H>H_0$ and $\tilde H<H_0$, then using (\ref{eq:solve2}) it follows that
		$\tilde\beta_{T,n}\rightarrow\infty$ and $\tilde\beta_{T,n}\rightarrow 0$, respectively, almost surely, as $T\rightarrow\infty$ and $n\rightarrow\infty$.

		Again, note that the conclusion does not depend upon compactness or non-compactness of the parameter space $\mathfrak B$, since in neither case
		$\tilde\beta_{T,n}$ can converge to $\beta_0$, for all $\beta_0\in\mathfrak B$.
	\end{proof}

\section{Proof of Theorem \ref{theorem:bayesian_consistency}}
\begin{proof}
	As argued in Remark \ref{remark:classical}, the $MLE$ $\hat\btheta_{T,n}$ is strongly consistent for $\btheta_0$ even when the parameter space of $\beta$ is
	$\mathbb R$. Hence, the proof follows by a simple application of Doob's theorem (see, for
	example, Theorem 7.78 of \ctn{Schervish95}) which is easily seen to hold valid in our case as $T\rightarrow\infty$, $n\rightarrow\infty$, 
	such that $n/T^2\rightarrow\infty$.
\end{proof}

\section{Proof of Theorem \ref{theorem:bayesian_consistency2}}
\begin{proof}
Again, the proof is a simple application of Doob's theorem.
\end{proof}

\section{Proof of Theorem \ref{theorem:bayesian_consistency_KL}}
\begin{proof}
Note that 
\begin{align}
	&\pi\left(\beta\in N_0|H=H_0,\tilde x_{t_i};~i=1,\ldots,n\right)\notag\\
	&=\frac{\int_{N_0}\exp\left\{-\frac{1}{2}\left(\ell_{T,n}(\beta,H_0)-\ell_{T,n}(\beta_0,H_0)\right)\right\}d\pi(\beta)}
	{\int_{N_0}\exp\left\{-\frac{1}{2}\left(\ell_{T,n}(\beta,H_0)-\ell_{T,n}(\beta_0,H_0)\right)\right\}d\pi(\beta)
	+\int_{N^c_0}\exp\left\{-\frac{1}{2}\left(\ell_{T,n}(\beta,H_0)-\ell_{T,n}(\beta_0,H_0)\right)\right\}d\pi(\beta)}
	\label{eq:posterior1}
\end{align}
If we can show that 
\begin{equation}
	\frac{\int_{N_0}\exp\left\{-\frac{1}{2}\left(\ell_{T,n}(\beta,H_0)-\ell_{T,n}(\beta_0,H_0)\right)\right\}d\pi(\beta)}
	{\int_{N^c_0}\exp\left\{-\frac{1}{2}\left(\ell_{T,n}(\beta,H_0)-\ell_{T,n}(\beta_0,H_0)\right)\right\}d\pi(\beta)}
	\rightarrow\infty,~\mbox{almost surely, as}~T\rightarrow\infty;~n\rightarrow\infty;~\frac{n}{T^2}\rightarrow\infty,
	\label{eq:posterior_asymp1}
\end{equation}
then it would prove that 
	\begin{equation}	
	\pi\left(\beta\in N_0|H=H_0,\tilde x_{t_i};~i=1,\ldots,n\right)\stackrel{a.s.}{\longrightarrow}1.
		\label{eq:posterior_asymp2}
	\end{equation}
	Now, 
	\begin{equation}
	\pi\left(\beta\in N_0|H=H_0,\tilde x_{t_i};~i=1,\ldots,n\right)=E\left(I_{N_0}(\beta)|H=H_0,\tilde x_{t_i};~i=1,\ldots,n\right),
		\label{eq:posterior_asymp3}
	\end{equation}
	and the expression on the right hand side is a martingale. By the martingale convergence theorem (see, for example, Theorem B.118 of \ctn{Schervish95}),
	$E\left(I_{N_0}(\beta)|H=H_0,\tilde x_{t_i};~i=1,\ldots,n\right)\stackrel{a.s.}{\longrightarrow}
	 E\left(I_{N_0}(\beta)|H=H_0,\tilde x_t;~t\geq 0\right)$, as $T\rightarrow\infty$, $n\rightarrow\infty$ and $n/T^2\rightarrow\infty$.
	 This, combined with (\ref{eq:posterior_asymp2}) and (\ref{eq:posterior_asymp3}), shows that
	 \begin{equation}
		 E\left(I_{N_0}(\beta)|H=H_0,\tilde x_t;~t\geq 0\right)=1.
		 \label{eq:posterior_asymp4}
	 \end{equation}
	Since $$\pi\left(\beta\in N_0|H=\hat H_{T,n},\tilde x_{t_i};~i=1,\ldots,n\right)=E\left(I_{N_0}(\beta)|H=\hat H_{T,n},\tilde x_{t_i};~i=1,\ldots,n\right)$$
	is also a martingale, by the martingale convergence theorem, almost sure consistency of $\hat H_{T,n}$
	and (\ref{eq:posterior_asymp4}),
\begin{align}
&\pi\left(\beta\in N_0|H=\hat H_{T,n},\tilde x_{t_i};~i=1,\ldots,n\right)\stackrel{a.s.}{\longrightarrow}
E\left(I_{N_0}(\beta)|H=H_0,\tilde x_t;~t\geq 0\right)=1,\notag\\
	&\qquad~\mbox{as}~T\rightarrow\infty,~n\rightarrow\infty,~\frac{n}{T^2}\rightarrow\infty.
\end{align}
Thus, we only require to prove that (\ref{eq:posterior_asymp1}) holds.
	In this regard, first let $\bar C_{\epsilon}=\left\{\beta:(\beta-\beta_0)^2g(H_0)\leq\epsilon\right\}$, the closure of $C_{\epsilon}$ and note that
\begin{align}
	&\int_{N_0}\exp\left\{-\frac{1}{2}\left(\ell_{T,n}(\beta,H_0)-\ell_{T,n}(\beta_0,H_0)\right)\right\}d\pi(\beta)\notag\\
	&\qquad\geq
	\int_{\bar C_{\epsilon}}\exp\left\{-\frac{T^{4-2H_0}}{2}
	\frac{\underset{\beta\in \bar C_{\epsilon}}{\sup}~\left(\ell_{T,n}(\beta,H_0)-\ell_{T,n}(\beta_0,H_0)\right)}{T^{4-2H_0}}\right\}d\pi(\beta).
	\label{eq:p1}
\end{align}
	Now, from (\ref{eq:normal1}), we see that, point-wise with respect to $\beta\in \bar C_{\epsilon}$, 
\begin{equation}
	\frac{\left(\ell_{T,n}(\beta,H_0)-\ell_{T,n}(\beta_0,H_0)\right)}{T^{4-2H_0}}\stackrel{a.s.}{\longrightarrow}(\beta-\beta_0)^2g(H_0),
	~\mbox{as}~T\rightarrow\infty,~n\rightarrow\infty,~\frac{n}{T^2}\rightarrow\infty.
	\label{eq:p2}
\end{equation}
	Hence, for $\beta_1,\beta_2\in \bar C_{\epsilon}$, it follows from (\ref{eq:p2}) and continuity of the absolute function, that
\begin{equation}
	\frac{\left|\ell_{T,n}(\beta_1,H_0)-\ell_{T,n}(\beta_2,H_0)\right|}{T^{4-2H_0}}\stackrel{a.s.}{\longrightarrow}
	g(H_0)\left|\beta_1-\beta_2\right|\left|\beta_1+\beta_2-2\beta_0\right|,
	~\mbox{as}~T\rightarrow\infty,~n\rightarrow\infty,~\frac{n}{T^2}\rightarrow\infty.
	\label{eq:p3}
\end{equation}
Since $\beta_1,\beta_2\in \bar C_{\epsilon}$, they are bounded, and so, $\left|\beta_1+\beta_2-2\beta_0\right|\leq M$, for some $M>0$.
This, and (\ref{eq:p3}) implies that $\frac{\left(\ell_{T,n}(\beta,H_0)-\ell_{T,n}(\beta_0,H_0)\right)}{T^{4-2H_0}}$ 
is stochastically equicontinuous, and hence (see, for example, \ctn{Newey91}),
\begin{equation}
	\underset{\beta\in \bar C_{\epsilon}}{\sup}~\left|\frac{\left(\ell_{T,n}(\beta,H_0)-\ell_{T,n}(\beta_0,H_0)\right)}{T^{4-2H_0}}-(\beta-\beta_0)^2g(H_0)\right|
	\stackrel{a.s.}{\longrightarrow}0.
	\label{eq:p4}
\end{equation}
Now, due to (\ref{eq:p4}), for $\epsilon>0$, there exists $n_0(\epsilon)\geq 1$ and $T_0(\epsilon)>0$ such that $n\geq n_0(\epsilon)$ and $T\geq T_0(\epsilon)$,
\begin{align}
	&\underset{\beta\in \bar C_{\epsilon}}{\sup}~\frac{\left(\ell_{T,n}(\beta,H_0)-\ell_{T,n}(\beta_0,H_0)\right)}{T^{4-2H_0}}\notag\\
	&\qquad	\leq\underset{\beta\in \bar C_{\epsilon}}{\sup}~\left|\frac{\left(\ell_{T,n}(\beta,H_0)-\ell_{T,n}(\beta_0,H_0)\right)}{T^{4-2H_0}}-(\beta-\beta_0)^2g(H_0)\right|
	+\underset{\beta\in \bar C_{\epsilon}}{\sup}~(\beta-\beta_0)^2g(H_0)\notag\\
	&\qquad \leq\epsilon+(\beta^*-\beta_0)^2g(H_0)\notag\\
	&\qquad\leq 2\epsilon,
	\label{eq:p5}
\end{align}
where $\beta^*\in\bar C_{\epsilon}$ is such that $\underset{\beta\in \bar C_{\epsilon}}{\sup}~(\beta-\beta_0)^2g(H_0)=(\beta^*-\beta_0)^2g(H_0)\leq\epsilon$, thanks to
compactness of $\bar C_{\epsilon}$. It follows from (\ref{eq:p5}) and (\ref{eq:p1}), that
for $n\geq n_0(\epsilon)$ and $T\geq T_0(\epsilon)$,
\begin{align}
	&\int_{N_0}\exp\left\{-\frac{1}{2}\left(\ell_{T,n}(\beta,H_0)-\ell_{T,n}(\beta_0,H_0)\right)\right\}d\pi(\beta)\notag\\
	&\qquad\geq
	\int_{\bar C_{\epsilon}}\exp\left\{-\epsilon T^{4-2H_0}\right\}d\pi(\beta)
	=\pi(\bar C_{\epsilon})\exp\left\{-\epsilon T^{4-2H_0}\right\}.
	\label{eq:p6}
\end{align}
Now observe that
\begin{align}
	&\int_{N^c_0}\exp\left\{-\frac{1}{2}\left(\ell_{T,n}(\beta,H_0)-\ell_{T,n}(\beta_0,H_0)\right)\right\}d\pi(\beta)\notag\\
	&\qquad\leq
	\int_{C^c_{\epsilon}}\exp\left\{-\frac{T^{4-2H_0}}{2}
	\frac{\underset{\beta\in C^c_{\epsilon}}{\inf}~\left(\ell_{T,n}(\beta,H_0)-\ell_{T,n}(\beta_0,H_0)\right)}{T^{4-2H_0}}\right\}d\pi(\beta).
	\label{eq:p7}
\end{align}
Also, there exists $\tilde\beta\neq\beta_0$ such that 
\begin{align}
	&\frac{\underset{\beta\in C^c_{\epsilon}}{\inf}~\left(\ell_{T,n}(\beta,H_0)-\ell_{T,n}(\beta_0,H_0)\right)}{T^{4-2H_0}}
	\geq \frac{\left(\ell_{T,n}(\tilde\beta,H_0)-\ell_{T,n}(\beta_0,H_0)\right)}{T^{4-2H_0}}
	\rightarrow (\tilde\beta-\beta_0)^2g(H_0),\notag\\
	&\qquad~\mbox{almost surely, as}~T\rightarrow\infty;~n\rightarrow\infty;~\frac{n}{T^2}\rightarrow\infty.
	\label{eq:p8}
\end{align}
It follows due to (\ref{eq:p8}) that for any $\epsilon>0$, there exists
$n_1(\epsilon)\geq 1$ and $T_1(\epsilon)>0$ such that $n\geq n_1(\epsilon)$ and $T\geq T_1(\epsilon)$, 
\begin{equation}
	\frac{\underset{\beta\in C^c_{\epsilon}}{\inf}~\left(\ell_{T,n}(\beta,H_0)-\ell_{T,n}(\beta_0,H_0)\right)}{T^{4-2H_0}}
	\geq (\tilde\beta-\beta_0)^2g(H_0)-\epsilon.
	\label{eq:p9}
\end{equation}
Using (\ref{eq:p9}) in (\ref{eq:p6}) we obtain for $n\geq n_1(\epsilon)$ and $T\geq T_1(\epsilon)$, 
\begin{align}
	&\int_{N^c_0}\exp\left\{-\frac{1}{2}\left(\ell_{T,n}(\beta,H_0)-\ell_{T,n}(\beta_0,H_0)\right)\right\}d\pi(\beta)\notag\\
	&\qquad\leq
	\int_{C^c_{\epsilon}}\exp\left\{-\frac{T^{4-2H_0}\left\{(\tilde\beta-\beta_0)^2g(H_0)-\epsilon\right\}}{2}\right\}d\pi(\beta)\notag\\
	&\qquad=\pi(C^c_{\epsilon})\exp\left\{-\frac{T^{4-2H_0}\left\{(\tilde\beta-\beta_0)^2g(H_0)-\epsilon\right\}}{2}\right\}.
	\label{eq:p10}
\end{align}
From (\ref{eq:p6}) and (\ref{eq:p10}) we obtain for $n\geq\max\{n_0(\epsilon),n_1(\epsilon)\}$ and for $T\geq\max\{T_0(\epsilon),T_1(\epsilon)\}$, that
\begin{align}
	&\frac{\int_{N_0}\exp\left\{-\frac{1}{2}\left(\ell_{T,n}(\beta,H_0)-\ell_{T,n}(\beta_0,H_0)\right)\right\}d\pi(\beta)}
	{\int_{N^c_0}\exp\left\{-\frac{1}{2}\left(\ell_{T,n}(\beta,H_0)-\ell_{T,n}(\beta_0,H_0)\right)\right\}d\pi(\beta)}
	\geq\frac{\pi(\bar C_{\epsilon})}{\pi(C^c_{\epsilon})}\exp\left\{\frac{T^{4-2H_0}\left\{(\tilde\beta-\beta_0)^2g(H_0)-2\epsilon\right\}}{2}\right\}\notag\\
	&\qquad\qquad\rightarrow\infty,~\mbox{almost surely, as}~T\rightarrow\infty;~n\rightarrow\infty;~\frac{n}{T^2}\rightarrow\infty.\notag
\end{align}
This proves the theorem.
\end{proof}

\section{Proof of Theorem \ref{theorem:bayesian_consistency2_KL}}
\label{sec:bayesian_consistency2_KL}
\begin{proof}
The proof is similar to that of Theorem \ref{theorem:bayesian_consistency_KL}.
\end{proof}

\section{Posterior asymptotic normal setup and the main result in \ctn{Schervish95}}
\label{sec:schervish}
Let $\tilde\ell_n(\btheta)= \log L_n(\btheta)$ denote the log-likelihood associated with $n$ observations $x_1,\ldots,x_n$
having density $f(x_1,\ldots,x_n|\btheta)$ with parameter $\btheta\in\bTheta$, where $\bTheta\subseteq\mathbb R^p$ is the parameter space, for $p\geq 1$.
Now denoting 
by $\hat\btheta_n$ the $MLE$ associated with
$n$ observations, 
let
\begin{equation}
\Sigma^{-1}_n=\left\{\begin{array}{cc}
-\tilde\ell''_n(\hat\btheta_n) & \mbox{if the inverse and}\ \ \hat\btheta_n\ \ \mbox{exist}\\
\mathbb I_p & \mbox{if not},
\end{array}\right.
\label{eq:information1}
\end{equation}
where for any $\bt\in\mathbb R^p$,
\begin{equation}
\tilde\ell''_n(\bt)=\left(\left(\frac{\partial^2}{\partial\theta_i\partial\theta_j}\tilde\ell_n(\btheta)\bigg\vert_{\btheta=\bt}\right)\right),
\label{eq:information2}
\end{equation}
and $\mathbb I_p$ is the identity matrix of order $p$.
Thus, $\Sigma^{-1}_n$ is the observed Fisher's information matrix.
\subsubsection{Regularity conditions}
\label{subsubsec:regularity}
\begin{itemize}
\item[(1)] The true value $\btheta_0$ is a point interior to $\bTheta$.
\item[(2)] The prior distribution of $\btheta$ has a density with respect to Lebesgue measure
that is positive and continuous at $\btheta_0$.
\item[(3)] There exists a neighborhood $\mathcal N_0\subseteq\bTheta$ of $\btheta_0$ on which
$\tilde\ell_n(\btheta)= \log f(x_1,\ldots,x_n|\btheta)$ is twice continuously differentiable with 
respect to all co-ordinates of $\btheta$, 
almost surely with respect to the true, data-generating distribution $P_{\btheta_0}$.
\item[(4)] The largest eigenvalue of $\Sigma_n$ goes to zero in probability.
\item[(5)] For $\delta>0$, define $\mathcal N_0(\delta)$ to be the open ball of radius $\delta$ around $\btheta_0$.
Let $\rho_n$ be the smallest eigenvalue of $\Sigma_n$. If $\mathcal N_0(\delta)\subseteq\bTheta$, there exists
$K(\delta)>0$ such that
\begin{equation}
\underset{n\rightarrow\infty}{\lim}~P_{\btheta_0}\left(\underset{\btheta\in\bTheta\backslash\mathcal N_0(\delta)}{\sup}~
\rho_n\left[\tilde\ell_n(\btheta)-\tilde\ell_n(\btheta_0)\right]<-K(\delta)\right)=1.
\label{eq:extra1}
\end{equation}
\item[(6)] For each $\epsilon>0$, there exists $\delta(\epsilon)>0$ such that
\begin{equation}
\underset{n\rightarrow\infty}{\lim}~P_{\btheta_0}\left(\underset{\btheta\in\mathcal N_0(\delta(\epsilon)),\|\gamma\|=1}{\sup}~
\left\vert 1+\gamma^T\Sigma^{\frac{1}{2}}_n\tilde\ell''_n(\btheta)\Sigma^{\frac{1}{2}}_n\gamma\right\vert<\epsilon\right)=1.
\label{eq:extra2}
\end{equation}
\end{itemize}
\begin{theorem}[Theorem 7.89 of \ctn{Schervish95}]
\label{theorem:posterior_normality1}
Assume the above six regularity conditions.
Letting $\Psi_n=\Sigma^{-1/2}_n\left(\btheta-\hat\btheta_n\right)$ and denoting by $\pi_n$ the posterior distribution of $\Psi_n$ given $x_1,\ldots,x_n$, 
for each compact subset $B$ of $\mathbb R^p$ and each $\epsilon>0$, the following holds:
\begin{equation}
\lim_{n\rightarrow\infty}P_{\btheta_0}
\left(\sup_{\bpsi\in B}\left\vert\pi_n(\bpsi)-\phi(\bpsi)\right\vert>\epsilon\right)=0,
\label{eq:Bayesian_normality}
\end{equation}
where $\phi(\cdot)$ denotes the density of the $p$-dimensional standard normal distribution.
\end{theorem}

\section{Proof of Theorem \ref{theorem:post_normal1}}
\label{sec:post_normal1_supp}
\begin{proof}
Note that $\tilde\ell_{T,n}(\btheta)=-\ell_{T,n}(\btheta)/2$, where we have replaced $\tilde\ell_{n}(\btheta)$ with $\tilde\ell_{T,n}(\btheta)$
for obvious reasons. For our purpose, we shall deal with $\tilde\ell_{T,n}(\beta,H_0)$. For the relevant Taylor's series expansion considered in the
proof of Theorem 7.89 provided in \ctn{Schervish95}, we consider the expansion of $\tilde\ell_{T,n}(\beta,H_0)$ around the $MLE$ 
$\hat\btheta_{T,n}=(\hat\beta_{T,n},\hat H_{T,n})$. 

	To obtain $\tilde\ell''_{T,n}(\beta,H_0)$, we consider the above Taylor's series expansion of $\tilde\ell_{T,n}(\beta,H_0)$ 
	around $\hat\btheta_{T,n}$ consisting of at least the first four terms, that is, at least till the term consisting of three derivatives. 
	It is then easy to verify by direct calculation that for all $\beta\in\mathbb R$ and also for $\beta=\hat\beta_{T,n}$,
\begin{equation}
	\tilde\ell''_{T,n}(\beta,H_0)\sim 
	-a^{-2}_{\hat H_{T,n}}T^{4-2\hat H_{T,n}},
	~\mbox{almost surely, as}~T\rightarrow\infty;~n\rightarrow\infty;~\frac{n}{T^2}\rightarrow\infty,
	\label{eq:norm1}
\end{equation}
where, for any $H$, $a_{H}$ is given by (\ref{eq:norm2}).
	Thus,
	\begin{equation}
		\Sigma_{T,n}=-\left(\tilde\ell''_{T,n}(\hat\beta_{T,n},H_0)\right)^{-1} \sim a^2_{\hat H_{T,n}}T^{-4+2\hat H_{T,n}}, 
	~\mbox{almost surely, as}~T\rightarrow\infty;~n\rightarrow\infty;~\frac{n}{T^2}\rightarrow\infty.
		\label{eq:norm3}
	\end{equation}

Now observe that the first three regularity conditions of Section \ref{subsubsec:regularity} are trivially satisfied, and because of (\ref{eq:norm3}),
condition (4) is also satisfied. 

Now note that condition (5), associated with $\rho_{T,n}$, is required to control equation (7.97) of the proof of \ctn{Schervish95}.
In this regard, first observe that the aforementioned equation remains valid even if
we consider the Taylor's series expansion of $\tilde\ell_{T,n}(\beta,H_0)$ in equation (7.97) around $(\hat\beta_{T,n},H_0)$ instead of around $(\hat\beta_n,\hat H_{T,n})$.  
To obtain $\rho_{T,n}$ corresponding to this case, we expand $\tilde\ell_{T,n}(\beta,H_0)$ $(\hat\beta_{T,n},H_0)$ till at least the first four terms, as before.
This yields $\rho_{T,n}\sim a^2_{H_0}T^{-4+2H_0}$,
almost surely, as $T\rightarrow\infty;~n\rightarrow\infty;~\frac{n}{T^2}\rightarrow\infty$.

Now, since $\tilde \ell_{T,n}(\beta,H_0)$ is quadratic in $\beta$, the quantity
$\underset{|\beta-\beta_0|>\delta}{\sup}~\rho_{T,n}\left[\tilde \ell_{T,n}(\beta,H_0)-\tilde\ell_{T,n}(\beta_0,H_0)\right]$ 
is $\rho_{T,n}\left[\tilde \ell_{T,n}(\beta^*,H_0)-\tilde\ell_{T,n}(\beta_0,H_0)\right]$, where $\beta^*=\beta_0+\delta$ or $\beta_0-\delta$.
Hence, 
	\begin{align}
		&\underset{|\beta-\beta_0|>\delta}{\sup}~\rho_{T,n}\left[\tilde \ell_{T,n}(\beta,H_0)-\tilde\ell_{T,n}(\beta_0,H_0)\right]<-K(\delta)\notag\\
		&\Leftrightarrow \rho_{T,n}\left[\tilde \ell_{T,n}(\beta^*,H_0)-\tilde\ell_{T,n}(\beta_0,H_0)\right]<-K(\delta).
		\label{eq:norm4}
	\end{align}
The left hand side of (\ref{eq:norm4}) converges almost surely to $-(\beta^*-\beta_0)^2/2=-\delta^2/2$, 
	as $T\rightarrow\infty;~n\rightarrow\infty;~\frac{n}{T^2}\rightarrow\infty$. Hence, for $K(\delta)<\delta^2/2$, condition (5) holds.

To verify condition (6), note that by (\ref{eq:norm1}), it follows that for any $\epsilon>0$, 
$\left\vert 1+\gamma^T\Sigma^{\frac{1}{2}}_n\tilde\ell''_n(\btheta)\Sigma^{\frac{1}{2}}_n\gamma\right\vert<\epsilon$ as 
$T\rightarrow\infty;~n\rightarrow\infty;~\frac{n}{T^2}\rightarrow\infty$, for $\|\gamma\|=1$. Here $\gamma$ is a scalar
and so $\gamma^T=\gamma$ and $\|\gamma\|=|\gamma|$. Hence, condition (6) holds trivially.

Thus, the proof of Theorem \ref{theorem:post_normal1} is complete.
\end{proof}

\section{Proof of Theorem \ref{theorem:post_normal2}}
\label{sec:post_normal2}
\begin{proof}
The proof is similar to that of Theorem \ref{theorem:post_normal1}.
\end{proof}

\normalsize
\bibliographystyle{natbib}
\bibliography{irmcmc}

\end{document}